\documentclass[10pt]{article}
\usepackage{amsmath,amssymb,amsthm,xspace,amscd,latexsym,amssymb,epic,stmaryrd}
\usepackage[mathscr]{eucal}
\usepackage[all]{xy}
\usepackage{enumerate}
\date{}
\numberwithin{equation}{section}

\theoremstyle{plain}
\newtheorem{theorem}{Theorem}
\newtheorem{definition}[theorem]{Definition}
\newtheorem{remark}[theorem]{Remark}
\newtheorem{lemma}[theorem]{Lemma}
\newtheorem{proposition}[theorem]{Proposition}
\newtheorem{corollary}[theorem]{Corollary}

\DeclareMathOperator{\Ob}{Ob}

\DeclareMathOperator{\Hom}{Hom}
\DeclareMathOperator{\HOM}{Hom}

\DeclareMathOperator{\Ext}{Ext}

\newcommand{\Mod}{\mathrm{Mod}}
\renewcommand{\mod}{\mathrm{mod}}
\newcommand{\gMod}{\mathrm{gMod}}
\newcommand{\gmod}{\mathrm{gmod}}
\newcommand{\fgmod}{\mathrm{gfmod}}
\newcommand{\fmod}{\mathrm{fmod}}
\newcommand{\fdgmod}{\mathrm{gfdmod}}
\newcommand{\fdmod}{\mathrm{fdmod}}
\newcommand{\mg}{\mathfrak{g}}
\newcommand{\mb}{\mathfrak{b}}
\newcommand{\mh}{\mathfrak{h}}
\newcommand{\cU}{\mathcal{U}}
\newcommand{\cO} {\mathcal{O}}
\newcommand{\cL} {\mathcal{L}}

\newcommand{\Z}{\mathbb{Z}}
\newcommand{\mZ}{\mathbb{Z}}
\newcommand{\mC}{\mathbb{C}}
\newcommand{\defis}{\text{-}}
\newcommand{\la}{\lambda}
\newcommand{\op}{\operatorname}
\newcommand{\END}{\op{End}}
\newcommand{\p}{\mathfrak{p}}
\renewcommand{\mp}{\mathfrak{p}}
\newcommand{\PK}{\mathbb{P}_{\text{\tiny $k$}}^{\bullet}}
\newcommand{\PKK}{\mathbb{P}_{\text{\tiny $k$}}}
\newcommand{\rhom}{\mathcal{R}\operatorname{Hom}}

\title{Quadratic duals, Koszul dual functors, \\ and applications\footnote{AMS-classification 16S37, 18E30,16G20,17B67}}
\author{Volodymyr Mazorchuk\footnote{partially supported by the
Swedish Research Council}, Serge Ovsienko\footnote{partially
supported by the Royal Swedish Academy of Sciences and The Swedish
Foundation for International Cooperation in Research and Higher
Education (STINT)} \\ and Catharina Stroppel\footnote{Supported by
The Engineering and Physical Sciences Research Council (EPSRC)}}

\begin{document}

\maketitle

\begin{abstract}
The paper studies quadratic and Koszul duality for modules over
positively graded categories. Typical examples are modules over a
path algebra, which is graded by the path length, of a not
necessarily finite quiver with relations. We present a very
general definition of quadratic and Koszul duality functors backed
up by explicit examples. This generalizes \cite{BGS} in two
substantial ways: We work in the setup of graded categories, i.e.
we allow infinitely many idempotents and also define a ``Koszul''
duality {\it functor} for not necessarily Koszul categories. As an
illustration of the techniques we reprove the Koszul duality
(\cite{Steen}) of translation and Zuckerman functors for the
classical category $\cO$ in a quite elementary and explicit way.
From this we deduce a conjecture of \cite{BFK}. As applications we
propose a definition of a ``Koszul'' dual category for integral
blocks of Harish-Chandra bimodules and for blocks outside the
critical hyperplanes for the Kac-Moody category $\cO$.
\end{abstract}

\tableofcontents

\section{Introduction}\label{s1}

This paper deals with (categories of) modules over positively graded
categories, defines quadratic duality and studies Koszul duality.
The first motivation behind this is to get a generalized Koszul or
quadratic duality which also works for module categories over not
necessarily finite dimensional, not necessarily unital algebras.
In our opinion, the language of {\it modules over (positively) graded
categories} is very well adapted to this task. Our second motivation
is to provide a  definition of quadratic duality {\it functors} for any
quadratic algebra. These functors give rise to the usual Koszul
duality functors for Koszul algebras. Remind yourself that a positively
graded algebra is Koszul if all simple modules have a linear
projective resolution. Despite this definition and the vast amount
of literature on Koszul algebras and Koszul duality (see, for example,
\cite{BGS,GK,GRS,Ke2} and references therein), and, in particular, its
relation to linear resolutions (see, for example,
\cite{GMVRSZ,HI,MVZ} and references therein), it seems that
(apart from \cite{MVS}) there are no attempts to study Koszul {\em duality}
by working seriously with the category of linear complexes of
projective modules. The intention of the paper is to provide the following:
\begin{itemize}
\item a very general definition of a quadratic dual category in terms
of the category of linear complexes of projectives, and a detailed
study of the latter (Section~\ref{section2}, Section~\ref{se3}),
\item a general setup of quadratic/Koszul duality for positively
graded categories instead of positively graded algebras
(Section~\ref{s2}) with an abstract definition of a duality
functor (Section~\ref{se3}) and a Koszul duality theorem for
Koszul categories. Using the word {\em duality} here might be too
optimistic, in particular, since the functors are not even equivalences in
general (see Theorem~\ref{koszul}). However, later on we will see many
``duality-like'' effects in our situation, which, from our point of view, justify this usage.
\item an illustration of our techniques in form of an alternative
proof of the statement that translation functors and Zuckerman functors
for the classical Bernstein-Gelfand-Gelfand category $\cO$ are Koszul
dual to each other (Theorem~\ref{Steen}). This fact is well-known,
was conjectured in \cite{BGS} and proved in \cite{Steen} using
dg-algebras. Our approach seems to be more elementary and more
explicit. As a consequence we prove that twisting/completion
and shuffling/coshuffling functors are Koszul dual. Although these
functors are well-studied their Koszul duality was a surprise even
for specialists. It clarifies the connection of the two categorifications
of \cite{BFK} and establishes a direct connection between the main
result of \cite{StDuke} and a result in \cite{Josh} (Section~\ref{s5}).
\item an elementary description of the Koszul complex as a
complex of $\mathbf{C}\defis\mathbf{C}^!$-bimodules
(see Section~\ref{s6}). This provides a connection to the quite recent article
\cite{Fl}.
\item A complimentary approach to Lef{\`e}vre's and Keller's
generalization of the Koszul duality from \cite{Ke3}.
The combination of these two approaches provides a sort
of quadratic homological duality.
\end{itemize}

A (positively) graded category $\mathbf{C}$ is a small category
with (positively) graded morphism spaces. To our knowledge, the study
of modules over categories was initiated by Bredon (\cite{Bredon}) and
tom Dieck (\cite{Dieck}) in the obstruction theory for finite groups
and appears now in different variations, see for example \cite{Ga,Mi,Ke2}.
In their setup, the categories of (right) modules over a category
$\mathbf{C}$ play an important role, where by definition a
{\it module} is a covariant functor from $\mathbf{C}$ to the category of
finite dimensional vector spaces.  From our perspective, modules over
a positively graded category $\mathbf{C}$ should be thought of
as representations of a (not necessarily unital) positively graded
algebra. One could consider $\mathbf{C}$-modules as representations
of the (not necessarily finite) quiver with relations associated with $\mathbf{C}$.
The vertices of the quiver with relations correspond to the objects in $\mathbf{C}$
and the path algebra $A$ is just the direct sum of all morphism spaces
$\bigoplus_{\la,\mu}\mathbf{C}(\la,\mu)$, where the sum runs over all
pairs of objects $\la$, $\mu$ from $\mathbf{C}$. Then, a
$\mathbf{C}$-module is a functor which associates to each object
in the category (i.e. to each vertex in the quiver) a finite dimensional
vector space and to each morphism (hence to each arrow in the quiver) a
linear map between the corresponding vector spaces. The functoriality
guarantees that each $\mathbf{C}$-module is exactly
what a representation of the associated quiver with relations should be. Note that,
in case $\mathbf{C}$ has only finitely many objects, the objects
from $\mathbf{C}$ are in bijection with a maximal set of pairwise
orthogonal primitive idempotents $e_\la$ of $A$ corresponding to the
identity elements in $\mathbf{C}(\la,\la)$. From the definitions of
graded categories (Definition~\ref{defpgc}) we get a correspondence as follows:
\begin{eqnarray}
\label{introequiv}
\begin{array}[thb]{ccccccc}
\left\{
\parbox{4cm}{positively graded categories with finitely many
objects and finite dimensional morphism spaces}
\right\}&
&\leftrightarrow&
&\left\{
\parbox{4cm}{finite dimensional positively
graded algebras}
\right\}
\end{array}
\end{eqnarray}
by mapping a category $\mathbf{C}$ to the graded algebra
$\oplus_{\la,\mu\in\op{Ob}(\mathbf{C})}\mathbf{C}(\la,\mu)$ of
morphisms. In the opposite direction, an algebra $A$ is mapped to the
category, whose  objects are a chosen system of pairwise orthogonal
primitive idempotents and morphisms are the morphisms between the
associated indecomposable projective modules.

Under this correspondence equivalent categories correspond to
Morita equivalent algebras and isomorphic categories to isomorphic
algebras. It is also easy to see that the notions of modules correspond.
If we remove the additional finiteness assumptions there is no such nice
correspondence, since there is no natural choice for a maximal set of pairwise
orthogonal idempotents. Therefore, one should think of graded categories as
the correct language to speak about algebras with a fixed set of pairwise
orthogonal primitive idempotents (see also \cite{BoGa}).
We want to illustrate the results of the paper in the
following two examples:

\subsection*{A (well-behaved) illustrating example}

We consider the $\mC$-algebra $A$ which is the path algebra of the quiver
\begin{eqnarray}
\label{Adef}
\parbox[b]{5cm}{\hspace*{1cm}{
\xymatrix{1 \ar@/^/[r]^{f}&2 \ar@/^/[l]^{g}}}}
\end{eqnarray}
modulo the relation $g\circ f=0$ (i.e. the loop starting at vertex
$1$ is zero). Putting the arrows in degree one defines a
non-negative $\mZ$-grading on $A$. We denote this graded algebra
by $\mathsf{A}$. Note that $\mathsf{A}_0$ is semi-simple. The
algebra $\mathsf{A}$ is quadratic and its quadratic dual is the
algebra $\mathsf{A}^!$ given as the path algebra of the same
quiver, but with the relation $f\circ g=0$. This algebra is again
graded by putting the arrows in degree one.  We get decompositions
$\mathsf{A}=\mathsf{P}(1)\oplus \mathsf{P}(2)$ and
$\mathsf{A}^!=\mathsf{P}^!(1)\oplus \mathsf{P}^!(2)$ into
indecomposable  (graded) projective $\mathsf{A}$-modules
corresponding to the vertices of  the quiver. Note that the
indecomposable projective graded $\mathsf{A}$-modules are all of
the form $\mathsf{P}(i)\langle j\rangle$, where $i\in\{1,2\}$,
$j\in\mZ$, and $\langle j\rangle$ shifts the grading of the module
down by $j$.

\subsubsection*{Linear complexes of projective $\mathsf{A}$-modules and the equivalence $\epsilon$}

To describe the category of finite dimensional, graded
$\mathsf{A}^!$-modules  we use a result of \cite{MVS} which says
that this category is  equivalent to $\mathscr{LC}(\mathsf{P})$,
the so-called category of linear complexes of projective modules,
i.e. complexes of projective $\mathsf{A}$-modules, where in
position $j$ we have a direct sum of projective modules of the
form $\mathsf{P}(i)\langle j\rangle$ for $i\in\{1,2\}$, $j\in\mZ$, each
occurring with finite multiplicity.
The category $\mathscr{LC}(\mathsf{P})$ is abelian with the usual
kernels and cokernels (Proposition~\ref{pcatlcomp}), the simple
objects are exactly the indecomposable objects
$\mathsf{P}(i)\langle j\rangle$, for $i\in\{1,2\}$, $j\in\mZ$,
considered as linear complexes with support concentrated in position
$j$. Let $\mathsf{S}(i)$ be the simple top of the graded
$\mathsf{A}$-module $\mathsf{P}(i)$ and let
$\mathsf{I}(i)$ be the injective hull of $\mathsf{S}(i)$.
Note that the simple objects in $\mathsf{A}\defis\fgmod$ are
exactly the $\mathsf{S}(i)\langle j\rangle$, where  $i\in\{1,2\}$,
$j\in\mZ$. Similarly we define $\mathsf{S}^!(i)$,
$\mathsf{P}^!(i)$, $\mathsf{I}^!(i)$ for the algebra $\mathsf{A}^!$.
Then the equivalence $\epsilon^{-1}: \mathsf{A}^!\defis\gMod\cong\mathscr{LC}(\mathsf{P})$
(Theorem~\ref{tqdual}) gives a correspondence as follows:
\small
\begin{eqnarray*}
{
  \begin{array}[th]{c|cccccccccc}
\mathsf{A}^!\defis\gMod&&&&\mathscr{LC}(\mathsf{P})\\
\hline\\
\mathsf{S}^!(1) &&&0&\rightarrow &\mathsf{P}(1)&\rightarrow& 0&& \\
\mathsf{S}^!(2) &&&0&\rightarrow &\mathsf{P}(2)&\rightarrow &0\\
\mathsf{I}^!(2)&0&\rightarrow&\mathsf{P}(1)\langle -1\rangle&\rightarrow
   &\mathsf{P}(2)&\rightarrow&0&&\\
\mathsf{I}^!(1)&\mathsf{P}(1)\langle -2\rangle&\rightarrow&
\mathsf{P}(2)\langle -1\rangle&\rightarrow
&\mathsf{P}(1)&\rightarrow&0&& \\
\mathsf{P}^!(2)&&&0&\rightarrow&\mathsf{P}(2)&\rightarrow
   &\mathsf{P}(1)\langle 1\rangle&&\\
\mathsf{P}^!(1)&&&0&\rightarrow &\mathsf{P}(1)&\rightarrow&
\mathsf{P}(2)\langle 1\rangle&\rightarrow
&\mathsf{P}(1)\langle 2\rangle&\nonumber
 \end{array}
}
\end{eqnarray*}
\normalsize
where the maps in the complexes are the obvious ones and the not shown parts of the complexes are just trivially zero. Note that
the indecomposable projective module $\mathsf{P}(i)\langle -j\rangle[j]$
occurs exactly $[\mathsf{M}:\mathsf{S}(i)^!\langle j\rangle]$ times in
the complex associated to $\mathsf{M}$. The maps in the complexes
are naturally obtained from the action of $\mathsf{A}^!$ on
$\mathsf{M}$. This equivalence
$\mathsf{A}^!\defis\fgmod\cong\mathscr{LC}(\mathsf{P})$ will be
explained in Theorem~\ref{tqdual} in the general setup of locally
finite dimensional modules over a quadratic graded category $\mathbf{C}$.
In Proposition~\ref{pltcproj} we will describe the indecomposable
injective objects in $\mathscr{LC}(\mathsf{P})$. It turns out
that the {\it injective} hull of the simple module $\mathsf{S}^!(i)$,
is nothing else than the maximal linear part of a minimal
{\it projective} resolution of $\mathsf{S}(i)$. Since the algebra
$\mathsf{A}$ from our example above is in fact Koszul, the minimal
projective resolution of $\mathsf{S}(i)$ is automatically linear.
In Proposition~\ref{pltcproj} we also describe how to get the
indecomposable projective objects: we take a minimal injective resolution,
(for $\mathsf{S}(2)$ we get $\mathsf{I}(2)\rightarrow \mathsf{I}(1)
\langle 1\rangle$), then we apply the inverse of  the Nakayama functor
(we get $\mathsf{P}(2)\rightarrow \mathsf{P}(1)\langle 1\rangle$),
finally we take the maximal linear part of the result (since the
resolution in our example is already linear, we are done).

\subsubsection*{The Koszul self-duality}

The algebra $\mathsf{A}$ from our example is very special, since it is
Koszul self-dual, i.e. $\mathsf{A}$ is isomorphic to its quadratic dual $\mathsf{A}^!$ (\cite[Theorem~18]{Sperv} for $\mg=\mathfrak{sl}_2$). An isomorphism is of course given by identifying $\mathsf{P}(1)$
with $\mathsf{P}^!(2)$ and $\mathsf{P}(2)$ with $\mathsf{P}(1)^!$.
In general,  the quadratic dual $\mathsf{A}^!$ could be very different
from $\mathsf{A}$. In Proposition~\ref{corext} we give a homological
characterization of the quadratic dual of a positively graded category. In the
example it gets reduced to the fact that
$\mathsf{A}^!\cong\op{Ext}^\bullet\big(\mathsf{S}(1)\oplus
\mathsf{S}(2), \mathsf{S}(1)\oplus \mathsf{S}(2)\big)$, which is
the usual Koszul dual. Note that if $\mathsf{A}$ is any, not
necessarily a finite dimensional, positively graded algebra (in the
sense of Definition~\ref{defpga}) of finite global dimension, then
its quadratic dual $\mathsf{A}^!$ is finite-dimensional
(Corollary~\ref{bred}, see also \cite[Section~10.4]{Ke2}).

\subsubsection*{The Koszul dual functors}

In Section~\ref{section3} we define a generalization of
(the pair of adjoint) Koszul dual functors (Theorem~\ref{tkfunctors}).
Using the category of linear complexes of projectives, it is easy to
describe the (inverse) Koszul functor $\mathrm{K}_{A}'$
(Lemma~\ref{funnyd2}): Given a graded $\mathsf{A}^!$-module, the
equivalence $\epsilon^{-1}$ maps this module to a linear complex of
projective graded $\mathsf{A}$-modules. This can be considered as an
object in the bounded derived category $\mathcal{D}^b(\mathsf{A})$
of graded $A$-modules. Hence we have a functor
$\mathrm{K}'_{\mathsf{A}(0)}:\mathsf{A}^!\defis\fgmod\rightarrow
\mathcal{D}^b(\mathsf{A})$ which can easily be extended to a
functor defined on $\mathcal{D}^b(\mathsf{A}^!)$. For example, the
simple graded $\mathsf{A}^!$-modules $\mathsf{S}(i)^!$ are
mapped to the complexes with $\mathsf{P}(i)$ concentrated in degree zero,
hence $\mathsf{S}(i)^!$ is mapped to $\mathsf{P}(i)$
(see Theorem~\ref{tkfunctors}\eqref{tkfunctors.3} for a
general result). Since the algebra $\mathsf{A}$ in the example is
Koszul, so is $\mathsf{A}^!$ and the Koszul functors are inverse to
each other (see \cite{BGS}). This statement will be generalized
in Theorem~\ref{koszul}.

\subsubsection*{Connection to representation theory of Lie algebras}

In Section~\ref{s4} we consider a special case of \cite{BGS} and
\cite{ErikKoszul}, namely the Koszul duality functor for blocks
of the Bernstein-Gelfand-Gelfand category $\cO$ associated to any
semisimple Lie algebra. For the principal block of $\cO$ corresponding
to the Lie algebra $\mathfrak{sl}_2$ we get exactly the category of
finitely generated $\mathsf{A}$-modules as described in the examples above. We
have the  Zuckerman functor $\mathrm{Z}$ which maps an $\mathsf{A}$-module
to its maximal quotient containing only simple composition factors of
the form $\mathsf{S}(1)$. If $e$ is the primitive idempotent corresponding
to $\mathsf{P}(2)$ then $\mathrm{Z}$ is a functor from the category of
finitely generated $\mathsf{A}$-modules to the category of finitely
generated $\mathsf{A}/\mathsf{A}e\mathsf{A}$-modules. The functor is
right exact and has the  obvious right adjoint (exact) functor
$\mathrm{i}$. If we take the  left derived functor of the composition we
get  $\cL(\mathrm{i} \mathrm{Z})\mathsf{S}(2)=\mathsf{S}(1)\langle -1\rangle[1]$ and
$\cL(\mathrm{i}  \mathrm{Z}) \mathsf{S}(1)\cong \mathsf{S}(1)\oplus
\mathsf{S}(1)\langle -2\rangle[2]$ (see the  fourth line of the table above). From the
results described above, the Koszul dual functor has to map injective
modules to injective modules. It turns out that this is exactly the
well-known so-called translation functor through the wall.
Section~\ref{s4} provides an alternative proof of the statement that
derived Zuckerman functors and translation functors are Koszul dual.
This was conjectured in \cite{BGS} and proved in the setup of
dg-algebras in \cite{Steen}.  Our proof avoids the use of
dg-algebras, but illustrates again the power of the equivalence
$\epsilon$. Theorem~\ref{Twist} finally shows that the left derived
functor of $\mathsf{A}e\mathsf{A}\otimes_\mathsf{A}\bullet$,
shifted by $\langle 1\rangle$, is Koszul dual to Irving's shuffling
functor (\cite{Irvingshuffle}).

\subsection*{Another (less well-behaved) illustrating example}

Consider the following quiver:
\begin{eqnarray}
\label{quiv1}
\parbox[b]{5cm}{{\xymatrix{1&2 \ar@{->}[l]&3 \ar@{->}[l]&4\ar@{->}[l]&\ar@{->}[l]\cdots}}}
\end{eqnarray}
This defines a positively graded category $\mathbf{C}$ where the objects
are the positive integers and the morphisms $\mathbf{C}(m,n)$ are just
the linear span of the paths from $m$ to $n$. The length of the path
defines a positive grading on $\mathbf{C}$.

\subsubsection*{The problem with projective covers}

The category $\mathbf{C}{\defis}\fgmod$ contains the simple
modules $\mathsf{L}(n)\langle k\rangle$ (concentrated in degree
$-k$, $k\in\mZ$) for any object $n$, and their projective
covers $\mathsf{P}(n)\langle k\rangle$ which has $n$ composition
factors, namely $\mathsf{L}(j)$ occurs in degree $n-j-k$ for any
$1\leq j\leq n$. The injective hull $\mathsf{I}(j)\langle k\rangle$ of
$\mathsf{L}(j)\langle k\rangle$ does not have a finite composition
series. The composition factors of $\mathsf{I}(n)\langle k\rangle$ are the
$\mathsf{L}(j)$ for $j\geq n$, each appearing once, namely in
degree $n-j-k$. Note that this category does not have enough
projectives, since, for example, the indecomposable injective
modules do not have projective covers. This makes life much
more complicated, but it turns out that for any positively graded
category we have at least projective covers and injective hulls for
any module of finite length, in particular for simple modules (see
Lemma~\ref{abcgmold}) and enough projectives in a certain truncated
category (Lemma~\ref{abcgm}). Similar problems can be found e.g. in
\cite{AR}.

\subsubsection*{The quadratic dual via linear complexes of projectives}

Let us look at the category $\mathscr{LC}(\mathsf{P})$, which describes
the quadratic dual. The indecomposable injective objects are the
linear complexes of the form
\begin{eqnarray*}
\cdots\rightarrow0\rightarrow\mathsf{P}(1)\rightarrow
0\rightarrow\cdots&{   \text or    }&
\cdots\rightarrow0\rightarrow\mathsf{P}(i-1)\langle -1\rangle\rightarrow\mathsf{P}(i)\rightarrow 0\rightarrow\cdots,
\end{eqnarray*}
for $i\geq 2$ and their $\langle k\rangle[-k]$-shifts for any $k\in\mZ$, since
they are just the maximal linear parts of the minimal {\it projective}
resolutions of the simple modules (Proposition~\ref{pltcproj}).  The
indecomposable projective objects are the
linear complexes of the form
\begin{eqnarray*}
\cdots\rightarrow0\rightarrow\mathsf{P}(i)\rightarrow\mathsf{P}(i+1)\langle 1\rangle\rightarrow 0\rightarrow\cdots,
\end{eqnarray*}
for any $i\geq 1$ and together with all their $\langle k\rangle[-k]$-shifts.
From the equivalence $\epsilon$ (Theorem~\ref{tqdual}) we get that
the quadratic dual $\mathbf{C}^!$ is the positively graded category
given by the following quiver:
\begin{eqnarray*}
\parbox[b]{5cm}{{\xymatrix{1 \ar@{->}[r]&2 \ar@{->}[r]&3 \ar@{->}[r]&4\ar@{->}[r]&\cdots}}}
\end{eqnarray*}
with the relation that the composition of two consecutive arrows is
always zero. The indecomposable projective module $\mathsf{P}^!(i)$ in
$\mathbf{C}^!{\defis}\fgmod$ has therefore the composition factors
$\mathsf{L}^!(j)$ for $j=i,i+1$ appearing in degree $0$ and $1$
respectively.  The indecomposable injective module $\mathsf{I}^!(i)$
in $\mathbf{C}^!{\defis}\fgmod$ is simple for $i=1$ and has the
composition factors $\mathsf{L}^!(j)$ for $j=i,i-1$ appearing in degree
$0$ and $-1$ respectively.

\subsubsection*{The quadratic dual of the quadratic dual}

Let us consider the category $\mathscr{LC}(\mathsf{P}^!)$. The
indecomposable injectives objects are the linear complexes of the form
\begin{eqnarray}
\cdots\rightarrow\mathsf{P}(i+2)\langle-2\rangle\rightarrow
\mathsf{P}(i+1)\langle-1\rangle\rightarrow\mathsf{P}(i)
\rightarrow 0\rightarrow\cdots
\end{eqnarray}
for $i\geq1$, and their $\langle k\rangle[-k]$-shifts for any $k\in\mZ$.
Since the projective resolutions of the simple modules are linear, these
are nothing else than the projective resolutions of the simple
$\mathbf{C}^!$-modules, and the category $\mathbf{C}^!$ is Koszul.
The indecomposable injective objects are the linear complexes of the form
\begin{eqnarray}
\cdots\rightarrow0\rightarrow\mathsf{P}(i)\rightarrow\mathsf{P}(i-1)\langle
1\rangle\rightarrow \cdots\rightarrow \mathsf{P}(1)\langle i-1\rangle\rightarrow0\rightarrow\cdots,
\end{eqnarray}
for any $i\geq 1$ together with their $\langle k\rangle[-k]$ shifts. From Theorem~\ref{tqdual} we get that the quadratic dual $\mathbf{C}^!$ is
the positively graded category given by the quiver \eqref{quiv1}.
Note that $\mathbf{C}$ and its quadratic dual are both Koszul. The
(inverse) Koszul duality functor $\mathrm{K}'_\mathbf{C}$ is again
nothing else than extending $\epsilon$ to a functor defined on
the corresponding derived category mapping a complex of
locally finite-dimensional graded $\mathbf{C}^!$-modules to a complex of
linear complexes of projectives. Taking the total complex we get a
complex of locally finite-dimensional graded $\mathbf{C}$-modules. This
description of the (inverse) Koszul duality functor can be found in
Proposition~\ref{funnyd2}.

\subsection*{A (classical) family with the same quadratic duals}
Consider the algebra $B(\infty)=\mC[x]$ or $B(n)=\mC[x]/(x^n)$ for any integer
$n\geq 3$. Putting $x$ in degree $1$ we get a graded algebra $\mathsf{B}(n)$
for $n\geq 3$ or $n=\infty$. The maximal linear part of a minimal projective
resolution of the trivial $\mathsf{B}(n)$-module is just the complex
$\mathsf{B}(n)\langle
-1\rangle\stackrel{x\cdot}{\rightarrow}\mathsf{B}(n)$. By
Proposition~\ref{pltcproj} and Theorem~\ref{tqdual} we get a description of the (only) indecomposable
injective $\mathsf{B}(n)^!$-module. In particular, $\mathsf{B}(n)^!\cong
\mC[x]/(x^2)$, independent of
$n$. Proposition~\ref{pltcproj} and Theorem~\ref{tqdual} also imply
$(\mathsf{B}(n)^!)^!\cong(\mC[x]/(x^2))^!\cong \mC[x]$ which is the classical
example of Koszul duality from \cite{BGGKoszul} for $n=2$.

\subsection*{A (too badly behaved) illustrating example}

Consider the path algebra $A$ of the following quiver:
\begin{equation}\label{quiv101}
\xymatrix{
1 \ar[d]& 2\ar[dl] & \dots & n \ar[dlll] &\dots\\
0
}
\end{equation}
(i.e. vertices are $\{0,1,2,\dots\}$ and for each $i>0$ there is an
arrow $i\to 0$). Putting the arrows in degree one defines a non-negative
grading $\mathsf{A}$ on $A$. However, this example is different from
the previous ones because there are infinitely many arrows pointing to the
vertex $0$. Hence the morphism space from $\mathsf{P}(0)$ to $\oplus_{i\geq
  0}\mathsf{P}(i)$ is infinite-dimensional in degree one. This
infinite-dimensionality makes some of our arguments inapplicable. Hence we will
avoid such situations in our paper by considering locally
bounded categories (condition (C-\ref{pgcond.5}) in
Subsection~\ref{s2.1}, cf. e.g. \cite[2.1]{BoGa}).

\subsubsection*{User's manual}

The following section contains basic definitions and results on graded
categories which are crucial for the general approach but quite
technical. Therefore, at the first reading attempt, we suggest to skip all the
details from Section~\ref{s2} and carry on with Section~\ref{section2}. Since our paper is rather long and contains lots of notation
for objects of rather different nature, we tried, for the readers
convenience, to organize our notation in a way as unified as possible via
different fonts. Of course there are exceptions due to already
well-established notation in the literature, but otherwise the
general convention for notation in the paper is as follows:

\begin{center}
\begin{tabular}{|l|l|}
\hline
Object & Notation\\
\hline
\hline
Algebras: & $A$, $B$, $C$,\dots\\
\hline
Graded algebras: & $\mathsf{A}$, $\mathsf{B}$, $\mathsf{C}$,\dots\\
\hline
Categories: & $\mathscr{A}$, $\mathscr{B}$, $\mathscr{C}$,\dots\\
\hline
Graded categories: & $\mathbf{A}$, $\mathbf{B}$, $\mathbf{C}$,\dots\\
\hline
Modules: & $M$, $N$, $L$,\dots\\
\hline
Graded modules: & $\mathsf{M}$, $\mathsf{N}$, $\mathsf{L}$,\dots\\
\hline
Complexes: & $\mathcal{X}^{\bullet}$, $\mathcal{Y}^{\bullet}$, $\mathcal{Z}^{\bullet}$,\dots\\
\hline
Functors: & $\mathrm{F}$,  $\mathrm{G}$, $\mathrm{K}$,\dots\\
\hline
Derived categories and functors: & $\mathcal{D}$,
$\mathcal{L}$, $\mathcal{R}$,\dots\\
\hline
Dualities: & $\mathbb{D}$,  $\mathbf{d}$,\dots\\
\hline
Objects in categories: & $\lambda$,  $\mu$, $\nu$,\dots\\
\hline
Idempotents: & $e$,  $e_{\lambda}$,\dots\\
\hline
\end{tabular}
\end{center}

\subsection*{Acknowledgments}

We would like to thank Joseph Chuang for posing a question (at a
conference in May 2005) which became the starting point for the
present paper. We are deeply grateful to Bernhard Keller for
many helpful suggestions and discussions. Special thanks to the referee 
for all her/his work and suggestions to improve the paper, and for 
fruitful additional mathematical discussions which we highly appreciate.

\section{Preliminaries}\label{s2}

For the whole paper we fix an arbitrary field $\Bbbk$. Throughout the
paper {\em graded} means {\em $\Z$-graded},  and {\em algebra} means,
if not otherwise stated, a {\em unital $\Bbbk$-algebra} with unit $1$;
$\dim$ means $\dim_{\Bbbk}$, and a {\em category} means a
{\em small category}. For any category $\mathscr{A}$ we
denote by $\op{Ob}(\mathscr{A})$ the set of objects of $\mathscr{A}$
and often just write $\la\in\mathscr{A}$ if $\la\in\op{Ob}(\mathscr{A})$.
For $\la$, $\mu\in\mathscr{A}$ the morphisms from $\la$ to $\mu$ are
denoted $\mathscr{A}(\la,\mu)$. We denote by $\mathscr{A}^{\mathrm{op}}$
the opposite category, that is $\mathscr{A}^{op}(\la,\mu)=\mathscr{A}(\mu,\la)$.
If not stated otherwise, functors are always covariant.

\subsection{Graded algebras and graded categories}\label{s2.1}

Let $\mathbf{C}$ be a $\Bbbk$-linear
category with set of objects $\Ob(\mathbf{C})$. Let
$e_\la\in\mathbf{C}(\lambda,\lambda)$ be the identity morphism.
Recall that the category $\mathbf{C}$ is called {\em graded}
provided that the morphism spaces are graded, that is
$\mathbf{C}(\lambda,\mu)=\oplus_{i\in\Z}
\mathbf{C}_i(\lambda,\mu)$ such that
$\mathbf{C}_i(\mu,\nu)\mathbf{C}_j(\lambda,\mu)\subseteq
\mathbf{C}_{i+j}(\lambda,\nu)$ for all
$\lambda,\mu,\nu\in\Ob(\mathbf{C})$ and $i,j\in\Z$.

{\it A standard example:} A standard example of a graded category is the category with objects
finite-dimensional graded $\Bbbk$-vector spaces and morphisms the
$\Bbbk$-linear maps. The $i$-th graded part is given then by graded maps which
are homogeneous of degree $i$.

{\it A rather naive example:} To any graded $\Bbbk$-algebra $\mathsf{A}=\oplus_{i\in\Z}\mathsf{A}_i$ one can associate, in a rather naive
way, the graded category $\mathbf{C}^{\mathsf{A}}$ containing one single
object, namely $\mathsf{A}$. The morphisms in this category are given by
putting $\mathbf{C}^{\mathsf{A}}_i(\mathsf{A},\mathsf{A})=\mathsf{A}_i$
for all $i\in \Z$ with compositions given by the multiplication
in $\mathsf{A}$. This example will not be very important for us,
although it appears quite often in the literature, namely whenever a
discrete group is considered as a category (groupoid) with one object
and morphisms given by the elements of the group and composition given
by the group multiplication. If $\mathbf{C}$ is a graded
$\Bbbk$-linear category, one can consider
$\oplus_{\lambda,\mu\in\Ob(\mathbf{C})}\mathbf{C}(\lambda,\mu)$, which
is a graded $\Bbbk$-algebra, however without a unit element if
$|\Ob(\mathbf{C})|=\infty$. As already mentioned in the introduction, this
procedure does not have a uniquely defined inverse in general.

\subsubsection*{From graded categories to quotient categories and vice versa}

Graded categories appear as quotient categories modulo
free $\Z$-actions. There is even a correspondence
\begin{eqnarray}
\label{quots}
\begin{array}{ccc}
\left\{\text{categories with a free $\mZ$-action}\right\}&\leftrightarrow&
\left\{\text{graded categories}\right\}\\
\mathscr{C}&\mapsto&\mathscr{C}/\mZ\\
\mathscr{C}^\mZ&\mapsfrom&\mathscr{C}
\end{array}
\end{eqnarray}
constructed in the following way: Let $\mathscr{C}$ be a
$\Bbbk$-category. Assume that the group $\Z$ acts freely on
$\mathscr{C}$ via automorphisms (here {\em freely} means that the
stabilizer of every object is trivial). In this case we
can define the quotient category $\mathscr{C}/\Z$, whose
objects are the orbits of $\Z$ on  $\Ob(\mathscr{C})$, and for
$\lambda,\mu\in \Ob(\mathscr{C})$ the  morphism set
$\mathscr{C}/\Z(\Z\lambda,\Z\mu)$ is defined as the quotient of
\begin{displaymath}
\bigoplus_{\substack{\lambda'\in\Z \lambda\\ \mu'\in\Z\mu}}
\mathscr{C}(\lambda',\mu')
\end{displaymath}
modulo the subspace, generated by all expressions $f-i\cdot f$,
where $i\in \Z$. The product of morphisms is defined in the obvious way. For
any $\lambda,\mu\in\Ob(\mathscr{C})$ we have
a canonical isomorphism of vector spaces,
\begin{displaymath}
\mathscr{C}/\Z(\Z\lambda,\Z\mu)\cong
\bigoplus_{\lambda'\in \Z\lambda}
\mathscr{C}(\lambda',\mu),
\end{displaymath}
which turns $\mathscr{C}/\Z$ into a graded category.
Conversely, let $\mathscr{C}$ be a graded $\Bbbk$-category.
Then we can consider the category $\mathscr{C}^{\Z}$ such that
$\Ob(\mathscr{C}^{\Z})=\Ob(\mathscr{C})\times{\Z}$, and
for $\lambda,\mu\in\Ob(\mathscr{C})$ and
$i,j\in\Z$ we have $\mathscr{C}^{\Z}((\lambda,i),(\mu,j))=
\mathscr{C}_{j-i}(\lambda,\mu)$. Then $\Z$ acts freely on
$\mathscr{C}^{\Z}$ in the obvious way and we have
$\mathscr{C}^{\Z}/\Z\cong \mathscr{C}$ as graded categories
(for details we refer the reader, for instance, to
\cite[Section~2]{CM}).

\subsubsection*{Positively graded algebras and categories}

In the following it will be useful to strengthen the definition of a graded
category and replace it by the notion of a {\it positively} graded category
defined as follows (compare with the notion of locally bounded
categories in \cite[2.1]{BG}):

\begin{definition}\label{defpgc}\rm
 A graded $\Bbbk$-category $\mathbf{C}$ is
 said to be {\em positively graded} provided that the
 following conditions are satisfied:
\begin{enumerate}[({\text C}-i)]
\item\label{pgcond.1}
$\mathbf{C}_i(\lambda,\mu)=0$ for all
$\lambda,\mu\in\Ob(\mathbf{C})$ and $i<0$.
\item\label{pgcond.2}
$\mathbf{C}_0(\lambda,\mu)=
\begin{cases}
0,& \text{ if } \lambda\neq \mu,\\
\Bbbk e_{\lambda},& \text{ if }\lambda=\mu.
\end{cases}
$
\item\label{pgcond.4}
$\dim\mathbf{C}_i(\lambda,\mu)<\infty$ for all
$\lambda,\mu\in\Ob(\mathbf{C})$ and $i\in\Z$.
\item \label{pgcond.5}
For any $\la\in\Ob(\mathbf{C})$, $i\in\mZ$, the sets $\{\mu\mid\mathbf{C}_i(\lambda,\mu)\not=\{0\}\}$
and $\{\mu\mid\mathbf{C}_i(\mu,\lambda)\not=\{0\}\}$ are finite.
\end{enumerate}
\end{definition}

A semi-simple category is always positively graded,
whereas the category of all finite-dimensional
graded $\Bbbk$-vector spaces is not positively graded
(both (C-\ref{pgcond.1}) and (C-\ref{pgcond.2}) fail).
For other examples of positively graded categories we refer to
the introduction, where also one finds an example of
a category,  which does not satisfy the condition~(C-\ref{pgcond.5}).
We remark that a positively graded category is in reality
{\em non-negatively graded} (since (C-\ref{pgcond.1}) only says that
all negatively graded components are zero), however, the use of the term
{\em positively graded} in this context is now commonly accepted
(see for example \cite[2.3]{BGS} or \cite[Introduction]{MVS}).
Positively graded categories with finitely many objects come along with
positively graded algebras:

\begin{definition}\label{defpga}\rm
A graded algebra, $\mathsf{A}=\oplus_{i\in\Z}\mathsf{A}_i$,
is said to be {\em positively graded} provided that the following
conditions are satisfied:
\begin{enumerate}[(\text{A}-i)]
\item\label{pgcond.7} $\dim\mathsf{A}_i<\infty$ for all $i\in\Z$.
\item\label{pgcond.8} $\mathsf{A}_i=0$ for all $i<0$.
\item\label{pgcond.9} $\mathsf{A}_0=\oplus_{\lambda\in\Lambda}\Bbbk e_{\lambda}$, where $1=\sum_{\lambda\in\Lambda}e_{\lambda}$ is a
(fixed) decomposition of the unit element $1$ into a finite sum of
pairwise  orthogonal primitive idempotents.
\end{enumerate}
\end{definition}

As already mentioned in the introduction, to any positively
graded algebra $\mathsf{A}$ one associates a positively
graded $\Bbbk$-category, which we denote by $\mathbf{A}$.
The objects of this category are $\Ob(\mathbf{A})=\Lambda$
(one can also interpret these objects as indecomposable projective
right $\mathsf{A}$-modules), and  the morphisms are defined by setting $\mathbf{A}_i(\mu,\lambda)=e_{\la}\mathsf{A}_ie_{\mu}$ for all
$\lambda,\mu\in\Ob(\mathbf{A})$ and $i\in\Z$ (in other words, the
morphisms are just the homomorphisms between the
corresponding projective modules).  The composition of morphisms
in $\mathbf{A}$ is induced by the multiplication in $\mathsf{A}$.
The condition~(C-\ref{pgcond.5}) is satisfied, since we have
finitely many objects.  Conversely, for any positively graded
$\Bbbk$-linear category $\mathbf{C}$ with finitely many objects the space
$\oplus_{\lambda,\mu\in\Ob(\mathbf{C})}\mathbf{C}(\lambda,\mu)$ is a
positively graded $\Bbbk$-algebra. These two processes
restrict naturally to the correspondence described in
\eqref{introequiv} which will always be in the background of
our considerations.  However, the setup of graded categories is
more general, since we also allow $|\Ob(\mathbf{C})|=\infty$.

\subsection{Modules over graded categories}\label{s2.2}

We have seen that positively graded categories correspond to positively
graded algebras in the sense of \eqref{introequiv} and should be
thought of being the correct framework to deal with not necessarily unital
algebras equipped with some fixed complete set of pairwise commuting
idempotents (for example path algebras of not necessarily
finite quivers). We therefore also introduce the notion of modules
over graded categories which provides the usual definition of modules
over an algebra under the correspondence~\eqref{introequiv} as
explained in the introduction. We denote
\begin{itemize}
\item by $\Bbbk\defis\Mod$ the category of {\em all} $\Bbbk$-vector spaces;
\item by $\Bbbk\defis\mod$ the category of all {\em finite-dimensional}
$\Bbbk$-vector spaces;
\item by $\Bbbk\defis\gMod$ the category of all {\em graded}
$\Bbbk$-vector spaces;
\item by $\Bbbk\defis\gmod$ the category of all {\em finite-dimensional
graded} $\Bbbk$-vector spaces.
\item by $\Bbbk\defis\fgmod$ the category of all graded
$\Bbbk$-vector spaces with {\em finite-di\-men\-si\-o\-nal graded components}.
\end{itemize}
Let $\mathbf{C}$ be a graded category. A $\Bbbk$-linear  functor
$\mathrm{F}:\mathbf{C}\rightarrow\Bbbk\defis\gMod$ is called
{\it homogeneous of degree $d$} if it maps morphisms of degree $k$ to
morphisms of  degree $k+d$ for all $k\in\Z$. In particular,
homogeneous functors  of degree $0$ preserve the degree of morphisms.
A natural transformation  between homogeneous functors is by definition
grading preserving. We define
\begin{itemize}
\item the category $\mathbf{C}\defis\Mod$ of all
{\em $\mathbf{C}$-modules}, as the category of all $\Bbbk$-linear
functors from $\mathbf{C}$ to $\Bbbk\defis\Mod$;
\item the category $\mathbf{C}\defis\gMod$ of all {\em graded
$\mathbf{C}$-modules}, as the category of all $\Bbbk$-linear
homogeneous functors of degree $0$ from $\mathbf{C}$ to
$\Bbbk\defis\gMod$.
\item the category $\mathbf{C}\defis\fmod$ of
{\em locally finite-dimensional $\mathbf{C}$-modules}, as the
category of all $\Bbbk$-linear functors from $\mathbf{C}$ to
$\Bbbk\defis\mod$.
\item the category $\mathbf{C}\defis\fgmod$ of
{\em locally finite-dimensional graded $\mathbf{C}$-modules},
as the category of all $\Bbbk$-linear homogeneous functors
of degree $0$ from $\mathbf{C}$ to $\Bbbk\defis\fgmod$.
\item the category $\mathbf{C}\defis\fdmod$ of
{\em finite-dimensional $\mathbf{C}$-modules}, as the
category of all $\Bbbk$-linear functors from $\mathbf{C}$ to
$\Bbbk\defis\mod$ satisfying the condition that the value of
such functor is non-zero only on finitely many objects from 
$\mathbf{C}$.
\item the category $\mathbf{C}\defis\fdgmod$ of
{\em finite-dimensional graded $\mathbf{C}$-modules},
as the category of all $\Bbbk$-linear homogeneous functors
of degree $0$ from $\mathbf{C}$ to $\Bbbk\defis\gmod$ 
satisfying the condition that the value of such functor is 
non-zero only on finitely many objects from  $\mathbf{C}$.
\end{itemize}
Similarly, we define the corresponding categories of right
${\mathbf C}$-modules via the opposite category
${\mathbf C}^{\mathrm{op}}$ (or, equivalently, using
contravariant functors instead of covariant).

\subsubsection*{Graded modules over graded categories and modules over
  quotient categories}
There is (see for example \cite[Section~2]{CM}) an equivalence
of categories
\begin{eqnarray}
  \label{eq:EC}
\mathrm{E}_{\mathbf{C}}:\quad\mathbf{C}\defis
\gMod&\tilde\longrightarrow&\mathbf{C}^{\Z}\defis\Mod
\end{eqnarray}
which is induced by the correspondence \eqref{quots} and explicitly given as follows: For a graded $\mathbf{C}$-module, that is a functor
$\mathsf{M}:\mathbf{C}\to \Bbbk\defis\gMod$, and for each object
$(\lambda,i)\in \Ob(\mathbf{C}^\Z)$ we set
$\mathrm{E}_{\mathbf{C}}(\mathsf{M})(\lambda,i)=\mathsf{M}(\lambda)_i$,
where $\mathsf{M}(\lambda)=
\oplus_{i\in\Z}\mathsf{M}(\lambda)_i$. For every
$f\in\mathbf{C}^{\mathbb{Z}}((\lambda,i),(\mu,j))=
\mathbf{C}_{j-i}(\lambda,\mu)$ we define $\mathrm{E}_{\mathbf{C}}(\mathsf{M})(f)=\mathsf{M}(f)_i$,
which is a map from
$\mathrm{E}_{\mathbf{C}}(\mathsf{M})(\lambda,i)=\mathsf{M}(\lambda)_i$ to
$\mathsf{M}_{i+(j-i)}(\mu)=
\mathsf{M}_j(\mu)=\mathrm{E}_{\mathbf{C}}(\mathsf{M})(\mu,j)$.
This defines a functor,
$\mathrm{E}_{\mathbf{C}}(\mathsf{M}):\mathbf{C}^\mZ\rightarrow
\Bbbk\defis\Mod$,
or, in other words, an object in $\mathbf{C}^\mZ\defis\Mod$. The
assignment $\mathsf{M}\mapsto\mathrm{E}_{\mathbf{C}}(\mathsf{M})$
specifies what the functor $\mathrm{E}_{\mathbf{C}}$ does on the
level of objects. If $\varphi:\mathsf{M}\to \mathsf{N}$ is a
homomorphism of graded modules, for every $\lambda\in \Ob(\mathbf{C})$ and
$i\in\mathbb{Z}$ we define $\mathrm{E}_{\mathbf{C}}(\varphi)(\lambda,i)$
to be $\varphi_{\lambda,i}:\mathsf{M}(\lambda)_i\to \mathsf{N}(\lambda)_i$
which is the restriction of the map  $\varphi_\lambda$ to the the
$i$-th graded component. This defines the functor $\mathrm{E}_{\mathbf{C}}$.
For the inverse functor $\mathrm{E}_{\mathbf{C}}^{-1}$ and
a $\mathbf{C}^{\Z}$-module $M$ we have
$\mathrm{E}_{\mathbf{C}}^{-1}(M)(\lambda)=
\oplus_{i\in\Z} M\big((\lambda,i)\big)$
for any $\lambda\in\Ob(\mathbf{C})$. If $f\in\mathbf{C}(\la,\mu)$ is
homogeneous of degree $j$ then for $i\in\Z$
we use the identification
$\mathbf{C}_{j}(\la,\mu)=\mathbf{C}^\Z\big((\la,i),(\mu,j+i)\big)$
to get the element
$f(i)$, corresponding to $f$. Then we have
$\mathrm{E}_{\mathbf{C}}^{-1}(M)(f)=
\oplus_{i\in\Z}M(f(i))$. If
$\varphi:M\to N$ is a natural transformation, we put
$\big(\mathrm{E}_{\mathbf{C}}^{-1}(\varphi)\big)_{\la}=
\bigoplus_{i\in\Z}\varphi_{\la,i}$. It is straight-forward to check
that these assignments define inverse equivalences of categories. For more details we refer the reader to  \cite{CM}. Obviously, the functor~\eqref{eq:EC}
restricts to an equivalence of categories
\begin{eqnarray}
  \label{eq:EC2}
\mathrm{E}_{\mathbf{C}}:\quad\mathbf{C}\defis
\fgmod&\tilde\longrightarrow&\mathbf{C}^{\Z}\defis\fmod.
\end{eqnarray}

Let $\mathsf{A}$ be a positively graded algebra and $\mathbf{A}$
the corresponding positively graded category. We define the following categories of
$\mathsf{A}$-modules and leave it as an exercise for the reader to check that
they coincide with our previous definitions under the
equivalence~\eqref{introequiv}.
\begin{itemize}
\item $\mathsf{A}\defis\Mod:=\mathbf{A}\defis\Mod$;
\item $\mathsf{A}\defis\gMod:=\mathbf{A}\defis\gMod$;
\item $\mathsf{A}\defis\mod$ as the category of all
{\em finitely-generated} $\mathsf{A}$-modules;
\item $\mathsf{A}\defis\gmod$ as the category of all
{\em finitely-generated graded} $\mathsf{A}$-modules;
\item $\mathsf{A}\defis\fgmod$ as the category of all
graded $\mathsf{A}$-modules with {\em finite-\-di\-men\-si\-o\-nal
graded  components}.
\end{itemize}

For a positively graded category $\mathbf{C}$ and
$i\in\Z$ we denote by $\langle i\rangle:\mathbf{C}\defis\gMod\to
\mathbf{C}\defis\gMod$ the functor of shifting the grading,
defined as follows: for objects $\lambda\in\Ob(\mathbf{C})$ we have
$\mathsf{M}\langle i\rangle(\lambda)_j=\mathsf{M}(\lambda)_{i+j}$
for all $j\in\Z$. On morphisms, the functor $\langle i\rangle$ is
defined in the obvious (trivial) way.

For a positively graded category $\mathbf{C}$ the category
$\mathbf{C}^{\mathrm{op}}$ inherits a positive
grading in the natural way, namely
$\mathbf{C}_i^{\mathrm{op}}(\la,\mu)=\mathbf{C}_i(\mu, \la)$
for any $\la$, $\mu\in\mathbf{C}$, and $i\in\mathbb{Z}$. If $f\in\mathbf{C}_i^{\mathrm{op}}(\mu,\la)=\mathbf{C}_i(\la,\mu)$, and
$g\in\mathbf{C}_j^{\mathrm{op}}(\nu,\mu)=\mathbf{C}_i(\mu,\nu)$ then
$f\circ^{\mathrm{op}} g=g\circ f\in\mathbf{C}_{i+j}(\la,\nu)=
\mathbf{C}_{i+j}^{\mathrm{op}}(\nu,\la)$.

\subsubsection*{Bimodules, tensor products, $\mathrm{Hom}$ functors,
and dualities}

If $\mathscr{A}$ and $\mathscr{B}$ are two $\Bbbk$-linear
categories, then an {\em $\mathscr{A}\defis\mathscr{B}$-bimodule}
is by definition an
$\mathscr{A}\otimes_{\Bbbk}\mathscr{B}^{\mathrm{op}}$-module,
where
\begin{displaymath}
\begin{array}{rcl}
\Ob(\mathscr{A}\otimes_{\Bbbk}\mathscr{B}^{\mathrm{op}})&=&
\Ob(\mathscr{A})\times \Ob(\mathscr{B}^{\mathrm{op}}),\\
\mathscr{A}\otimes_{\Bbbk}\mathscr{B}^{\mathrm{op}}((\lambda,\mu),
(\lambda',\mu'))&=&
\mathscr{A}(\lambda,\lambda')\otimes_{\Bbbk}
\mathscr{B}(\mu',\mu)
\end{array}
\end{displaymath}
for all $\lambda,\lambda'\in\Ob(\mathscr{A})$ and
$\mu,\mu'\in\Ob(\mathscr{B})$.

Given an $\mathscr{B}^{\mathrm{op}}$-module $X$ and
and a $\mathscr{B}$-module $Y$ we define the tensor product
$X\otimes_\mathscr{B}Y$ as the vector space
$\oplus_{\la,\mu\in\mathscr{B}}X(\la)\otimes_\Bbbk{Y}(\mu)$
modulo the subspace $W$, which is generated by all the elements
${X}(b)(v)\otimes w-v\otimes {Y}(b)(w)$, where
$v\in{X}(\la)$, $w\in{Y}(\mu)$ and
$b\in\mathscr{B}(\mu,\la)=\mathscr{B}^{\mathrm{op}}(\la,\mu)$.
If ${X}$ was an
$\mathscr{A}\defis\mathscr{B}$-bimodule,
then the tensor product $X\otimes_\mathscr{B}Y$ is the
$\mathscr{A}$-module, which assigns to $a\in\Ob(\mathscr{A})$ the
vector space ${X}(a,_-)\otimes_\mathscr{B}Y$ (and the obvious
assignment on morphisms). One can easily check that this corresponds
exactly to the usual tensor product of (bi)modules under
the correspondence~\eqref{introequiv}.

For two $\mathscr{A}$-modules $X$ and $Y$ the set
$\mathscr{A}\defis\Mod({X},{Y})$ is obviously a vector space.
If ${X}$ is an $\mathscr{A}\defis\mathscr{B}$-bimodule, we
define the $\mathscr{B}$-module $\mathscr{A}\defis\Mod({X},{Y})$
in the following way: to any object $b$ from $\mathscr{B}$ we
assign the vector space $\mathscr{A}\defis\Mod({X}({}_-, b),{Y})$,
and to each $f\in \mathscr{B}(b,b')$ we assign
the map, which maps $g=(g_a)_{a\in \Ob(\mathscr{A})}
\in \mathscr{A}\defis\Mod({X}({}_-, b),{Y})$ to
$h=(h_a)_{a\in \Ob(\mathscr{A})}
\in \mathscr{A}\defis\Mod({X}({}_-, b'),{Y})$,
where $h_a=g_a\circ {X}(e_a, f)$.
Again, one checks that this corresponds exactly to the usual
homomorphism construction under the correspondence~\eqref{introequiv}.

It is straightforward to check that for
a $\mathscr{B}$-module $Y$, an $\mathscr{A}$-module $Z$, and
an $\mathscr{A}\defis\mathscr{B}$-bimodule $X$
we have the usual functorial adjunction isomorphism
\begin{eqnarray*}
\mathscr{A}\defis\Mod({X}\otimes_\mathscr{B}{Y},
{Z})&\cong&\mathscr{B}\defis\Mod
\big({Y},\mathscr{A}\defis\Mod({X},{Z})\big)\\
\varphi=\{\varphi_a\}_{a\in\Ob(\mathscr{A})}&\mapsto&\hat{\varphi}=
\{\hat{\varphi}_b\}_{b\in\Ob(\mathscr{B})},\\
\check{\psi}=\{\check{\psi}_a\}_{a\in\Ob(\mathscr{A})}&
\mapsfrom&{\psi}=\{{\psi}_b\}_{b\in\Ob(\mathscr{B})},
\end{eqnarray*}
where for any $a\in\Ob(\mathscr{A})$, $b\in\Ob(\mathscr{B})$,
$m\in{Y}(b)$ and  $x\in {X}(a,b)$ we have
$\hat{\varphi}_b(m)(x)=\varphi_a(x\otimes
m)\in{Z}(a)$, and  $\check{\psi}_a(x\otimes m)=\psi_b(x)(m)\in{Z}(a)$.

Let $\mathbf{d}:\Bbbk\defis\Mod\to\Bbbk\defis\Mod$ be the usual
duality functor $\Bbbk\defis\Mod({}_-,\Bbbk)$. We also have the graded
duality $\mathbb{D}:\Bbbk\defis\gMod\to\Bbbk\defis\gMod$ for which
$(\mathbb{D}\mathsf{V})_i=\mathbf{d}(\mathsf{V}_{-i})$ for any
$\mathsf{V}\in\Bbbk\defis\gMod$ and which acts as the usual duality on
morphisms. Note that if $\mathsf{M}$ is a graded $\mathbf{C}$-module, then
$\mathbb{D}\mathsf{M}$ (defined simply as composition of functors) 
is a $\mathbf{C}^{\mathrm{op}}$-module.

\subsection{The abelian category $\mathbf{C}\defis\fgmod$}\label{s3.1}

Let $\mathbf{C}$ be a positively graded $\Bbbk$-category. The following
statement is obvious, but crucial:
\begin{lemma}
\label{abelian}
For any positively graded $\Bbbk$-category $\mathbf{C}$ the categories
$\mathbf{C}\defis\fgmod$ and $\mathbf{C}\defis\gmod$ are abelian categories.
\end{lemma}

\begin{proof}
  The abelian structure is inherited form the abelian structure of
  $\Bbbk\defis\fgmod$ and $\Bbbk\defis\gmod$. For details we refer to \cite[p.104]{Schubert}.
\end{proof}

\begin{lemma}\label{cover}
Let $\mathbf{C}$ be a positively graded category and
$\la\in\Ob(\mathbf{C})$. Then $\mathsf{P}_{\mathbf{C}}(\lambda)=\mathbf{C}(\lambda,{}_-)$
is an indecomposable projective object in both $\mathbf{C}\defis\gMod$
and $\mathbf{C}\defis\fgmod$.
\end{lemma}

\begin{proof}
By definition we have
$\mathsf{P}(\lambda)=\mathsf{P}_{\mathbf{C}}(\lambda)=
\mathbf{C}(\lambda,{}_-)\in\mathbf{C}\defis\gMod$.  Because of the
assumption~(C-\ref{pgcond.4}), it is even an object of
$\mathbf{C}\defis\fgmod$. It is indecomposable, since the only non-trivial
idempotent of its endomorphism ring is the identity (by the
assumption~(C-\eqref{pgcond.2}) and using the Yoneda lemma).
To see that it is projective, let
$\varphi: \mathsf{F}\rightarrow \mathsf{G}$
and $\alpha:\mathsf{P}(\lambda)\rightarrow \mathsf{G}$ be
morphisms between graded $\mathbf{C}$-modules, where $\varphi$ is
surjective.  We have to show that there is a morphism
$\Phi:\mathsf{P}(\lambda)\rightarrow \mathsf{F}$
such that $\varphi\circ\Phi=\alpha$. Choose $b\in \varphi^{-1}(\alpha(e_{\lambda}))\subset \mathsf{F}(\la)$ and define
$\Phi(f)=\mathsf{F}(f)(b)$ for any
$f\in\mathsf{P}(\lambda)(\mu)=\mathbf{C}(\lambda,\mu)$. Then we have
\begin{eqnarray*}
\varphi(\Phi(f))=(\varphi\circ
\mathsf{F}(f))(b)=(\mathsf{G}(f)\circ\varphi)(b)=
\mathsf{G}(f)(\alpha(e))=\alpha(fe)=\alpha(f).
\end{eqnarray*}
Hence, $\mathsf{P}_{\mathbf{C}}(\lambda)$ is projective and we are done.
\end{proof}

Factoring out the unique maximal graded submodule of $\mathsf{P}(\lambda)$,
that is the submodule given by all elements of positive degree, we obtain
the graded  simple module $\mathsf{L}(\lambda)$. The duality $\mathbb{D}$
maps projective objects to injective objects and preserves
indecomposability, hence we have the graded indecomposable injective
envelope $\mathsf{I}(\lambda)=\mathbb{D}
\mathbf{C}^{\mathrm{op}}(\lambda,{}_-)$ of $\mathsf{L}(\lambda)$.
If we forget the grading, we obtain the ungraded $\mathbf{C}$-modules
$P(\lambda)$, $L(\lambda)$ and $I(\lambda)$ respectively. Note that they are
still indecomposable, and, of course, $L(\lambda)$ is simple.
We define
\begin{eqnarray}
\label{P}
\begin{array}[ht]{lclcl}
\mathsf{P}&=&\mathsf{P}_{\mathbf{C}}&=
&\bigoplus_{\lambda\in\Ob(\mathbf{C})}\mathsf{P}_{\mathbf{C}}(\lambda),\\
\mathsf{I}&=&\mathsf{I}_{\mathbf{C}}&=
&\bigoplus_{\lambda\in\Ob(\mathbf{C})}\mathsf{I}_{\mathbf{C}}(\lambda),\\
\mathsf{L}&=&\mathsf{L}_{\mathbf{C}}&=
&\bigoplus_{\lambda\in\Ob(\mathbf{C})}\mathsf{L}_{\mathbf{C}}(\lambda).
\end{array}
\end{eqnarray}

\begin{lemma}\label{abcgmold}
Let $\mathbf{C}$ be a positively graded category, without necessarily
satisfying (C-\ref{pgcond.5}). The simple objects in $\mathbf{C}\defis\fgmod$
are exactly the modules of the form $\mathsf{L}(\la)\langle i\rangle$,
$\la\in\mathbf{C}$, $i\in\mZ$.  Any object in $\mathbf{C}\defis\fgmod$ of
finite length has a projective cover and an injective hull.
\end{lemma}

\begin{proof}
From Lemma~\ref{abelian} we know that  $\mathbf{C}\defis\fgmod$ is an abelian
category. Let $\mathsf{M}\in\mathbf{C}\defis\fgmod$ be simple. Let
$0\not= v\in\mathsf{M}(\la)_i$ for some $\la$, $i$. Then there is a
non-trivial, hence a surjective,
morphism $\mathsf{P}(\la)\langle -i\rangle\rightarrow \mathsf{M}$ sending
$e_\la$ to $v$. From the positivity of the grading we get
$\mathsf{M}\cong\mathsf{L}(\la)\langle -i\rangle$.
\end{proof}

The assumption (C-\ref{pgcond.5}) in Definition~\ref{defpgc} was introduced to have the following result available:

\begin{lemma}\label{abcgm}
Let $\mathbf{C}$ be a positively graded category.
\begin{enumerate}[(a)]
\item \label{one}
Let $\mathsf{M}$ be an object in $\mathbf{C}\defis\fgmod$.
Assume there exist some
$k\in\mZ$ with the following property: if $j<k$ then
$\mathsf{M}(\mu)_j=\{0\}$ for any
$\mu\in\Ob(\mathbf{C})$. Then $\mathsf{M}$ has a
projective cover in $\mathbf{C}\defis\fgmod$.
\item In particular, any simple object $\mathsf{L}(\lambda)$ has a minimal
  projective resolution.
\item Dually, any simple object $\mathsf{L}(\lambda)$ has a minimal
  injective coresolution.
\end{enumerate}
\end{lemma}

\begin{proof}
Let $\mathscr{A}$ denote the full subcategory of $\mathbf{C}\defis\fgmod$,
which consists of all modules satisfying the the conditions of 
the statement \eqref{one}. The statement \eqref{one} would follow from the
general theory of projective covers, see for example \cite[Proposition~1]{Sh},
provided that we prove two things. Namely, that each object from 
$\mathscr{A}$ is a quotient of some projective object from $\mathbf{C}\defis\fgmod$; and that for any epimorphism 
$f:\mathsf{X}\twoheadrightarrow\mathsf{Y}$ in $\mathscr{A}$ 
there is a minimal submodule $\mathsf{Z}$
of $\mathsf{X}$ with respect to the condition $f(\mathsf{Z})=\mathsf{Y}$.

First we prove that each object from $\mathscr{A}$ is a quotient of some
projective object from $\mathbf{C}\defis\fgmod$
Let $M=\mathrm{E}_{\mathbf{C}}(\mathsf{M})\in\mathbf{C}^\mZ\defis\fmod$.
Define $N\in\mathbf{C}^\mZ\defis\Mod$ as follows:
\begin{eqnarray}
\label{N}
{N}=\bigoplus_{\mu\in\Ob(\mathbf{C}),r\in\mZ}
\bigoplus_{s=1}^{\dim_\Bbbk{{M}\big((\mu,r)\big)}}
\mathbf{C}^\mZ\big((\mu,r),_-\big).
\end{eqnarray}
Note that the second sum just indicates that we take a certain number of
copies of $\mathbf{C}\big((\mu,r),_-\big)$.
Since $\mathsf{M}\in\mathbf{C}\defis\fgmod\cong\mathbf{C}^\mZ\defis\fmod$,
the space ${M}(\mu,r)$ is always finite dimensional, hence the second sum
of \eqref{N} is finite. By the assumption on $\mathsf{M}$ it is
enough to take $r\geq k$.
Since $\mathbf{C}$ is positively graded, we can have
$\mathbf{C}\big((\mu,r),(\la,i)\big)=\mathbf{C}(\mu,\la)_{i-r}\not=\{0\}$
only if $i-r\geq 0$, that is $r\leq i$. Hence we get
\begin{eqnarray}
\label{N2}
{N}\big((\la,i)\big)
&=&\bigoplus_{r=k}^{i}\bigoplus_{\mu\in\Ob(\mathbf{C})}
\bigoplus_{s=1}^{\dim_\Bbbk{{M}\big((\mu,r)\big)}}
\mathbf{C}\big((\mu,r),(\la,i)\big).
\end{eqnarray}
Because of the condition (C-\ref{pgcond.5}), the second sum
appearing in \eqref{N2} in fact produces only a finite number of
non-zero summands. Hence ${N}\big((\la,i)\big)$
is finite dimensional, so
${N}\in\mathbf{C}^\mZ\defis\fmod\cong\mathbf{C}\defis\fgmod$. By
construction, ${N}$ is projective and surjects onto $M$. 

The fact that for any epimorphism $f:\mathsf{X}\twoheadrightarrow\mathsf{Y}$ in 
$\mathscr{A}$ there is a minimal submodule $\mathsf{Z}$
of $\mathsf{X}$ with respect to the condition $f(\mathsf{Z})=\mathsf{Y}$
follows from the definition of $\mathscr{A}$ using the standard arguments
involving Zorn's lemma. Hence the statement \eqref{one} now 
follows \cite[Proposition~1]{Sh}.

The second statement of the lemma follows from the first one and
the remark that the condition on $\mathsf{M}$ in \eqref{one} is
also satisfied for the kernel of the projective
cover of $\mathsf{M}$, constructed above.
The last statement follows by duality.
\end{proof}

For $\lambda\in\Ob(\mathbf{C})$ let $\mathcal{Q}_{\lambda}^{\bullet}$
and $\mathcal{J}_{\lambda}^{\bullet}$ denote a fixed minimal projective
resolution and a minimal injective coresolution of
$\mathsf{L}(\lambda)$ in ${\mathbf C}\defis\fgmod$, respectively. It is easy
to check that such (co)resolutions are unique up to isomorphism. It is not
difficult to see that $\mathbf{C}\defis\gMod$ has enough projectives, whereas $\mathbf{C}\defis\fgmod\cong\mathbf{C}^\mZ\defis\fmod$ does not
need to have enough projectives in general (see the example
of the quivers \eqref{quiv1} and \eqref{quiv101}).

\subsection{Some general notation}\label{s2.3}

In the following we will sometimes write $M_{\mathbf{C}}$ to
indicate that ${M}$ is a (left!) $\mathbf{C}$-module. If $\mathcal{X}^{\bullet}$ is a complex of modules with
differential $d^{\bullet}$,
then $d^i:\mathcal{X}^i\to \mathcal{X}^{i+1}$ for all $i\in\Z$.
If ${M}$ is a module, ${M}^{\bullet}$ will denote the complex
where ${M}^i=0$, $i\neq 0$, and  ${M}^0=M$ with
the trivial differential. For $i\in \Z$ we denote by $[i]$ the functor of
shifting the position in a complex, defined for any complex $\mathcal{X}^{\bullet}$ as follows:
$\mathcal{X}[i]^j=\mathcal{X}^{i+j}$ for all $j\in \Z$.  We denote by
$\mathcal{H}^{i}\mathcal{X}^{\bullet}$ the $i$-th cohomology
of $\mathcal{X}^{\bullet}$. In the hope to avoid confusions we will
use the word {\em degree} for the degree in the grading, and the word
{\em position} for the degree in a complex. An example: if a graded
module $\mathsf{M}$ is concentrated in degree $0$, then
$\mathsf{M}^{\bullet}[i]\langle j\rangle$ is concentrated in
position $-i$ and degree $-j$.

For an abelian category, $\mathscr{A}$, we denote
by $\mathcal{C}(\mathscr{A})$ the category of complexes of
objects from $\mathscr{A}$, by $\mathcal{K}(\mathscr{A})$ its homotopy category, and by  $\mathcal{D}(\mathscr{A})$ the corresponding derived
category. We will use the standard upper indices $b$, $+$, and
$-$, to denote the corresponding categories of bounded,
right bounded and left bounded complexes.
If $\mathscr{A}$ has enough projectives and
$\mathrm{F}:\mathscr{A}\to\mathscr{A}$ is a right exact functor,
we denote by $\mathcal{L}\mathrm{F}$ its left derived functor
and by $\mathcal{L}_i\mathrm{F}$, the $i$-th cohomology
functor of $\mathrm{F}$. Analogously we define
$\mathcal{R}\mathrm{F}$ and $\mathcal{R}^i\mathrm{F}$,
if $\mathrm{F}$ is left exact and $\mathscr{A}$ has enough injectives.
The symbol $\mathrm{ID}$ denotes the identity functor.

For a graded vector space, $V=\oplus_{i\in\Z} V_i$, and
for $j\in\Z$ we denote by $\mathrm{Lev}_j$ the operation
of taking the $j$-th graded component of $V$, that is
$\mathrm{Lev}_j(V)=V_j$.

\begin{figure}[tb]
\special{em:linewidth 0.4pt}
\unitlength 0.80mm
\linethickness{0.4pt}
\begin{center}
\begin{picture}(120.00,40.00)
\drawline(50.00,40.00)(50.00,10.00)
\drawline(51.00,13.00)(50.00,10.00)
\drawline(49.00,13.00)(50.00,10.00)
\drawline(00.00,25.00)(100.00,25.00)
\drawline(97.00,26.00)(100.00,25.00)
\drawline(97.00,24.00)(100.00,25.00)
\drawline(40.00,30.00)(40.00,10.00)
\drawline(40.00,30.00)(10.00,10.00)
\drawline(60.00,20.00)(60.00,40.00)
\drawline(60.00,20.00)(90.00,40.00)
\put(50.00,05.00){\makebox(0,0)[cc]{(grading) degree}}
\put(110.00,25.00){\makebox(0,0)[cc]{position}}
\put(30.00,15.00){\makebox(0,0)[cc]{$\mathscr{A}^{\downarrow}$}}
\put(70.00,35.00){\makebox(0,0)[cc]{$\mathscr{A}^{\uparrow}$}}
\end{picture}
\end{center}
\caption{The supports of objects from the categories $\mathscr{A}^{\downarrow}$ and $\mathscr{A}^{\uparrow}$.}\label{fig1}
\end{figure}
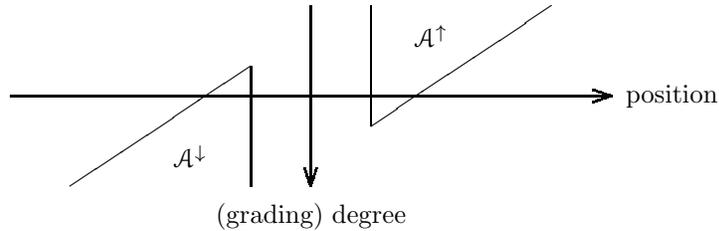

Let $\mathscr{A}$ be an abelian category, whose objects are some
graded modules. Following \cite[2.12]{BGS},
we denote by $\mathcal{C}^{\downarrow}(\mathscr{A})$ the
category of complexes of graded modules from $\mathscr{A}$, which 
consists of all complexes $\mathcal{X}^{\bullet}\in \mathscr{A}$, such 
that there exist integers $N_1(\mathcal{X}^{\bullet})$ and
$N_2(\mathcal{X}^{\bullet})$ satisfying
\begin{equation}\label{eqdar}
\mathrm{Lev}_j(\mathcal{X}^{i})=0
\quad\text{for all}\quad i>N_1(\mathcal{X}^{\bullet})\quad
\text{and all}\quad i+j<N_2(\mathcal{X}^{\bullet});
\end{equation}
and by $\mathcal{C}^{\uparrow}(\mathscr{A})$ the category of complexes of
graded modules from $\mathscr{A}$, such that there exist integers
$N_1(\mathcal{X}^{\bullet})$ and $N_2(\mathcal{X}^{\bullet})$ satisfying
\begin{equation}\label{equar}
\mathrm{Lev}_j(\mathcal{X}^{i})=0
\quad\text{for all}\quad i<N_1(\mathcal{X}^{\bullet})\quad
\text{and all}\quad i+j>N_2(\mathcal{X}^{\bullet}).
\end{equation}
Thus the non-zero components of the objects from
$\mathcal{C}^{\downarrow}(\mathscr{A})$ and 
$\mathcal{C}^{\uparrow}(\mathscr{A})$
are concentrated in regions as depicted in Figure~\ref{fig1}. 
We further denote by $\mathcal{K}^{\downarrow}(\mathscr{A})$,
$\mathcal{K}^{\uparrow}(\mathscr{A})$, $\mathcal{D}^{\downarrow}(\mathscr{A})$,
and $\mathcal{D}^{\uparrow}(\mathscr{A})$ the corresponding homotopy and 
derived categories. Our notation is exactly opposite to the one in \cite[2.12]{BGS}. We made this change, since we think our choice is better
adjusted to the usual terminology that an indecomposable projective module
has a simple head (or top), which in our picture indeed corresponds to 
the highest part of a depicted module.

For a complex $\mathcal{X}^{\bullet}$ and $i\in\Z$ we denote by
$(\mathfrak{t}_i\mathcal{X})^{\bullet}$ the naively
$i$-truncated complex, defined as follows:
$(\mathfrak{t}_i\mathcal{X})^{j}=
\mathcal{X}^{j}$ for all  $j\leq i$, and
$(\mathfrak{t}_i\mathcal{X})^{j}=0$ for all  $j>i$,
with the differential on $(\mathfrak{t}_i\mathcal{X})^{\bullet}$
induced from that on $\mathcal{X}^{\bullet}$.

\section{Categories of linear complexes}\label{section2}

In this section we will introduce one of the main players, the category of
linear complexes (usually of projective modules) as they appear
for example in \cite{MVS}. For Koszul algebras the categories of
linear complexes appeared already in \cite[Corollary~2.13.3]{BGS}
(as cores of non-standard $t$-structure).
Let $\mathsf{M}$ be a graded-structures $\mathbf{C}$-module. We denote by $\mathscr{LC}(\mathsf{M})$ the {\em category of linear complexes
associated with $\mathsf{M}$}, which is defined as follows:
the objects of $\mathscr{LC}(\mathsf{M})$ are all complexes
$\mathcal{X}^{\bullet}$ such that for every $i\in\mathbb{Z}$ every
indecomposable summand of the module $\mathcal{X}^{i}$ occurs with
finite multiplicity and has the form $\mathsf{N}\langle i\rangle$, where
$N$ is an  indecomposable summand of $M$; the morphisms in $\mathscr{LC}(\mathsf{M})$
are all possible morphisms of complexes of graded modules.
In the special case when $\mathsf{M}=\mathsf{P}$ (as defined in \eqref{P}), the category $\mathscr{LC}(\mathsf{M})=\mathscr{LC}(\mathsf{P})$
is called {\em the category of linear complexes of projective modules}.
In the case $\mathsf{M}=\mathsf{I}$, the category
$\mathscr{LC}(\mathsf{M})=\mathscr{LC}(\mathsf{I})$
is called {\em the category of linear complexes of injective modules}.

For $k\in\mZ$ let $\mathscr{LC}(\mathsf{P})^{\geq k}$ be
the full subcategory of $\mathscr{LC}(\mathsf{P})$ given by all complexes $\mathsf{M}$ satisfying $\mathsf{M}^j=\{0\}$ for $j<k$. Obviously,
$\mathscr{LC}(\mathsf{P})=\varprojlim\mathscr{LC}(\mathsf{P})^{\geq k}$,
where the inverse system is given by truncation functors.
Let us recall some basic facts about the categories of linear complexes:

\begin{proposition}\label{pcatlcomp}
\begin{enumerate}[(i)]
\item\label{pcatlcomp.1} Both, $\mathscr{LC}(\mathsf{P})$ and
$\mathscr{LC}(\mathsf{I})$, are abelian categories with the
usual kernels and cokernels for complexes.
\item\label{pcatlcomp.2} The simple objects of
$\mathscr{LC}(\mathsf{P})$ (resp. $\mathscr{LC}(\mathsf{I})$)
are exactly the complexes of the form $\mathsf{P}(\lambda)^{\bullet}\langle -i\rangle[i]$
(resp. $\mathsf{I}(\lambda)^{\bullet}\langle -i\rangle[i]$), where
$\lambda\in\Ob(\mathbf{C})$ and $i\in\Z$.
\item\label{pcatlcomp.3} The Nakayama functor
$\mathrm{N}=\mathrm{N}_{\mathbf{C}}=
(\mathbb{D}\mathbf{C}({}_-,{}_-))\otimes_{\mathbf{C}}{}_- $
induces an equivalence between $\mathscr{LC}(\mathsf{P})$ and
$\mathscr{LC}(\mathsf{I})$. The Nakayama functor satisfies
\begin{displaymath}
\mathrm{N}\big(\mathsf{P}(\la)\langle
i\rangle[-i]\big)\cong\mathrm{N}\big(\mathsf{P}(\la)\big)\langle
i\rangle[-i]\cong \mathsf{I}(\la)\langle i\rangle[-i]
\end{displaymath}
for any $\la\in\Ob(\mathbf{C})$ and $i\in\mZ$.
\end{enumerate}
\end{proposition}

\begin{proof}
The statements \eqref{pcatlcomp.1} and \eqref{pcatlcomp.2} are
proved in \cite[Lemma~5]{MO}. The existence of the equivalence from part~\eqref{pcatlcomp.3} follows from the standard fact that
$\mathrm{N}$ induces an equivalence between the additive closures
of $\mathsf{P}$ and $\mathsf{I}$ with finite multiplicities
(see for example \cite[I.4.6]{Ha}).  The formulas hold by definition.
\end{proof}

\subsection{Projective and injective objects in
$\mathscr{LC}(\mathsf{P})$}\label{s3.3}

The purpose of this section is to give an explicit constructible description
of the indecomposable projective covers and injective hulls of simple objects in the
category $\mathscr{LC}(\mathsf{P})$. These projective and injective objects
exist, although the category does not have enough projectives or enough injectives in general. The analogous results for
$\mathscr{LC}(\mathsf{I})$ can be obtained by applying the Nakayama
automorphism from Proposition~\ref{pcatlcomp}. Recall that for
$\lambda\in\Ob(\mathbf{C})$ we denote by $\mathcal{Q}_{\lambda}^{\bullet}$
(and $\mathcal{J}_{\lambda}^{\bullet}$ respectively) a fixed minimal projective
resolution (and a fixed minimal injective coresolution) of
$\mathsf{L}(\lambda)$, considered as an object of ${\mathbf C}\defis\fgmod$
(see Lemma~\ref{abcgm}). We will show in Proposition~\ref{pltcproj} below how injective (respectively projective)
objects in $\mathscr{LC}(\mathsf{P})$ can be considered as maximal linear
parts of the $\mathcal{Q}_{\lambda}^{\bullet}$'s (respectively of the
images under the inverse Nakayama functor applied to the $\mathcal{J}_{\lambda}^{\bullet}$'s).
We start with some preparation.\\

We will call a complex {\em minimal} provided that it does not contain
any direct summands of the form
\begin{displaymath}
\dots\to 0\to\mathsf{M}\overset{\sim}{\to}\mathsf{M}\to 0\to\dots.
\end{displaymath}
Consider the full subcategory $\overline{\mathcal{C}}(\mathbf{C})$
of the category of complexes of graded  $\mathbf{C}$-modules, whose
objects are all possible minimal complexes $\mathcal{X}^{\bullet}$
such that for every $j\in \mathbb{Z}$ every indecomposable direct
summand of $\mathcal{X}^{j}$ is isomorphic to
$\mathsf{P}(\lambda)\langle k\rangle$ for some
$\lambda\in\Ob(\mathbf{C})$ and some $k\in\Z$.
Denote by $\overline{\mathcal{K}}(\mathbf{C})$ the corresponding
homotopy category.

Fix for the moment $i\in\Z$ and let $\mathcal{X}^{\bullet}\in
\overline{\mathcal{K}}(\mathbf{C})$ with the differential $d^{\bullet}$.
For every $j\in\Z$ we have the following canonical decomposition of
the $\mathcal{X}^{j}$:
$\mathcal{X}^{j}=\mathcal{X}\{>i\}^{j}\oplus
\mathcal{X}\{=i\}^{j}\oplus\mathcal{X}\{<i\}^{j}$,
where all the indecomposable direct summands
\begin{itemize}
\item[] of $\mathcal{X}\{>i\}^{j}$ are isomorphic to
$\mathsf{P}(\lambda)\langle k\rangle$ for some
$\lambda\in\Ob(\mathbf{C})$ and $k>i$,
\item[] of $\mathcal{X}\{=i\}^{j}$ are isomorphic to
$\mathsf{P}(\lambda)\langle i\rangle$ for some
$\lambda\in\Ob(\mathbf{C})$,
\item[] of $\mathcal{X}\{<i\}^{j}$ are isomorphic to
$\mathsf{P}(\lambda)\langle k\rangle$ for some
$\lambda\in\Ob(\mathbf{C})$ and $k<i$.
\end{itemize}

\begin{lemma}\label{l5}
Let $\mathcal{X}^{\bullet}\in\overline{\mathcal{K}}(\mathbf{C})$.
For any  $i,j\in\Z$ we have the following inclusion: $d^j(\mathcal{X}\{>i\}^{j})\subset\mathcal{X}\{>i+1\}^{j+1}$.
\end{lemma}

\begin{proof}
We have of course $d^j(\mathcal{X}^{j})\subset\mathcal{X}^{j+1}$.
Let now $\mathsf{P}(\la)^\bullet\langle k\rangle$ be a summand of
$\mathcal{X}\{>i\}^{j}$, that is $k>i$. Let
$\mathsf{P}(\mu)^\bullet\langle l\rangle$ be a summand of $
\mathcal{X}^{j+1}$ such that $d^j$ induces a non-trivial morphism
\begin{displaymath}
\alpha\in\mathbf{C}\defis\gMod\big(\mathsf{P}(\la)\langle
k\rangle, \mathsf{P}(\mu)\langle
l\rangle)=\mathbf{C}\defis\Mod\big(\mathsf{P}(\la),
\mathsf{P}(\mu))_{l-k}=\mathbf{C}\big(\mu,\la)_{l-k}.
\end{displaymath}
Since $\mathbf{C}$ is positively graded, we have $l\geq k$, hence
\begin{displaymath}
d^j(\mathcal{X}\{>i\}^{j})\subseteq\mathcal{X}\{>i\}^{j+1}=
\mathcal{X}\{>i+1\}^{j+1}\oplus\mathcal{X}\{=i+1\}^{j+1}.
\end{displaymath}
The positivity of the grading also implies that the only indecomposable
direct summands of $\mathcal{X}\{>i\}^{j}$ which can be mapped to
$\mathcal{X}\{=i+1\}^{j+1}$, are the ones isomorphic to
$\mathsf{P}(\lambda)\langle i+1\rangle$ for some
$\lambda\in\Ob(\mathbf{C})$, in which case the corresponding
map must be an isomorphism. This is impossible because of the
minimality of $\mathcal{X}^{\bullet}$. The claim follows.
\end{proof}

Lemma~\ref{l5} allows us to define, depending on some fixed $i\in\mathbb{Z}$,
the following functor (which picks out the part ``supported above the
$i$-shifted diagonal'')
\begin{displaymath}
\begin{array}{rlcl}
\mathrm{S}_i=\mathrm{S}^{\mathbf{C}}_i:&
\overline{\mathcal{K}}(\mathbf{C})&\to&
\overline{\mathcal{K}}(\mathbf{C})\\
& \mathcal{X}^{\bullet}& \mapsto & \mathcal{X}\{>i+\bullet\}^{\bullet},
\end{array}
\end{displaymath}
where the differential on $\mathcal{X}\{>i+\bullet\}^{\bullet}$ is
induced from that on $\mathcal{X}^{\bullet}$ by restriction.
By definition, there is a natural inclusion of functors
$\mathrm{S}_i\hookrightarrow\mathrm{ID}$. We denote by
$\mathrm{Q}_i=\mathrm{Q}^{\mathbf{C}}_i$ the quotient functor.

\begin{lemma}\label{l6}
Let $\mathcal{X}^{\bullet}\in \overline{\mathcal{K}}(\mathbf{C})$
be such that for every $j\in\Z$ each indecomposable summand
of $\mathcal{X}^{j}$ occurs with finite
multiplicity. Then $\mathrm{S}_{-1}\mathrm{Q}_{0}\mathcal{X}^{\bullet}$
is a linear complex of projectives,
hence an object in $\mathscr{LC}(\mathsf{P})$.
\end{lemma}

\begin{proof}
The statement follows directly from the definitions, because, at the position
$j$, the functor $\mathrm{S}_{-1}\mathrm{Q}_{0}\mathcal{X}^{\bullet}$ picks
out the summands of the form $\mathsf{P}(\lambda)\langle k\rangle$,
where $\lambda\in \Ob(\mathbf{C})$ and  $j-1<k\leq j$.
\end{proof}

Note that the functor $\mathrm{S}_{-1}\mathrm{Q}_{0}$ is exactly picking out
the (maximal) linear part of a complex.
Denote by $\mathcal{K}^{\vee}_{\mathbf{C}}$ the full subcategory
of $\overline{\mathcal{K}}(\mathbf{C})$, which consists of all
complexes $\mathcal{X}^{\bullet}\in \overline{\mathcal{K}}(\mathbf{C})$,
such that each indecomposable direct summand occurs with a finite
multiplicity in $\mathcal{X}^j$ for any $j$, and ${\mathbf{C}}\defis\Mod(\mathcal{X}^{i},
\mathsf{L}\langle j\rangle)\neq 0$ implies $j\leq i$ for all
$i,j\in\mathbb{Z}$. Then we have the natural inclusion
$$\op{incl}:\mathscr{LC}(\mathsf{P})\to
\mathcal{K}^{\vee}_{\mathbf{C}}.$$

\begin{lemma}\label{ladjdlc}
The functor $\mathrm{S}_{-1}\mathrm{Q}_{0}:\mathcal{K}^{\vee}_{\mathbf{C}}\to
\mathscr{LC}(\mathsf{P})$ is right adjoint to $\op{incl}$.
\end{lemma}

\begin{proof}
Let $\mathcal{X}^{\bullet}\in \mathscr{LC}(\mathsf{P})$
and $\mathcal{Y}^{\bullet}\in \mathcal{K}^{\vee}_{\mathbf{C}}$.
Since $\mathbf{C}$ is positively graded, using the same arguments
as in the proof of Lemma~\ref{l5}, we have
\begin{displaymath}
{\overline{\mathcal{K}}(\mathbf{C})}
\big(\op{incl}\mathcal{X}^{\bullet},\mathcal{Y}^{\bullet}\big)\cong
{\overline{\mathcal{K}}(\mathbf{C})}
\big(\op{incl}\mathcal{X}^{\bullet},\mathrm{S}_{-1}\mathrm{Q}_{0}\mathcal{Y}^{\bullet}\big)
\overset{\text{\tiny (Lemma~\ref{l6})}}{\cong}
{\mathscr{LC}(\mathsf{P})}
\big(\mathcal{X}^{\bullet},\mathrm{S}_{-1}\mathrm{Q}_{0}\mathcal{Y}^{\bullet}\big).
\end{displaymath}
The claim follows.
\end{proof}

\begin{proposition}\label{pltcproj}
Let $\lambda\in\Ob(\mathbf{C})$.
  \begin{enumerate}[(a)]
  \item \label{existence} The simple object
$\mathsf{P}(\lambda)^{\bullet}$ of
$\mathscr{LC}(\mathsf{P})$ has a projective cover
$\mathcal{P}_{\lambda}^{\bullet}$ and an injective hull
$\mathcal{I}_{\lambda}^{\bullet}$. Hence, any simple object in $\mathscr{LC}(\mathsf{P})$ has a projective cover and
  an injective hull.
\item \label{concrete} There are isomorphisms, of objects from
$\mathscr{LC}(\mathsf{P})$, as follows:
\begin{enumerate}[(i)]
\item\label{pltcproj.1}
$\mathcal{I}_{\lambda}^{\bullet}\cong
\mathrm{S}_{-1}\mathrm{Q}_{0}\mathcal{Q}_{\lambda}^{\bullet}$.
\item\label{pltcproj.2}
$\mathcal{P}_{\lambda}^{\bullet}\cong\mathrm{S}_{-1}\mathrm{Q}_{0}\mathrm{N}^{-1}{\mathcal{J}}_{\lambda}^\bullet$,
where $\mathrm{N}$ is the Nakayama functor from Proposition~\ref{pcatlcomp}.
\end{enumerate}
\end{enumerate}
\end{proposition}

Since it is quite easy to prove Proposition~\ref{pltcproj}\eqref{concrete} assuming the existence of the projective covers and injective hulls as claimed,
we will first give a separate proof for this part. In this proof we
will compare the functors
${\mathscr{LC}(\mathsf{P})}\big({}_-,\mathcal{I}_{\lambda}^{\bullet}\big)$
and ${\mathscr{LC}\big(\mathsf{P})}({}_-,
\mathrm{S}_{-1}\mathrm{Q}_{0}\mathcal{Q}_{\lambda}^{\bullet}\big)$ and show
that they are isomorphic. The second proof is more technical, but provides the
existence as well. It characterizes
$\mathrm{S}_{-1}\mathrm{Q}_{0}\mathcal{Q}_{\lambda}^{\bullet}$ as the unique
object having simple socle $\mathsf{P}(\la)^\bullet$ and being
injective.

\begin{proof}[Proof of Proposition~\ref{pltcproj}\eqref{concrete} assuming
  the existence part~\eqref{existence}.]
Consider the functors
\begin{displaymath}
\mathrm{F}_1:=
{\mathscr{LC}(\mathsf{P})}\big({}_-,\mathcal{I}_{\lambda}^{\bullet}\big),
\quad\quad
\mathrm{F}_2:={\mathscr{LC}(\mathsf{P})}\big({}_-,
\mathrm{S}_{-1}\mathrm{Q}_{0}\mathcal{Q}_{\lambda}^{\bullet}\big).
\end{displaymath}
Since $\mathcal{I}_{\lambda}^{\bullet}$ is the injective hull
of the simple object $\mathsf{P}(\lambda)^\bullet$ in
$\mathscr{LC}(\mathsf{P})$, we have for any $\mathcal{X}^{\bullet}\in
\mathscr{LC}(\mathsf{P})$ the isomorphism
\begin{displaymath}
\mathrm{F}_1(\mathcal{X}^{\bullet})\cong{\mathbf{C}\defis\fgmod}(\mathcal{X}^{0},\mathsf{P}(\lambda))
\cong
\mathbf{d}(\mathrm{Lev}_0(\mathcal{X}^{0}(\lambda))).
\end{displaymath}
On the other hand, since
$\mathcal{X}^{\bullet}$ is a linear complex of projective modules,
applying Lemma~\ref{ladjdlc} we have
\begin{displaymath}
\mathrm{F}_2(\mathcal{X}^{\bullet})\cong
{\mathcal{K}^{\vee}_{\mathbf{C}}}
(\op{incl}\mathcal{X}^{\bullet},\mathcal{Q}_{\lambda}^{\bullet})\cong
\mathcal{K}({\mathbf{C}}\defis\fgmod)
(\op{incl}\mathcal{X}^{\bullet},\mathsf{L}(\lambda)^{\bullet})\cong
\mathbf{d}(\mathrm{Lev}_0(\mathcal{X}^{0}(\lambda))).
\end{displaymath}
Since all the isomorphisms are natural, it follows
that the functors $\mathrm{F}_1$ and
$\mathrm{F}_2$ are isomorphic. Therefore, there must be an isomorphism $\mathcal{I}_{\lambda}^{\bullet}\cong
\mathrm{S}_{-1}\mathrm{Q}_{0}\mathcal{Q}_{\lambda}^{\bullet}$.
The second statement follows then by applying Proposition~\ref{pcatlcomp}.
\end{proof}

\begin{proof}[Proof of Proposition~\ref{pltcproj} including the existence.]

We first note that the implication in part
\eqref{existence} is clear, since if $\mathcal{P}_{\lambda}^{\bullet}$ is a projective cover and
$\mathcal{I}_{\lambda}^{\bullet}$ is an injective hull of the simple
object $\mathsf{P}(\lambda)^{\bullet}$, then $\mathcal{P}_{\lambda}^{\bullet}\langle -i\rangle[i]$ is a projective
cover and  $\mathcal{I}_{\lambda}^{\bullet}\langle -i\rangle[i]$ is
an injective hull of the simple object
$\mathsf{P}(\lambda)^{\bullet}\langle -i\rangle[i]$ and we are
done by Proposition~\ref{pcatlcomp}.

Set $\mathcal{X}^{\bullet}=
\mathrm{S}_{-1}\mathrm{Q}_{0}\mathcal{Q}_{\lambda}^{\bullet}$.
This is an object of ${\mathscr{LC}(\mathsf{P})}$ by Lemma~\ref{l6}.
Using Lemma~\ref{ladjdlc}, we calculate:
\begin{eqnarray*}
{\mathscr{LC}(\mathsf{P})}\big(\mathsf{P}(\mu)^\bullet
\langle -i\rangle[i],\mathcal{X}^{\bullet}\big)
&=&
{\mathscr{LC}(\mathsf{P})}
\big(\mathsf{P}(\mu)^\bullet\langle
-i\rangle[i],\mathrm{S}_{-1}\mathrm{Q}_{0}
\mathcal{Q}_{\lambda}^{\bullet}\big)\\
&\cong&
\mathcal{K}^{\vee}_{\mathbf{C}}
\big(\op{incl}\mathsf{P}(\mu)^\bullet\langle
-i\rangle[i],\mathcal{Q}_{\lambda}^{\bullet}\big).
\end{eqnarray*}
From the definition of $\mathcal{Q}_{\lambda}^{\bullet}$ we
therefore get the
following: For $\mu\in\Ob(\mathbf{C})$ and $i\in\Z$ we have
\begin{displaymath}
{\mathscr{LC}(\mathsf{P})}
\big(\mathsf{P}(\mu)^\bullet\langle -i\rangle[i],\mathcal{X}^{\bullet}\big)=
\begin{cases}
\Bbbk, & \text{ if $\mu=\lambda$, $i=0$};\\
0, & \mathrm{otherwise}.
\end{cases}
\end{displaymath}
This implies that $\mathcal{X}^{\bullet}$ has, as an object of
$\mathscr{LC}(\mathsf{P})$, simple socle, namely
$\mathsf{P}(\lambda)^\bullet$. Thus, to complete the proof we just
have to show that $\mathcal{X}^{\bullet}$ is an injective object of
$\mathscr{LC}(\mathsf{P})$. We claim that it is even enough to show that
\begin{equation}\label{eqpltcproj.1}
\Ext_{\mathscr{LC}(\mathsf{P})}^1
(\mathsf{P}(\mu)\langle -i\rangle[i],\mathcal{X}^{\bullet})=0
\end{equation}
for all $\mu\in\Ob(\mathbf{C})$ and $i\in\Z$. Indeed, if we fix $k\in\mZ$
then the formula~\eqref{eqpltcproj.1} implies that
$\Ext_{\mathscr{LC}(\mathsf{P})}^1
(\mathcal{Y}^{\bullet},\mathcal{X}^{\bullet})=0$
for any $\mathcal{Y}^{\bullet}\in\mathscr{LC}(\mathsf{P})^{\geq k}$.
Since $\mathscr{LC}(\mathsf{P})=\varprojlim\mathscr{LC}(\mathsf{P})^{\geq k}$
we are done.

For $i\leq 0$ the formula
\eqref{eqpltcproj.1} is clear. Let us assume $i>0$. Let
$d^{\bullet}$ be the differential in $\mathcal{X}^{\bullet}$, and
$f:\mathsf{P}(\mu)\langle -i\rangle[i]\to \mathcal{X}^{-i+1}$
be a non-zero map such that $d^{-i+1}\circ f=0$.
Let $\mathcal{Y}^{\bullet}=\mathrm{Cone}(f)$ be the cone of $f$.
Let $\mathrm{V}$ denote the kernel of $d^{-i+1}$, restricted
to $\mathrm{Lev}_i(\mathcal{X}^{-i+1})$ and $v\in \mathrm{V}$.
Since
$\mathrm{Lev}_i(\mathcal{H}^{-i+1}\mathcal{Q}_{\lambda}^{\bullet})=0$,
there exists $w\in \mathrm{Lev}_i(\mathcal{Q}_{\lambda}^{-i})$
such that $d^{-i}(w)=v$. However,
$\mathrm{Lev}_i(\mathcal{Q}_{\lambda}^{-i})=
\mathrm{Lev}_i(\mathcal{X}^{-i})$ by construction,
which implies that there exists an indecomposable direct
summand, say $\mathsf{M}$, of $\mathcal{X}^{-i}$, such that
$d^{-i}(\mathsf{M})\cong f(\mathsf{P}(\mu)\langle -i\rangle[i])$.
It follows that $\mathsf{M}\cong \mathsf{P}(\mu)\langle -i\rangle[i]$
and one can find generators, $a\in \mathsf{M}$,
$b\in \mathsf{P}(\mu)\langle -i\rangle[i]$, such that
$d^{-i}(a)=f(b)$. The element $a-b$ thus generates in
$\mathcal{Y}^{-i}$ a $\mathbf{C}$-submodule,
isomorphic to $\mathsf{P}(\mu)\langle -i\rangle[i]$. The
latter belongs to the socle of the complex
$\mathcal{Y}^{\bullet}\in \mathscr{LC}(\mathsf{P})$. Hence
$\mathcal{Y}^{\bullet}$ splits.
This proves \eqref{eqpltcproj.1} for $i>0$. Hence,
$\mathcal{I}_{\lambda}^{\bullet}$ exists and has the required form. The remaining statements follow then by applying Proposition~\ref{pcatlcomp}~\eqref{pcatlcomp.3}.
\end{proof}

\section{Quadratic duality for positively graded categories}\label{se3}

In this section we develop the abstract theory of quadratic duality in
terms of linear complexes. This approach has its origins in \cite{MVS}
and  \cite{MO}.

Recall that a positively graded category $\mathbf{C}$
is said to be  {\em generated in degree one} if any morphism in
$\mathbf{C}$ is a linear combination of either scalars or
compositions of homogeneous morphisms of degree one. Further,
$\mathbf{C}$  is called {\em quadratic} if it is
generated in degree one and any relation for
morphisms in $\mathbf{C}$ follows from relations in degree two.
The purpose of this section is to describe locally finite
dimensional modules over the quadratic dual category in terms of linear
complexes of projectives in the original category. We start by defining
the quadratic dual.

\subsection{The quadratic dual of a positively graded category via
linear complexes of projectives}\label{s3.2}

Let still $\mathbf{C}$ be a positively graded $\Bbbk$-linear
category. Let $\mathbf{C}_0$ be the subcategory  of $\mathbf{C}$ with the
same set of objects but only homogeneous morphisms of degree $0$. Then $\mathbf{C}_1({}_-,{}_-)$ becomes a
$\mathbf{C}_0$-bimodule in the natural way, which also induces
a $\mathbf{C}_0$-bimodule structure on
$\mathbf{V}=\mathbf{d}(\mathbf{C}_1({}_-,{}_-))$.
Therefore, one can define the free tensor bimodule
\begin{displaymath}
\mathbf{C}_0[\mathbf{V}]({}_-,{}_-)=\mathbf{C}_0({}_-,{}_-)\oplus
\mathbf{V}({}_-,{}_-)\oplus
\left(\mathbf{V}({}_-,{}_-)
\otimes_{\mathbf{C}_0}\mathbf{V}({}_-,{}_-)\right)\oplus\dots
\end{displaymath}
and the corresponding category
$\mathbf{F}$, where
$\Ob(\mathbf{F})=\Ob(\mathbf{C})$, and for
$\lambda,\mu\in \Ob(\mathbf{F})$ we have
$\mathbf{F}(\lambda,\mu)=\mathbf{C}_0[\mathbf{V}](\lambda,\mu)$.

For $\lambda,\mu,\nu\in \Ob(\mathbf{C})$ consider the multiplication map
\begin{displaymath}
\mathbf{m}_{\la,\mu}^\nu:\quad\mathbf{C}_1(\nu,\mu)\otimes
\mathbf{C}_1(\la,\nu) \longrightarrow\mathbf{C}_2(\la,\mu),
\end{displaymath}
which gives rise to the dual map
\small
\begin{equation}
\label{dual}
\mathbf{d}(\mathbf{m}_{\lambda,\mu}^{\nu}):
\mathbf{d}(\mathbf{C}_2(\la,\mu))\to
\mathbf{d}\big(\mathbf{C}_1(\nu,\mu)
\otimes\mathbf{C}_1(\la,\nu)\big)\cong
\mathbf{d}(\mathbf{C}_1(\la,\nu))\otimes
\mathbf{d}(\mathbf{C}_1(\nu,\mu)).
\end{equation}
\normalsize
Note that the canonical isomorphism as indicated in \eqref{dual} exists by
Property (C-\ref{pgcond.4}) and \cite[Page~147]{McL}.
We denote by $\mathbf{J}({}_-,{}_-)$ the
subbimodule of $\mathbf{C}_0[\mathbf{V}]({}_-,{}_-)$,
generated by the images of all these maps, and
define the (positively) graded category $\mathbf{C}^!$, called
the {\em quadratic dual} of $\mathbf{C}$, as follows: We just have
$\Ob(\mathbf{C}^!)=\Ob(\mathbf{C})=\Ob(\mathbf{F})$, and for
$\lambda,\mu\in \Ob(\mathbf{C}^!)$ we set
\begin{displaymath}
\mathbf{C}^!(\lambda,\mu)=\mathbf{F}(\lambda,\mu)/\mathbf{J}(\lambda,\mu).
\end{displaymath}
By definition, the quadratic dual is quadratic.

The following statement was proved originally in
\cite[Theorem~2.4]{MVS} for unital algebras, an alternative proof
was given in \cite[Theorem~8]{MO}. The latter one can be adjusted to
the setup of the graded categories:

\begin{theorem}\label{tqdual}
There is an equivalence of categories,
\begin{displaymath}
\epsilon=\epsilon_\mathbf{C}:\quad \mathscr{LC}(\mathsf{P})\cong
\mathbf{C}^!\defis\fgmod
\end{displaymath}
such that $\epsilon\langle i\rangle [-i]\cong \langle-i\rangle \epsilon$.
\end{theorem}

\begin{proof}
We will use the identification \eqref{eq:EC} and define an equivalence
$\epsilon': \mathscr{LC}(\mathsf{P})\cong(\mathbf{C}^!)^\mZ\defis\fmod$.
We start by defining the inverse functor. Let ${X}$ be an object
from $(\mathbf{C}^!)^\mZ\defis\fmod$. In
particular, for any $(\la, i)\in\Ob(\mathbf{C}^\mZ)$ we have
$\op{\dim}{X}(\la,i)<\infty$.
For $i\in\mZ$ let $(\mathsf{M}_{X})^{i}$ be the graded
$\mathbf{C}$-module
\begin{eqnarray}
(\mathsf{M}_{X})^{i}=\bigoplus_{\la\in\Ob(\mathbf{C})}\mathsf{P}(\la)\langle
  i\rangle\otimes{X}(\la,i)
\end{eqnarray}
(this means we just take $\op{\dim}{X}(\la,i)$ many copies
of $\mathsf{P}(\la)\langle i\rangle$).
We consider the graded $\mathbf{C}$-module $(\mathsf{M}_{X})^{i}$ as a
  $\mathbf{C}^\mZ$-module via
  \begin{eqnarray}
\label{DefP}
    \big(\mathsf{P}(\la)\langle
  i\rangle\otimes{X}(\la,i)\big)
  (\nu,k)=\mathbf{C}^\mZ\big((\la,i),(\nu,k)\big)\otimes {X}(\la,i)
  \end{eqnarray}
for any $(\nu,k)\in\Ob(\mathbf{C}^\mZ)$.
We want to construct an object
 $\mathsf{M}_{X}$ in $\mathscr{LC}(\mathsf{P})$ with $i$-component
 $(\mathsf{M}_{X})^{i}$. Any object ${X}$ in
 $(\mathbf{C}^!)^\mZ\defis\fmod$ is uniquely defined by the following data
 describing the module structure:
\begin{enumerate}[(D1)]
\item a collection of finite dimensional vector spaces ${X}(\la,j)$
for any $\la\in\Ob(\mathbf{C})$, $j\in\mZ$; and
\item certain elements
$$f'_{\la,\mu,j}\in\Bbbk\defis\fmod
\Big((\mathbf{C}^!)^\mZ((\la,j),(\mu,j+1)),
\Bbbk\defis\fmod\big({X}(\la,j),{X}(\mu,j+1)\big)\Big)$$
for any $\la, \mu\in\Ob(\mathbf{C})$ and $j\in\mZ$.
\end{enumerate}
Note that it is enough to consider just the action of morphisms of
degree one, since $\mathbf{C}_0[\mathbf{V}]({}_-,{}_-)$ is generated in
degrees zero and one. By the definition of the quadratic dual
we have fixed isomorphisms
\begin{eqnarray*}
(\mathbf{C}^!)^\mZ\big((\la,i),(\mu,i+1)\big)\cong\mathbf{d}
\big(\mathbf{C}\big((\mu,i+1),(\la,i)\big)\big).
\end{eqnarray*}
We get natural isomorphisms as follows
\begin{eqnarray}
\label{adjointnessetc}
&&\Bbbk\defis\fmod\Big((\mathbf{C}^!)^\mZ\big((\la,i),(\mu,i+1)\big),
\Bbbk\defis\fmod\big({X}(\la,i),{X}(\mu,i+1)\big)\Big)\nonumber\\
&\cong&\Bbbk\defis\fmod\Big({X}(\la,i)
\otimes(\mathbf{C}^!)^\mZ
\big((\la,i),(\mu,i+1)\big),{X}(\mu,{i+1})\Big)\nonumber\\
&\cong&\Bbbk\defis\fmod\Big({X}(\la,i),\mathbf{C}^\mZ
\big((\mu,i+1),(\la,i)\big)\otimes {X}(\mu,i+1)\Big).
\end{eqnarray}
We denote by
$f_{\la,\mu,i}\in\Bbbk\defis\fmod
\Big({X}(\la,i),\mathbf{C}^\mZ\big((\mu,i+1),(\la,i)\big)\otimes
{X}(\mu,i+1)\Big)$ the image of $f'_{\la,\mu,i}$ under \eqref{adjointnessetc}.
Hence, ${{X}}$ comes along with this collection $f_{\la,\mu,i}$ of maps
and is uniquely determined by this collection. For any
$(\nu,k)\in\Ob(\mathbf{C^\mZ})$ the map $f_{\la,\mu,i}$ induces a
$\Bbbk$-linear map
\begin{eqnarray}
\label{Ps}
\mathbf{C}^\mZ\big((\la,i),(\nu,k)\big)\otimes {X}(\la,i)&\longrightarrow&
\mathbf{C}^\mZ\big((\mu,i+1),(\nu,k)\big)\otimes {X}(\mu,i+1)\nonumber\\
c\otimes x&\longmapsto& (c\otimes\op{id})(f_{\la,\mu,i}(x)).
\end{eqnarray}
This construction is obviously natural in $(\nu,k)$. Together with the
formula~\eqref{DefP} we therefore get a natural transformation of functors
(that is a morphism of $\mathbf{C}$-modules):
\begin{eqnarray*}
\op{d}_{\la,\mu}^i:\quad\mathsf{P}(\la)\langle
 i\rangle\otimes{{X}}(\la,i)&\longrightarrow&
\mathsf{P}(\mu)\langle i+1\rangle\otimes{X}(\mu,i+1).
\end{eqnarray*}
Taking the direct sum defines a morphism of graded $\mathbf{C}$-modules
$\op{d}^i:(\mathsf{M}_{X})^{i}\rightarrow (\mathsf{M}_{X})^{i+1}$.

We claim that we in fact constructed a complex. For this we have to
consider the compositions
$\op{d}_{\mu,\sigma}^{i+1}\circ \op{d}_{\la,\mu}^i$ for any
$\la,\mu,\sigma\in\Ob(\mathbf{C})$ and $i\in\mZ$. We have to show
that the composition
\small
\begin{eqnarray}
{X}(\la,i)
&\stackrel{f_{\la,\mu,i}}{\longrightarrow}&
\mathbf{C}_1(\mu,\la)\otimes{X}(\mu,{i+1})\nonumber\\
&\stackrel{\op{id}\otimes f_{\mu,\sigma,i+1}}{\longrightarrow}&
\mathbf{C}_1(\mu,\la)\otimes\mathbf{C}_1(\sigma,\mu)
\otimes{X}(\sigma,{i+2})\nonumber\\
&\stackrel{{\bf m}\otimes\op{id}}{\longrightarrow}&\mathbf{C}_2(\sigma,\la)
\otimes{X}(\sigma,{i+2})
\label{form1}
\end{eqnarray}
\normalsize
is zero.
Via the isomorphisms~\eqref{adjointnessetc} it is enough to show that the
following composition
\small
\begin{eqnarray}
\mathbf{d}(\mathbf{C}_2(\sigma,\la))
&\stackrel{\mathbf{d}(\mathbf{m}_{\sigma,\la}^{\mu})}{\longrightarrow}&
\mathbf{d}\big(\mathbf{C}_1(\mu,\la)
\otimes\mathbf{C}_1(\sigma,\mu)\big)\cong\mathbf{C}^!_1(\mu,\sigma)\otimes
\mathbf{C}^!_1(\la,\mu)\nonumber\\
&\stackrel{f'_{\mu,\sigma,j}\otimes f'_{\la,\mu, j}}\longrightarrow&
\Bbbk\defis\fmod({X}(\la,j),{X}(\mu,j+1))\otimes\nonumber\\
&&\hspace{2cm}\Bbbk\defis\fmod({X}(\mu,{j+1}),{X}(\sigma,{j+2}))\nonumber\\
&\stackrel{{\text composition}}\longrightarrow&
\Bbbk\defis\fmod({X}(\la,j),{X}(\sigma,{j+2}))
\label{form2}
\end{eqnarray}
\normalsize
is zero.  The latter is obviously satisfied by the definition of the
quadratic dual category $\mathbf{C}^!$ and the fact that ${X}$ is a
$(\mathbf{C}^!)^{\mathbb{Z}}$-module. Altogether, we defined a functor
$\eta:(\mathbf{C}^!)^{\mathbb{Z}}-\fmod\rightarrow\mathscr{LC}
(\mathsf{P})$.\\

Let $\mathcal{P}_{\lambda}^{\bullet}$ be the projective cover of
$\mathsf{P}(\lambda)^{\bullet}$ in $\mathscr{LC}(\mathsf{P})$ (see
Proposition~\ref{pltcproj}). We define a functor $$\epsilon':\quad\mathscr{LC}(\mathsf{P})
\longrightarrow(\mathbf{C}^!)^\mZ\defis\fmod$$ as follows:
If $\mathcal{M}^{\bullet}$ is an object from
$\mathscr{LC}(\mathsf{P})$ then we define
$\epsilon'(\mathcal{M}^{\bullet})(\la,i)=
\mathscr{LC}(\mathcal{P}_{\lambda}^{\bullet}\langle
-i\rangle[i],\mathcal{M}^{\bullet})$.
From the definitions it follows that
$\op{dim}\mathscr{LC}(\mathcal{P}_{\lambda}^{\bullet}\langle-i\rangle[i],
\mathcal{M}^{\bullet})$
is the multiplicity of $\mathsf{P}(\lambda)^{\bullet}\langle -i\rangle[i]$
in $\mathcal{M}^{\bullet}$. This number is finite by the definition of $\mathscr{LC}(\mathsf{P})$. Since $\mathcal{M}^{\bullet}$ is a complex of
graded $\mathbf{C}$-modules, its differential $d^j$
induces a map
\begin{eqnarray*}
  \Psi:\mathbf{C}\big((\la,i),_-\big)\otimes
   \epsilon'(\mathcal{M}^{\bullet})(\la,i)\rightarrow
  \mathbf{C}\big((\mu,i+1),_-\big)\otimes
   \epsilon'(\mathcal{M}^{\bullet})(\mu,i+1)
\end{eqnarray*}
for any $\la$, $\mu$, $i$, in particular,
\begin{eqnarray*}
  \Psi_{(\la,i)}:\mathbf{C}\big((\la,i),(\la,i)\big)\otimes
  \epsilon'(\mathcal{M}^{\bullet})(\la,i)\rightarrow
   \mathbf{C}\big((\mu,i+1),(\la,i)\big)\otimes
    \epsilon'(\mathcal{M}^{\bullet})(\mu,i+1).
\end{eqnarray*}
Since  $\Psi_{(\la,i)}$ is a morphism of $\mathbf{C}$-modules, it is uniquely
determined by the induced $\Bbbk$-linear map
\begin{eqnarray*}
\epsilon'(\mathcal{M}^{\bullet})(\la,i)\rightarrow
 \mathbf{C}\big((\mu,i+1),(\la,i)\big)\otimes
  \epsilon'(\mathcal{M}^{\bullet})(\mu,i+1).
\end{eqnarray*}
Using formula~\eqref{adjointnessetc} we get a possible data $(D2)$ defining a
$(\mathbf{C}^!)^\mZ$-module structure on $\epsilon'(M)$. Using again the
formulas~\eqref{form2} and \eqref{form1} we get that this is in fact a module
structure. Hence $\epsilon'(\mathcal{M}^{\bullet})$ becomes an object in
$(\mathbf{C}^!)^\mZ\defis\fmod$. From the naturality of the construction
it follows that this defines a functor
$\epsilon':\mathscr{LC}(\mathsf{P})\cong(\mathbf{C}^\mZ)^!\defis\fmod$.
Together with the identification from \eqref{eq:EC} we get an equivalence
$\epsilon$ as asserted in the theorem.

By definition we have $\epsilon\langle i\rangle [-i])\cong \langle
-i\rangle \epsilon$. By construction we
have $\epsilon\eta(\mathsf{X})\cong\mathsf{X}$ and
$\eta \epsilon(\mathcal{M}^{\bullet})\cong\mathcal{M}^{\bullet}$.
Hence $\epsilon$ and $\eta$ are dense. Moreover, by
construction, they are both faithful, hence automatically full as
well. Therefore, $\epsilon$ and $\eta$ are equivalences of categories. (In
fact they are mutually inverses. To see this one has to fix  a minimal
system of representatives for the isomorphism classes of the
indecomposable projective $\mathbf{C}$-modules and work only with projectives
from this system.) The theorem follows.
\end{proof}

For $k\in\mZ$ let $(\mathbf{C}^!)^{\geq k}$ denote the full subcategory of
$\mathbf{C}^!$, whose objects are $(\lambda,i)$, where
$\la\in\Ob(\mathbf{C})$ and $i\geq k$. The
$(\mathbf{C}^!)^{\geq k}\defis\fmod$ can be considered as full subcategories of
$(\mathbf{C}^!)^{\mZ}\defis\fmod$. The inclusions
$(\mathbf{C}^!)^{\geq k+1}\hookrightarrow (\mathbf{C}^!)^{\geq k}$
induce  an inverse system on $(\mathbf{C}^!)^{\geq k}\defis\fmod$
via truncations, and we have
$(\mathbf{C}^!)^{\mZ}\defis\fmod=\varprojlim\big((\mathbf{C}^!)^{\geq
  k}\defis\fmod\big)$.

\begin{corollary}
\label{Pk}
  \begin{enumerate}[(a)]
  \item Let $k\in\mZ$. The equivalence $\epsilon$ restricts to an equivalence
    \begin{eqnarray*}
      \epsilon^{\geq k}: \mathscr{LC}(\mathsf{P})^{\geq
      k}\cong(\mathbf{C}^!)^{\geq k}\defis\fmod.
    \end{eqnarray*}
   \item For any $k\in\mZ$, the category $\mathscr{LC}(\mathsf{P})^{\geq k}$
   has  enough projectives. Moreover, $\mathcal{P}_{\lambda}^{\bullet}\langle
   -i\rangle[i]\in\mathscr{LC}(\mathsf{P})^{\geq k}$ for any $i\leq k$ and
   $\lambda\in\Ob(\mathbf{C})$.
  \end{enumerate}
\end{corollary}
\begin{proof}
  This follows directly from Theorem~\ref{tqdual} and Lemma~\ref{abcgm}.
\end{proof}

Let ${\PK}$ be a fixed minimal projective
generator of $\mathscr{LC}(\mathsf{P})^{\geq k}$.
This is by definition a complex of graded $\mathbf{C}$-modules.
If $\mathcal{M}^{\bullet}$ is a complex in ${\mathscr{LC}(\mathsf{P})}$,
then there is a $(\mathbf{C}^!)^{\geq k}$-module structure on
${\mathscr{LC}(\mathsf{P})}({\PK},\mathcal{M}^{\bullet})$
as follows: to each $\la\in \Ob(\mathbf{C})$ and  $i\geq k$ we assign
the space ${\mathscr{LC}(\mathsf{P})}(\mathcal{P}^{\bullet}_{\lambda}
\langle -i\rangle[i],\mathcal{M}^{\bullet})$
(recall that $\mathcal{P}^{\bullet}_{\lambda}
\langle -i\rangle[i]$ is the direct summand of ${\PK}$,
which corresponds to these $\lambda$ and $i$).

\begin{proposition}
  The functor ${\mathscr{LC}(\mathsf{P})}({\PK},_-)$
  defines an equivalence of categories $\mathscr{LC}(\mathsf{P})^{\geq
      k}\cong(\mathbf{C}^!)^{\geq k}\defis\fmod$.
\end{proposition}
\begin{proof}
  This follows directly form Theorem~\ref{tqdual} and Corollary~\ref{Pk}.
\end{proof}

If we choose the ${\PK}$ such that they give rise to a directed system
we directly get the following result
\begin{corollary}\label{corn543}
The functor
$\mathcal{C}(\mathbf{C}\defis\gMod)(\varinjlim
{\PK},_-)=
\varprojlim{\mathscr{LC}(\mathsf{P})}
({\PK},_-)$ defines an equivalence of categories
$\mathscr{LC}(\mathsf{P})\cong(\mathbf{C}^!)\defis\fmod$.
\end{corollary}

\subsection{The complex $\mathbb{P}^\bullet$}
\label{s3.2new}

We denote $\mathbb{P}^\bullet=\varinjlim{\PK}$. This should be
thought of playing the role of a minimal projective generator of
$\mathscr{LC}(\mathsf{P})$, see Corollary~\ref{corn543}. Proposition~\ref{pltcproj} gives us at least some information about the
structure of $\mathbf{C}$-direct summands of $\mathbb{P}^{\bullet}$. We would
like to describe the components $\mathbb{P}^{l}$ of $\mathbb{P}$ as well:

\begin{proposition}\label{pcomp}
Let $k\in\mZ$.
\begin{enumerate}[(a)]
\item \label{bba} ${\PK}$ is a complex of
$\mathbf{C}^\mZ-(\mathbf{C}^!)^{\geq k}$-bimodules, which is
projective both as a left and as a right module.
\item \label{bb}
$\mathbb{P}^{\bullet}$ is a complex of
$\mathbf{C}^\mZ\defis(\mathbf{C}^!)^\mZ$-bimodules, which is
projective both as a left and as a right module.
\end{enumerate}
\end{proposition}

\begin{proof}
Let $\mathbf{C}_0$ be the subcategory of $\mathbf{C}$ from
Subsection~\ref{s3.2} and  $l\in\mZ$. Then the
$\mathbf{C}^\mZ-(\mathbf{C}^!)^{\geq k}$-bimodule
structure on the component $\PKK^{l}$ is given by the following:
$\PKK^{l}\big((\la,i),(\mu,j)\big)=\mathcal{P}_{\mu}^l
\langle -j\rangle[j](\lambda,i)$ with the obvious assignments on
morphisms. We even claim that
\begin{displaymath}
{\PKK^{l}}\cong \bigoplus_{\la\in\Ob(\mathbf{C})}
\big(\mathbf{C}((\lambda,l),_-)
\otimes_{\mathbf{C_0}}{(\mathbf{C}^!)}^{\geq k}({}_-,(\lambda,l))\big)
\end{displaymath}
if $l\geq k$, which would imply the projectivity. For such $l$
and each $\lambda\in\Ob(\mathbf{C})$ we
can choose  $0\not=v_\la\in{\PKK^l}((\la,l),(\la,l))$ and
$$0\not=w_\la\in \big(\mathbf{C}((\lambda,l),(\la,l))\big)
\otimes_{\mathbf{C_0}}\big({(\mathbf{C}^!)}^{\geq k}((\lambda,l),(\la,l)\big).$$
Then, sending  $w_\la\mapsto v_\la$ (for all $\lambda$) defines a
homomorphism  of bimodules, which is surjective.
Since the bimodules have the same composition factors, the
surjection is an isomorphism.  This implies \eqref{bba} and \eqref{bb}
follows by  taking limits.
\end{proof}

\subsection{A homological description of the quadratic dual of a category}

Given a finite dimensional Koszul algebra $A$, it's Koszul dual is characterized or, depending on the author, even defined, as the Ext-algebra corresponding to the direct sum of all simple modules concentrated in degree zero. In this section we describe an extension of this characterization which applies to our more general setup.\\ 

Let $\Ext^{\mathrm{lin}}_{\mathbf{C}}(\mathsf{L})$ denote the
full subcategory of $\mathcal{D}(\mathbf{C}\defis\fgmod)$, objects of which are all complexes of the form
$\mathsf{L}(\lambda)^{\bullet}\langle -i\rangle[i]$,
$\lambda\in\Ob(\mathbf{C})$,
$i\in\Z$. Proposition~\ref{pltcproj} implies a
homological characterization of the category $\mathbf{C}^!$ as follows:

\begin{proposition}\label{corext}
There is an isomorphism of categories,
$$\Ext^{\mathrm{lin}}_{\mathbf{C}}(\mathsf{L})\cong
((\mathbf{C}^!)^{\Z})^{\mathrm{op}},$$
compatible with the natural $\mZ$-actions on both sides. In particular, $\Ext^{\mathrm{lin}}_{\mathbf{C}}(\mathsf{L})$ is generated by the elements of degree zero and one.
\end{proposition}

\begin{proof}
For each $\lambda\in\Ob(\mathbf{C})$ and $i\in\Z$ set
$\mathcal{Q}^{\bullet}_{\lambda,i}=\mathcal{Q}^{\bullet}_{\lambda}\langle
-i\rangle[i]$. Denote by $\mathscr{A}$ the
full subcategory of $\mathcal{D}(\mathbf{C}\defis\fgmod)$,
whose objects are all complexes
$\mathcal{Q}^{\bullet}_{\lambda,i}$, where
$\lambda\in\Ob(\mathbf{C})$ and $i\in\Z$. There is an obvious functor
\begin{eqnarray*}
  \alpha:\quad\Ext^{\mathrm{lin}}_{\mathbf{C}}(\mathsf{L})&\longrightarrow&
  \mathscr{A}\\
\mathsf{L}(\lambda)\langle
  -i\rangle[i]&\longmapsto&\mathcal{Q}^{\bullet}_{\lambda,i},
\end{eqnarray*}
identifying an object with its projective resolution.
On the other hand, by Lemma~\ref{l6}, we get the functor
\begin{eqnarray*}
\beta:=\mathrm{S}_{-1}\mathrm{Q}_0:\quad\mathscr{A}
\longrightarrow\mathscr{LC}(\mathsf{P}).
\end{eqnarray*}
By Proposition~\ref{pltcproj} we have $\beta\,\alpha\,\big(\mathsf{L}(\lambda)\langle
-i\rangle[i]\big)\cong\mathcal{I}_{\lambda}^{\bullet}$.
By Theorem~\ref{tqdual}, to prove the statement of the
proposition, it is therefore enough to show that functor
$\beta\,\alpha$ is full
and faithful. From the definition of $\alpha$ and $\mathscr{A}$ it is in fact
enough to show that $\beta$ is fully faithful. Let $\lambda,\mu\in \Ob(\mathbf{C})$ and
$i,j\in\Z$. Since $\mathbf{C}$ is positively graded, it is easy to see that we have
\begin{equation}\label{mneqeq}
{\mathcal{D}(\mathbf{C}\defis\fgmod)}
(\mathrm{Q}_{-1}\mathcal{Q}^{\bullet}_{\lambda,i},
\mathsf{L}(\mu)\langle -j\rangle[j])=0.
\end{equation}
Hence we get the following chain of isomorphisms
\begin{displaymath}
\begin{array}{rclr}
&&{\mathscr{LC}(\mathsf{P})}\,
\big(\beta\mathcal{Q}^{\bullet}_{\lambda,i},\beta\mathcal{Q}^{\bullet}_{\mu,j}\big)\\
&=&{\mathscr{LC}(\mathsf{P})}\,
\big(\mathrm{S}_{-1}\mathrm{Q}_0\mathcal{Q}^{\bullet}_{\lambda,i},
\mathrm{S}_{-1}\mathrm{Q}_0\mathcal{Q}^{\bullet}_{\mu,j}\big)\\
&\stackrel{\text{(Lemma~\ref{ladjdlc})}}{\cong}&
{\mathcal{K}^{\vee}_{\mathbf{C}}}\,
\big(\op{incl}\mathrm{S}_{-1}\mathrm{Q}_0\mathcal{Q}^{\bullet}_{\lambda,i},
\mathcal{Q}^{\bullet}_{\mu,j}\big)&\\
&\cong&
{\mathcal{D}(\mathbf{C}\defis\fgmod)}\,
\big(\op{incl}\mathrm{S}_{-1}\mathrm{Q}_0\mathcal{Q}^{\bullet}_{\lambda,i},
\mathsf{L}(\mu)\langle -j\rangle[j]\big)\\
&\stackrel{(\text{Proposition~\ref{pltcproj} and \eqref{mneqeq}})}{\cong}&
{\mathcal{D}(\mathbf{C}\defis\fgmod)}\,
\big(\mathcal{Q}^{\bullet}_{\lambda,i},
\mathsf{L}(\mu)\langle -j\rangle[j]\big)
\\
&\cong&
{\mathcal{D}(\mathbf{C}\defis\fgmod)}\,
\big(\mathsf{L}(\lambda)\langle -i\rangle[i],
\mathsf{L}(\mu)\langle -j\rangle[j]\big)\\
&\cong&
\mathscr{A}
\big(\mathsf{L}(\lambda)\langle -i\rangle[i],
\mathsf{L}(\mu)\langle -j\rangle[j]\big).
\end{array}
\end{displaymath}
Hence, the functor $\beta=\mathrm{S}_{-1}\mathrm{Q}_0$ is fully faithful. The claim follows.
\end{proof}

\begin{corollary}\label{clinextgen}
The morphism space of the positively graded category
$\Ext^{\mathrm{lin}}_{\mathbf{C}}(\mathsf{L})$
can be generated by elements of degree $0$ or $1$ only.
\end{corollary}

\begin{proof}
We have
$\Ext^{\mathrm{lin}}_{\mathbf{C}}(\mathsf{L})\cong
((\mathbf{C}^!)^{\Z})^{\mathrm{op}}$ by Proposition~\ref{corext},
and the positively graded category
$((\mathbf{C}^!)^{\Z})^{\mathrm{op}}$
is quadratic by definition, hence its morphism space is generated
in degrees $0$ and $1$.
\end{proof}

As an immediate consequence of Proposition~\ref{corext} we also
obtain the following statement, which is obvious for Koszul algebras:

\begin{corollary}\label{bred}
If $A$ is a positively graded algebra of finite homological
dimension. Then $A^{!}$ is finite-dimensional.
\end{corollary}

\begin{proof}
By Proposition~\ref{corext}, $A^{!}$ is a subalgebra of the
ext-algebra $\Ext_A^*(\mathsf{L},\mathsf{L})$, which is finite-dimensional, because $A$
is assumed to have finite homological dimension.
\end{proof}

\section{The quadratic duality functor}\label{section3}

The purpose of this section is to introduce what we call the
quadratic duality functor. In some sense it is a generalization
of the Koszul duality functor for positively graded Koszul algebras.
We will start with stating some general abstract nonsense. For
details we refer for example to \cite{KellerinKZ} and \cite{Deligne}.

Let $\mathscr{A}$ and $\mathscr{B}$ be two arbitrary $\Bbbk$-linear
categories and $\mathcal{X}^\bullet$ be a complex of
$\mathscr{A}\defis\mathscr{B}$-bimodules. Then we have the {\em inner
$\mathrm{Hom}$ functor}
\begin{displaymath}
\mathrm{Hom}_{\mathscr{A}}^{\bullet}(\mathcal{X}^\bullet,{}_-):
\mathcal{C}(\mathscr{A})\to \mathcal{C}(\mathscr{B}),
\end{displaymath}
as defined in \cite[III.6.14]{GM}. At the same time, for any
complex $\mathcal{Z}^\bullet$ of $\mathscr{B}$-modules we have
the associated bicomplex
$\mathcal{X}^\bullet\otimes_{\mathscr{B}}\mathcal{Z}^\bullet$. Applying the functor $\mathrm{Tot}$ of taking the total
complex defines a functor,
\begin{displaymath}
\mathcal{X}^\bullet\otimes_{\mathscr{B}} {}_-:
\mathcal{C}(\mathscr{B})\to \mathcal{C}(\mathscr{A}).
\end{displaymath}
These functors form an adjoint pair
$(\mathcal{X}^\bullet\otimes_{\mathscr{B}} {}_-,
\mathrm{Hom}_{\mathscr{A}}^{\bullet}(\mathcal{X}^\bullet,{}_-))$.

\subsection{Definition and the main theorem}\label{section3.n1}

Recall from Proposition~\ref{pcomp} that $\mathbb{P}^\bullet$ is a complex of
$\mathbf{C}^\Z\defis{\mathbf{C}^!}^\Z$-bimodules. Hence by the general 
definition above we have the following pair of functors:
\begin{eqnarray}
\label{KK-01}
\xymatrix{
\mathcal{C}(\mathbf{C}\defis\gMod)
\ar@/^/[rrrrrr]^{\mathrm{F}:=
\mathrm{E}_{\mathbf{C}^!}^{-1}\,{\mathrm{Hom}}_{\mathbf{C}^\Z}
\big(\mathbb{P}^\bullet,\mathrm{E}_{\mathbf{C}}{}(_-)\big)}
&&&&&&
\ar@/^/[llllll]^{\mathrm{F}':=
\mathrm{E}_{\mathbf{C}}^{-1}\,\mathbb{P}^\bullet
{\otimes}_{\mathbf{C!}}{}\mathrm{E}_{\mathbf{C}^!}(_-)}
\mathcal{C}(\mathbf{C}^!\defis\gMod)
}.
\end{eqnarray}

\begin{proposition}\label{tkfunctors-p}
Let $\mathbf{C}$ be a positively graded category.
\begin{enumerate}[(i)]
\item\label{tkfunctors-p.1}
The functors $\mathrm{F}$ and $\mathrm{F}'$ as in \eqref{KK-01} form 
a pair $(\mathrm{F}',\mathrm{F})$ of adjoint functors.
\item\label{tkfunctors-p.2}
For every $\mathcal{X}^{\bullet}\in\mathcal{C}(\mathbf{C}\defis\fgmod)$,
$\mathcal{Y}^{\bullet}\in \mathcal{C}(\mathbf{C}^!\defis\fgmod)$,
$i,j\in\Z$, we have
\begin{displaymath}
\begin{array}{lcl}
\mathrm{F}(\mathcal{X}^{\bullet}\langle j\rangle[i])& =&
(\mathrm{F}\mathcal{X}^{\bullet})\langle -j\rangle[i+j],\\
\mathrm{F}'(\mathcal{Y}^{\bullet}\langle j\rangle[i])&=&
(\mathrm{F}'\mathcal{Y}^{\bullet})\langle -j\rangle[i+j].
\end{array}
\end{displaymath}
\item\label{tkfunctors-p.3}
For $\lambda\in\Ob(\mathbf{C})$ and $i,j\in\Z$ we have
\begin{eqnarray}
\label{provedal-1}
\begin{array}{lcl}
\mathrm{F}(\mathsf{L}_{\mathbf{C}}(\lambda)^\bullet\langle j\rangle[i])&\cong&
\mathsf{I}_{\mathbf{C}^!}(\lambda)^\bullet\langle -j\rangle[i+j],\\
\mathrm{F}'(\mathsf{L}_{\mathbf{C}^!}
(\lambda)^\bullet\langle j\rangle[i])&\cong&
\mathsf{P}_{\mathbf{C}}(\lambda)^\bullet\langle -j\rangle[i+j].
\end{array}
\end{eqnarray}
\item\label{tkfunctors-p.4}
We have:
\begin{displaymath}
\begin{array}{lcl}
\mathrm{F}(\mathcal{C}^{\downarrow}(\mathbf{C}\defis\gMod))&\subset&
\mathcal{C}^{\uparrow}(\mathbf{C}^!\defis\gMod),\\
\mathrm{F}'(\mathcal{C}^{\uparrow}(\mathbf{C}^!\defis\gMod))&\subset&
\mathcal{C}^{\downarrow}(\mathbf{C}\defis\gMod).
\end{array}
\end{displaymath}
\item\label{tkfunctors-p.5}
$\mathrm{F}$ sends acyclic complexes from 
$\mathcal{C}^{\downarrow}(\mathbf{C}\defis\gMod)$
to acyclic complexes; and $\mathrm{F}'$ sends acyclic complex from 
$\mathcal{C}^{\uparrow}(\mathbf{C}^!\defis\gMod)$ to acyclic complexes.
\end{enumerate}
\end{proposition}

\begin{proof}
The statement \eqref{tkfunctors-p.1} follows from \cite[III.6.14]{GM}.
The properties $\mathrm{F}(\mathcal{X}^{\bullet}[i])=
(\mathrm{F}\mathcal{X}^{\bullet})[i]$, and $\mathrm{F}'(\mathcal{Y}^{\bullet}[i])=
(\mathrm{F}'\mathcal{Y}^{\bullet})[i]$, as well as
$\mathrm{F}(\mathcal{X}^{\bullet}\langle
j\rangle)=(\mathrm{K}\mathcal{X}^{\bullet})\langle -j\rangle[j]$, and
$\mathrm{F}'(\mathcal{Y}^{\bullet}\langle j\rangle)=
(\mathrm{F}'\mathcal{Y}^{\bullet})\langle -j\rangle[j]$
follow immediately
from the definitions of $\mathrm{F}$, $\mathrm{F}'$, and
$\mathbb{P}^{\bullet}$. This proves \eqref{tkfunctors-p.2}.

From Theorem~\ref{tqdual} and Proposition~\ref{pcatlcomp}\eqref{pcatlcomp.2} 
we have 
\begin{eqnarray}\label{eq:GL}
\mathrm{F}'(\mathsf{L}_{\mathbf{C}^!}(\lambda))=
\mathsf{P}_{\mathbf{C}}(\lambda).  
\end{eqnarray}
From the definition of $\mathbb{cP}^{\bullet}$ it follows immediately that
$\mathsf{M}^{\bullet}:=\mathrm{F}(\mathsf{L}_{\mathbf{C}}(\lambda)^{\bullet})$
 is concentrated in position $0$ (i.e. it is in fact a $\mathbf{C}^!$-module).
We have to show that $\mathsf{M}$ is an indecomposable
injective module. To see this we
calculate
\begin{eqnarray}
&&\mathbf{C}^!\defis\fgmod(\mathsf{L}_{\mathbf{C}^!}(\mu)\langle k\rangle,\mathsf{M}^0)
\nonumber\\
&\cong&{\mathcal{C}}(\mathbf{C}^!\defis\fgmod)\;
\big(\mathsf{L}_{\mathbf{C}^!}(\mu)^\bullet\langle k\rangle, \mathsf{M}^{\bullet}\big)\nonumber\\
&\cong&{\mathcal{C}}(\mathbf{C}^!\defis\fgmod)\;
\big(\mathsf{L}_{\mathbf{C}^!}(\mu)^\bullet\langle k\rangle,
\mathrm{F}(\mathsf{L}_{\mathbf{C}}(\lambda)^\bullet)\big)\nonumber\\
&\cong&{\mathcal{C}}(\mathbf{C}\defis\fgmod)\;
\big(\mathrm{F}'(\mathsf{L}_{\mathbf{C}^!}(\mu)^\bullet\langle k\rangle),
\mathsf{L}_{\mathbf{C}}(\lambda)^\bullet\big)\label{adjoint.1}\\
&\cong&{\mathcal{C}}(\mathbf{C}\defis\fgmod)\;
\big(\mathsf{P}_\mathbf{C}(\mu)^\bullet\langle -k\rangle[k],\mathsf{L}_{\mathbf{C}}(\lambda)^\bullet\big)\label{part2.1}\\
&\cong&\begin{cases}\Bbbk,&\lambda=\mu\text{ and }k=0;
\\0,&\text{otherwise.}\end{cases}\nonumber
\end{eqnarray}
Here, the isomorphism \eqref{adjoint.1} follows by adjointness and the isomorphism~\eqref{part2.1} is given by \eqref{eq:GL}. This implies that 
that module $\mathsf{M}^0$ has the simple socle 
$\mathsf{L}_{\mathbf{C}^!}(\lambda)$. Now to prove that 
$\mathsf{M}^0\cong \mathsf{I}_{\mathbf{C}^!}(\lambda)$ it remains to compare 
the characters: Let $\mu\in \mathbf{C}$, $i\geq 0$, and $m_{\mu,\lambda}^i$
denote the multiplicity of $\mathsf{P}(\lambda)$ as a direct summand of
the zero component of the complex 
$\mathcal{P}_{\mu}^{\bullet}\langle -i\rangle[i]$.
From the definition of $\mathrm{F}$ we have that the dimension of $\mathsf{M}^0(\mu)_{-i}$ equals $m_{\mu,\lambda}^i$, that is the
composition multiplicity of $\mathsf{L}_{\mathbf{C}^!}(\lambda)$ in 
$\mathsf{P}_{\mathbf{C}^!}(\mu)\langle i\rangle$. The latter equals the 
composition multiplicity of $\mathsf{L}_{\mathbf{C}^!}(\mu)\langle i\rangle$ in 
$\mathsf{I}_{\mathbf{C}^!}(\lambda)$ (as both numbers equal the dimension of
$\mathrm{Hom}_{\mathbf{C}^!}(\mathsf{P}_{\mathbf{C}^!}(\mu)\langle i\rangle,
\mathsf{I}_{\mathbf{C}^!}(\lambda))$).
Now the statement \eqref{tkfunctors-p.3} follows from \eqref{tkfunctors-p.2}.

The statement \eqref{tkfunctors-p.4} follows from 
\eqref{tkfunctors-p.3} and \eqref{tkfunctors-p.2} by a direct calculation.

Finally, to prove \eqref{tkfunctors-p.5} we first note the following
simplification: if 
$\mathcal{X}^{\bullet}\in\mathcal{C}^{\downarrow}(\mathbf{C}\defis\gMod)$,
$\lambda\in\mathrm{Ob}(\mathbf{C})$, and $i\in\mathbb{Z}$, then from the
construction of $\mathcal{P}_{\lambda}^{\bullet}$ it follows that the bicomplex
of vector spaces $\mathrm{Hom}_{\mathbf{C}^{\mathbb{Z}}}(\mathcal{P}_{\lambda}^{\bullet}
\langle -i\rangle[i],\mathcal{X}^{\bullet})$ has only finitely many non-zero
components, moreover, they all are finite-dimensional. Hence in the definition 
of the functor $\mathrm{F}$ all direct products which occur are finite
direct products of finite dimensional spaces. Hence they coincide with the
corresponding direct sums.

Let now
$\mathcal{X}^{\bullet}\in\mathcal{C}^{\downarrow}(\mathbf{C}\defis\gMod)$ 
be an acyclic complex of graded $\mathbf{C}^!$-modules. Then we can write $\mathcal{X}^{\bullet}$ as a direct sum of acyclic complexes of vector spaces of the form $0\rightarrow \mathsf{V}\rightarrow \mathsf{W}\rightarrow 0$, 
where $\mathsf{V}\cong \mathsf{W}\cong  \Bbbk\langle i\rangle$ is such that
both $\mathsf{V}$ and $\mathsf{W}$ are annihilated by all but one
$e_{\lambda}$. Denote by $\mathsf{V}'$ and $\mathsf{W}'$ the subspace of 
$\mathrm{F}(\mathcal{X}^{\bullet})$, which consist of all those homomorphisms
in which the images of the generators of indecomposable projective summands 
of $\mathbb{P}^{\bullet}$ belong to $\mathsf{V}$ and $\mathsf{W}$ respectively.
From the definitions it is obvious that both $\mathsf{V}'$ and $\mathsf{W}'$
are in fact $\mathbf{C}^!$-modules. From \eqref{tkfunctors-p.3} it even follows
that both  $\mathsf{V}'$ and $\mathsf{W}'$ are indecomposable injective
$\mathbf{C}^!$-modules. The differential in $\mathcal{X}^{\bullet}$ induces an
isomorphism $\mathsf{V}'\cong \mathsf{W}'$ of these $\mathbf{C}^!$-modules.
This means that the complex $\mathrm{F}(\mathcal{X}^{\bullet})$ decomposes into
a direct sum of acyclic complexes of the form $0\rightarrow \mathsf{V}'\cong
\mathsf{W}'\rightarrow 0$, and hence is acyclic. For the functor $\mathrm{F}'$ 
the proof is similar (and even easier as we do not have any direct product in 
the definition at all). This completes the proof.
\end{proof}

Thanks to Proposition~\ref{tkfunctors-p} we can restrict the functors 
$\mathrm{F}$ and $\mathrm{F}'$ to subcategories 
$\mathcal{C}^{\downarrow}(\mathbf{C}\defis\gMod)$ and
$\mathcal{C}^{\uparrow}(\mathbf{C}^!\defis\gMod)$ respectively, and then 
derive the picture  \eqref{KK-01} in the following way:
\begin{eqnarray}
\label{KK}
\xymatrix{
\mathcal{D}^{\downarrow}(\mathbf{C}\defis\gMod)
\ar@/^/[rrrrrr]^{\mathrm{K}_{\mathbf{C}}:=
\mathrm{E}_{\mathbf{C}^!}^{-1}\,{\rhom}_{\mathbf{C}^\Z}
\big(\mathbb{P}^\bullet,\mathrm{E}_{\mathbf{C}}{}(_-)\big)}
&&&&&&
\ar@/^/[llllll]^{\mathrm{K}'_{\mathbf{C}}
:=\mathrm{E}_{\mathbf{C}}^{-1}\,\mathbb{P}^\bullet\stackrel{\mathcal L}{\otimes}_{\mathbf{C!}}{}\mathrm{E}_{\mathbf{C}^!}(_-)}
\mathcal{D}^{\uparrow}(\mathbf{C}^!\defis\gMod)
}.
\end{eqnarray}
The functors $\mathrm{K}=\mathrm{K}_{\mathbf{C}}$ and $\mathrm{K}'=\mathrm{K}'_{\mathbf{C}}$ are what
we call {\it quadratic duality functors}. We emphasize once more that we have
\begin{displaymath}
\mathrm{K}=\mathcal{R}(\mathrm{F}|_{\mathcal{C}^{\downarrow}
(\mathbf{C}\defis\gMod)}) \quad\text{ and }\quad
\mathrm{K}'=\mathcal{L}(\mathrm{F}'|_{\mathcal{C}^{\uparrow}
(\mathbf{C}^!\defis\gMod)})
\end{displaymath}

The following alternative description depicts clearly the importance
of the equivalence $\epsilon$ from Theorem~\ref{tqdual} (namely,
$\mathrm{K}'$ is just the equivalence $\epsilon^{-1}$ extended to
the derived category followed by taking the total complex):

\begin{proposition}
\label{funnyd2}
Up to isomorphism of functors, the following diagram commutes:
\begin{equation*}
\xymatrix{
&&\mathcal{D}^{\uparrow}(\mathscr{LC}(\mathsf{P}))
\ar[lld]_{\mathrm{Tot}}
&&\\
\mathcal{D}^{\downarrow}(\mathbf{C}\defis\fgmod)
&&&&
\ar[llll]^{\mathrm{K}'\mathrm{E}_{\mathbf{C}^!}^{-1}}
\ar[llu]_{
\quad\mathrm{E}_{\mathbf{C}}^{-1}\,\mathbb{P}^{\bullet}
\otimes_{(\mathbf{C}^!)^{\Z}}
{}_-}
\mathcal{D}^{\uparrow}((\mathbf{C}^!)^\mZ\defis\fmod).
}
\end{equation*}
\end{proposition}

\begin{proof}
This follows directly from the definitions and
Proposition~\ref{pcomp}\eqref{bb}.
\end{proof}

Our main statement here is the following:

\begin{theorem}[Quadratic duality]\label{tkfunctors}
Let $\mathbf{C}$ be a (positively graded) category.
\begin{enumerate}[(i)]
\item\label{tkfunctors.1}
$(\mathrm{K}',\mathrm{K})$ is a pair of adjoint functors.
\item\label{tkfunctors.2}
For every $\mathcal{X}^{\bullet}\in
\mathcal{D}^{\downarrow}(\mathbf{C}\defis\fgmod)$,
$\mathcal{Y}^{\bullet}\in
\mathcal{D}^{\uparrow}(\mathbf{C}^!\defis\fgmod)$,
$i,j\in\Z$, we have
\begin{displaymath}
\begin{array}{lcl}
\mathrm{K}(\mathcal{X}^{\bullet}\langle j\rangle[i])& =&
(\mathrm{K}\mathcal{X}^{\bullet})\langle -j\rangle[i+j],\\
\mathrm{K}'(\mathcal{Y}^{\bullet}\langle j\rangle[i])&=&
(\mathrm{K}'\mathcal{Y}^{\bullet})\langle -j\rangle[i+j].
\end{array}
\end{displaymath}
\item\label{tkfunctors.3}
For $\lambda\in\Ob(\mathbf{C})$ and $i,j\in\Z$ we have
\begin{eqnarray}
\label{provedal}
\begin{array}{lcl}
\mathrm{K}(\mathsf{L}_{\mathbf{C}}(\lambda)^\bullet\langle j\rangle[i])&\cong&
\mathsf{I}_{\mathbf{C}^!}(\lambda)^\bullet\langle -j\rangle[i+j],\\
\mathrm{K}'(\mathsf{L}_{\mathbf{C}^!}
(\lambda)^\bullet\langle j\rangle[i])&\cong&
\mathsf{P}_{\mathbf{C}}(\lambda)^\bullet\langle -j\rangle[i+j].
\end{array}
\end{eqnarray}
\end{enumerate}
\end{theorem}

\begin{proof}
The statement \eqref{tkfunctors.1} follows from general nonsense (see e.g 
\cite[8.1.4]{KellerinKZ} or \cite{Deligne}). The rest follows follows from the
definitions  and  Proposition~\ref{tkfunctors-p}.
\end{proof}

\subsection{The duality functors applied to modules}

Calling $\mathrm{K}$ a {\em duality} functor might be too
optimistic, in particular, since
$\mathrm{K}$ is not an equivalence in general (see e.g. \cite[Proposition
8.1.4]{KellerinKZ} and also Theorem~\ref{koszul}). However, later on we will see many ``duality-like'' effects in our situation,
which, from our point of view, justify this usage. Our first
general observation is the following:

\begin{proposition}\label{prcat508}
\begin{enumerate}[(i)]
\item\label{prcat508.1} Let
$\mathsf{X}\in\mathbf{C}\defis\fgmod\cap
\mathcal{D}^{\downarrow}(\mathbf{C}\defis\fgmod)$. Then
$$\mathrm{K}\mathsf{X}^{\bullet}\in
\mathscr{LC}(\mathsf{I}_{\mathbf{C}^!})\cap
\mathcal{D}^{\uparrow}(\mathbf{C}^!\defis\fgmod).$$
\item\label{prcat508.2} The functor $\mathrm{Tot}$ from
Proposition~\ref{funnyd2} sends non-zero objects to non-zero objects.
\end{enumerate}
\end{proposition}

\begin{proof}
Let $\mathsf{X}=\mathsf{X}_0\supset\mathsf{X}_1\supset\dots$
be a decreasing filtration of $\mathsf{X}$ such that
for every $i=0,1,\dots$ the module
$\mathsf{X}_i/\mathsf{X}_{i+1}$ is semi-simple and
concentrated in a single degree, say $k_i$. Since
$\mathbb{P}^{\bullet}$ is a complex of projective
$\mathbf{C}$-modules, analogously to the proof of
Theorem~\ref{tkfunctors}\eqref{tkfunctors.3}, the module
$\mathsf{X}_i/\mathsf{X}_{i+1}$ gives rise to an
injective $\mathbf{C}^!$-module, which is, however,
shifted by $\langle -k_i\rangle[k_i]$ because of
Theorem~\ref{tkfunctors}\eqref{tkfunctors.2}. The claim
\eqref{prcat508.1} follows.

Every element from $\mathcal{D}^{\uparrow}
(\mathscr{LC}(\mathsf{P}))$ is a double complex of
projective $\mathbf{C}$-mo\-du\-les, linear in one direction,
and  isomorphic to a direct sum of trivial complexes and
complexes of the form
\begin{equation}\label{eq875}
\dots\to 0\to \mathsf{M}\to 0\to \dots
\end{equation}
in the other direction. Moreover, it is acyclic in
$\mathcal{D}^{\uparrow}(\mathscr{LC}(\mathsf{P}))$
if and only if no direct summands of
the form \eqref{eq875} in the second direction occur.
Since the image of $\mathrm{Tot}$ is a complex of
projective $\mathbf{C}$-module, bounded from the
right, we obtain that the image is acyclic if and
only if the bicomplex we started with was acyclic.
The claim \eqref{prcat508.2} follows and the proof
of Proposition~\ref{prcat508} is complete.
\end{proof}

In case $\mathbf{C}$ is a quadratic category, the functors $\mathrm{K}$ and
$\mathrm{K}'$ are particularly well-behaved as we will illustrate now.
We first show that the functor $\mathrm{K}'$ for $\mathbf{C}$ can be realized
using the functor $\mathrm{K}$ for $(\mathbf{C}^!)^{\mathrm{op}}$
(which means that these two functors are in fact {\em dual} to each other).

\begin{proposition}\label{qdualq}
Assume that $\mathbf{C}$ is quadratic. Then
\begin{displaymath}
\mathrm{K}_{(\mathbf{C}^!)^{\mathrm{op}}}\cong
\mathbb{D}\,\mathrm{K}'_{\mathbf{C}}\,\mathbb{D}.
\end{displaymath}
\end{proposition}

\begin{proof}
First let $\mathsf{M}\in\mathbf{C}^!\defis\fgmod$ be such that 
$\mathsf{M}_i=0$ for all big enough $i$. Let further
$\mathsf{X}=\mathbf{C}\otimes_{\mathbf{C}_0}\mathbf{C}^!$. Then 
$\mathsf{X}$ is a graded projective ${\mathbf{C}}-\mathbf{C}^!$-bimodule,
which has a unique decomposition into a direct sum of 
indecomposable projective ${\mathbf{C}}-\mathbf{C}^!$-bimodules of the 
form $\mathsf{N}_{\lambda}=\mathbf{C}(\lambda,{}_-)\otimes_{\mathbf{C}_0}
\mathbf{C}^!({}_-,\lambda)$, where $\lambda\in \mathbf{C}$, each occurring 
with multiplicity one (this multiplicity is given by the dimension of the
homomorphism space to the appropriate simple bimodule).
Under these assumptions the graded left $\mathbf{C}$-module 
$\mathsf{X}\otimes_{\mathbf{C}^!}\mathbb{D}\mathsf{M}$ has 
finite-dimensional graded components (and we even have 
$(\mathsf{X}\otimes_{\mathbf{C}^!}\mathbb{D}\mathsf{M})_i=0$ for all 
small enough $i$). Then the module 
$\mathbb{D}(\mathsf{X}\otimes_{\mathbf{C}^!}\mathbb{D}\mathsf{M})$ is a
well-defined right ${\mathbf{C}}$-module with finite-dimensional graded
components, moreover,  this module is isomorphic to the module
$\HOM_\Bbbk(\mathsf{X}\otimes_{\mathbf{C}^!}\mathbb{D}\mathsf{M},\Bbbk)$
as graded ${\mathbf{C}}$-module by definition of $\mathbb{D}$ 
(note that we understand the Hom-functor 
using the definitions from the last part of  Subsection~\ref{s2.2}).
From the definition of $\mathsf{N}_{\lambda}$ and Section~\ref{s3.1} 
it follows that the graded components of the graded vector spaces in 
the formula \eqref{eq:new} below are finite dimensional, and hence the 
adjunction morphism defines an isomorphism of these components:
\begin{eqnarray}
\label{eq:new}
\Phi:\HOM_\Bbbk(\mathsf{N}_{\lambda}\otimes_{\mathbf{C}^!}
\mathbb{D}\mathsf{M},\Bbbk)& \cong& \HOM_{\mathbf{C}^!}(\mathbb{D}\mathsf{M},
\HOM_\Bbbk(\mathsf{N}_{\lambda},\Bbbk)). 
\end{eqnarray}
The space on the right hand side of \eqref{eq:new} is nothing else 
than 
$\HOM_{\mathbf{C}^!}(\mathbb{D}\mathsf{M},\mathbb{D}\mathsf{N}_{\lambda})$ 
which is naturally isomorphic to $\HOM_{({\mathbf{C}^!})^{\rm{op}}}(\mathsf{N}_{\lambda},\mathsf{M})$ by applying the duality $\mathbb{D}$. Now, from the definition of 
$\mathsf{X}$ and the assumptions on $\mathsf{M}$ it follows that there is
a natural isomorphism of graded right ${\mathbf{C}}$-modules with finite-dimensional graded components:
\begin{eqnarray*}
\Phi: \mathbb{D}(\mathsf{X}\otimes_{\mathbf{C}^!}\mathbb{D}\mathsf{M})&\cong
&\HOM_{({\mathbf{C}^!})^{\rm{op}}}(\mathsf{X},\mathsf{M}).
\end{eqnarray*}
Let now $\mathsf{M}^\bullet\in \mathcal{C}^{\uparrow}
(\mathbf{C}^!\defis\fgmod)$ and $\mathsf{X}^\bullet=\mathbb{P}^{\bullet}$.
Then, by Proposition~\ref{pcomp}, the components of $\mathsf{X}^\bullet$ 
are (up to shift) isomorphic to the bimodule $\mathsf{X}$ above. 
Hence it follows that  $\Phi$ induces a natural isomorphism
\begin{eqnarray*}
\Phi^{i,j}:\mathbb{D}(\mathsf{X}^i\otimes_{\mathbf{C}^!}
(\mathbb{D}\mathsf{M})^j)&\cong&
\HOM_{({\mathbf{C}^!})^{\rm{op}}}(\mathsf{X}^i,\mathsf{M}^j) 
\end{eqnarray*}
for any $i,j\in\mZ$. The naturality of $\Phi$ induces an isomorphism of bicomplexes  
\begin{eqnarray*}   \mathbb{D}(\mathsf{X}^\bullet\otimes_{\mathbf{C}^!}\mathbb{D}
\mathsf{M}^\bullet)&\cong&
\HOM_{({\mathbf{C}^!})^{\rm{op}}}(\mathsf{X}^\bullet,\mathsf{M}^\bullet).
\end{eqnarray*}
By the arguments from the proof of Proposition~\ref{tkfunctors-p}, taking
the total complex reduces to taking direct sums of finitely many 
non-zero spaces. Hence the above induces an isomorphism of the 
corresponding total complexes. The 
isomorphism $\mathbb{D}\mathrm{K}'_{\mathbf{C}}\mathbb{D}\cong
\mathrm{K}_{(\mathbf{C}^!)^{\mathrm{op}}}$ follows therefore from the definition of the involved functors. 
\end{proof}

For $\lambda\in\Ob(\mathbf{C})$ let $\mathcal{Q}_{\lambda^!}^{\bullet}$
denote a minimal  projective resolution
of $\mathsf{L}_{(\mathbf{C}^!)^{\mathrm{op}}}(\lambda)
\in(\mathbf{C}^!)^{\mathrm{op}}\defis\fgmod$.
The following result says that the images of indecomposable projective
(resp. injective) modules under the functor $\mathrm{K}$ (resp.
$\mathrm{K}'$) is nothing else than the linear part of a minimal
injective (projective) resolution of the corresponding simple module.

\begin{proposition}\label{primpi}
Let $\mathbf{C}$ be a positively graded category. Then there are
isomorphisms
\begin{enumerate}[(i)]
\item\label{primpi.2} $\mathrm{K}'
\mathsf{I}_{\mathbf{C}^!}(\lambda)^\bullet
\cong \mathcal{I}_{\lambda}^{\bullet}\cong
\mathrm{S}_{-1}\mathrm{Q}_{0}^{}
\mathcal{Q}_{\lambda}^{\bullet}$ of objects
in
$\mathscr{LC}(\mathsf{P}_{\mathbf{C}})\cap
\mathcal{D}^{\downarrow}(\mathbf{C}\defis\fgmod)$,
and
\item\label{primpi.1}
$\mathrm{K}\mathsf{P}_{\mathbf{C}}(\lambda)^\bullet\cong
\mathbb{D}\,\mathrm{S}_{-1}^{(\mathbf{C}^!)^{\mathrm{op}}}
\mathrm{Q}_{0}^{(\mathbf{C}^!)^{\mathrm{op}}}
\mathcal{Q}_{\lambda^!}^{\bullet}$
of objects in
$\mathscr{LC}(\mathsf{I}_{\mathbf{C}^!})\cap
\mathcal{D}^{\uparrow}(\mathbf{C}^!\defis\fgmod)$
in case $\mathbf{C}$ is quadratic.
\end{enumerate}
\end{proposition}

\begin{proof}
Let $\la\in\Ob(\mathbf{C})$. From Proposition~\ref{funnyd2},
Theorem~\ref{tqdual} and Proposition~\ref{pltcproj}\eqref{concrete}
we know that $\mathrm{K}'\mathsf{I}_{\mathbf{C}^!}(\lambda)\cong
\mathcal{I}_{\lambda}^{\bullet}\cong
\mathrm{S}_{-1}\mathrm{Q}_{0}\mathcal{Q}_{\lambda}^{\bullet}$. This proves \eqref{primpi.2}. From this \eqref{primpi.1} follows using
Proposition~\ref{qdualq}.
\end{proof}

For quadratic  $\mathbf{C}$ we have
$((\mathbf{C}^!)^{\mathrm{op}})^!=((\mathbf{C}^!)^!)^{\mathrm{op}}=
\mathbf{C}^{\mathrm{op}}$ canonically and from Theorem~\ref{tqdual} and
Proposition~\ref{pcatlcomp} it follows immediately  that the
categories $\mathbf{C}\defis\fgmod\cap
\mathcal{D}^{\downarrow}(\mathbf{C}\defis\fgmod)$ and
$\mathscr{LC}(\mathsf{I}_{\mathbf{C}^!})\cap
\mathcal{D}^{\uparrow}(\mathbf{C}^!\defis\fgmod)$ are equivalent.
This equivalence can also be realized in the following way:

\begin{proposition}\label{onemore}
Assume that $\mathbf{C}$ is quadratic. Then
\begin{displaymath}
\mathrm{K}:\mathbf{C}\defis\fgmod\cap
\mathcal{D}^{\downarrow}(\mathbf{C}\defis\fgmod)\longrightarrow
\mathscr{LC}(\mathsf{I}_{\mathbf{C}^!})\cap
\mathcal{D}^{\uparrow}(\mathbf{C}^!\defis\fgmod)
\end{displaymath}
is an equivalence.
\end{proposition}

\begin{proof}
From the arguments in the
proof of Proposition~\ref{prcat508} it follows
immediately that the functor
\begin{displaymath}
\mathrm{K}:\mathbf{C}\defis\fgmod\cap
\mathcal{D}^{\downarrow}(\mathbf{C}\defis\fgmod)\to
\mathscr{LC}(\mathsf{I}_{\mathbf{C}^!})\cap
\mathcal{D}^{\uparrow}(\mathbf{C}^!\defis\fgmod)
\end{displaymath}
is exact. By Proposition~\ref{primpi}\eqref{primpi.1}
and Proposition~\ref{pltcproj}\eqref{concrete},
$\mathrm{K}$ sends indecomposable projective objects from
$\mathbf{C}\defis\fgmod$ to the corresponding
indecomposable projective objects from
$\mathscr{LC}(\mathsf{I}_{\mathbf{C}^!})$. By
\cite[Lemma~6]{MO}, the induced map on the morphisms is an isomorphism when
restricted to the part of degree $1$. Hence
it is an isomorphism, since $\mathbf{C}$ is quadratic.
This completes the proof.
\end{proof}

\subsection{Quadratic dual functors}\label{s3.5}
Let $\mathbf{C}$ be a positively graded category and
$\Lambda\subset\Ob(\mathbf{C})$, $\Lambda\neq \varnothing$. We denote by
$\mathbf{C}_{\Lambda}$ the full subcategory of $\mathbf{C}$ such that
$\Ob(\mathbf{C}_{\Lambda})=\Lambda$. The category $\mathbf{C}_{\Lambda}$
obviously inherits a positive grading. Let $\mathbf{B}_{\Lambda}$ be the
$\mathbf{C}^\mZ\defis(\mathbf{C}_{\Lambda})^\mZ$-bimodule
$\mathbf{C}(_-,_-)$, which means that it maps the object $\big((\la,i),(\mu,j)\big)$ to $\mathbf{C}(\la,\mu)_{j-i}$,
where $\mu\in\Ob(\mathbf{C})$, $\la\in\Ob(\mathbf{C}_\Lambda)$.
Further, we define the category
${}_{\Lambda}\mathbf{C}$ as follows: $\Ob({}_{\Lambda}\mathbf{C})=\Lambda$,
and for $\lambda,\mu\in\Ob({}_{\Lambda}\mathbf{C})$,
the space ${}_{\Lambda}\mathbf{C}(\lambda,\mu)$
is the quotient of $\mathbf{C}(\lambda,\mu)$ modulo the subspace,
generated by all morphisms, which factor through some object
outside $\Lambda$. Let ${}_{\Lambda}\mathbf{D}$ be the
$(\mathbf{C}^!)^\mZ\defis ({}_{\Lambda}\mathbf{C}^!)^\mZ$-bimodule
such that ${}_{\Lambda}\mathbf{D}((\mu,i),(\la,j))=
{}_\Lambda(\mathbf{C}^!)(\mu,\la)_{j-i}$ if $\mu\in\Lambda$ and which is
the trivial vector space otherwise. The assignments for the maps are
the obvious ones. The following  observation (which was made in
\cite[3.1]{Ma1}) about the connections between
$\mathbf{C}_{\Lambda}$ and ${}_{\Lambda}(\mathbf{C}^!)$ is easy but crucial:

\begin{lemma}
\label{technical}
There is an isomorphism of categories
\begin{eqnarray*}
 \tau:(\mathbf{C}_\Lambda)^!\defis\fgmod&\cong&
{}_\Lambda(\mathbf{C}^!)\defis\fgmod,
\end{eqnarray*}
such that $\tau\mathsf{M}(\la)=\mathsf{M}(\la)$ for any $\la\in\Lambda$
  and $\tau\mathsf{M}(f)=\mathsf{M}(f)$ for any morphism $f$
  homogeneous of degree one.
\end{lemma}

\begin{proof}
Let $\la,\mu\in\Lambda$. The equalities
\begin{multline*}
((\mathbf{C}_\Lambda)^!)_1(\la,\mu)=
{\bf d}((\mathbf{C}_\Lambda)_1(\mu,\la))=
{\bf d}(\mathbf{C}_1(\mu,\la))=\\=
{\bf d}(({}_\Lambda\mathbf{C})_1(\mu,\la))=
(({}_\Lambda\mathbf{C})^!)_1(\la,\mu)
\end{multline*}
give rise to an identification of the morphisms of degree one.
It is easy to see that this gives rise to an identification
$((\mathbf{C}_\Lambda)^!)(\la,\mu)=(({}_\Lambda\mathbf{C})^!)(\la,\mu)$
and the statement follows.
\end{proof}

Motivated by the Koszul duality (as proved in \cite{Steen}) between
translation functors and Zuckerman functors for the classical
Bernstein-Gelfand-Gelfand category $\mathcal{O}$, we would like to
extend the above  equivalence to the following correspondence on
functors:

\begin{theorem}\label{tqdfunctors}
Let $\mathbf{C}$ be a positively graded category and
$\Lambda\subset \Ob(\mathbf{C})$, $\Lambda\neq \varnothing$.
Then the following diagrams commute up to isomorphism of functors:
\small
\begin{equation*}\label{eqdiag1}
\xymatrix{
\mathcal{D}^{\downarrow}(\mathbf{C}\defis\fgmod)
\ar[r]^{\mathrm K_{\mathbf C}}
\ar[dd]_{\mathrm{F}'}
&
\mathcal{D}^{\uparrow}(\mathbf{C}^!\defis\fgmod)
\ar[dd]^{\mathrm{G}'}
&\mathcal{D}^{\downarrow}(\mathbf{C}\defis\fgmod)
&
\mathcal{D}^{\uparrow}(\mathbf{C}^!\defis\fgmod)
\ar[l]_{\mathrm{K}'_{\mathbf{C}}}
\\
\\
\mathcal{D}^{\downarrow}(\mathbf{C}_{\Lambda}\defis\fgmod)
\ar[r]^{\tau\mathrm{K}_{\mathbf{C}_\Lambda}}
&
\mathcal{D}^{\uparrow}({}_{\Lambda}(\mathbf{C}^!)\defis\fgmod)
&
\mathcal{D}^{\downarrow}(\mathbf{C}_{\Lambda}\defis\fgmod)
\ar[uu]^{\mathrm{F}}
&
\mathcal{D}^{\uparrow}({}_{\Lambda}(\mathbf{C}^!)\defis\fgmod),
\ar[uu]_{\mathrm{G}}
\ar[l]_{\mathrm{K}'_{\mathbf{C}_\Lambda}\tau^{-1}}
}
\end{equation*}
where
\begin{eqnarray*}
\mathrm{F}=\mathrm{E}_{\mathbf{C}}^{-1}\,\mathbf{B}_{\Lambda}
\stackrel{\mathcal L}{\otimes}_{(\mathbf{C}_\Lambda)^\mZ}
\mathrm{E}_{\mathbf{C}_\Lambda}(-),&
\mathrm{G}=\mathrm{E}_{\mathbf{C}^!}^{-1}\,{}_{\Lambda}\mathbf{D}
\stackrel{\mathcal L}{\otimes}_{({}_\Lambda(\mathbf{C}^!))^\mZ}
\mathrm{E}_{_\Lambda(\mathbf{C}^!)}(-),\\
\mathrm{F}'=\mathrm{E}_{\mathbf{C}_\Lambda}^{-1}{\rhom}_{\mathbf{C}^\mZ}
\big(\mathbf{B}_{\Lambda},\mathrm{E}_{\mathbf{C}}(-)\big),&
\mathrm{G}'=\mathrm{E}_{{}_\Lambda(\mathbf{C}^!)}^{-1}
{\rhom}_{(\mathbf{C}^!)^\mZ}({}_\Lambda\mathbf{D},
\mathrm{E}_{\mathbf{C}^!}(-)).
\end{eqnarray*}
\end{theorem}
\normalsize
\begin{proof}
Since the diagrams are adjoint to each other, it is enough to prove the
commutativity of the second say. We have the natural restriction functor
$\op{res}:\mathbf{C}\defis\fgmod\rightarrow\mathbf{C}_\Lambda\defis\fgmod$
which is isomorphic to $\mathbf{C}\defis\fgmod(\mathbf{B}_\Lambda,_-)$ and
has the right adjoint $\op{ind}$ given by tensoring with $\mathbf{B}_\Lambda$.
On the other hand we have the natural functor
$\op{J}:{}_\Lambda(\mathbf{C}^!)\defis\fgmod
\rightarrow\mathbf{C}^!\defis\fgmod$, which is given by
$\op{J}(\mathsf{M})(\la)=\mathsf{M}(\la)$
if $\la\in\Lambda$ and $\op{J}(\mathsf{M})(\la)=\{0\}$
otherwise, and on morphisms $\op{J}(\mathsf{M})(f)=\mathsf{M}(f)$
if $\mathsf{M}(f)$ is defined and $\op{J}(\mathsf{M})(f)=0$ otherwise.
In other words: $\op{J}\cong {{}_{\Lambda}\mathbf{D}}
\otimes_{{}_\Lambda(\mathbf{C}^!)^\mZ}{}_-$. We claim that there is an
isomorphism of functors as follows:
\begin{eqnarray*}
\epsilon_{\mathbf{C}}^{-1}\op{J}\cong
\op{ind}\epsilon_{\mathbf{C}_\Lambda}^{-1}\,\tau^{-1}.
\end{eqnarray*}
This can be checked by an easy direct calculation. The commutativity
of the second diagram above follows then directly from Theorem~\ref{tqdual}
and Proposition~\ref{funnyd2}.
\end{proof}

\begin{remark}
{\rm
The statement of Theorem~\ref{tqdfunctors} resembles the
equivalence of categories, given by Auslander's approximation
functor from \cite[Section~5]{Au}. A substantial part of the
``easy direct calculation'' in the proof repeats the calculation,
used to establish the fact that Auslander's functor is an equivalences
of certain categories.
}
\end{remark}

\subsection{The Koszul duality theorem}\label{s3.6}

We call a positively graded category $\mathbf{C}$ {\em Koszul}
provided that the minimal projective resolution of
$\mathsf{L}(\lambda)\in\mathbf{C}\defis\gmod$ is linear for every
$\lambda$. This generalizes the usual definition of Koszul algebras
(see e.g. \cite[Section~2]{BGS}). It is of course not a big surprise
that for Koszul categories all our previous results can be seriously
strengthened. Our main result here is the following:

\begin{theorem}[Koszul duality]\label{koszul}
Let $\mathbf{C}$ be a positively graded category. The following conditions are equivalent:
\begin{enumerate}[(a)]
\item\label{koszul.1} $\mathbf{C}$ is Koszul.
\item\label{koszul.2}
\begin{displaymath}
\xymatrix{
\mathcal{D}^{\downarrow}(\mathbf{C}\defis\fgmod)
\ar@/^/[rrrr]^{\mathrm{K}_\mathbf{C}}
&&&&
\mathcal{D}^{\uparrow}(\mathbf{C}^!\defis\fgmod)
\ar@/^/[llll]^{\mathrm{K}'_{\mathbf{C}^!}}
}
\end{displaymath}
are mutually inverse equivalences of categories.
\item\label{koszul.3} $\mathrm{K}\mathsf{P}_{\mathbf{C}}(\lambda)^\bullet
\cong\mathsf{L}_{\mathbf{C}^!}(\lambda)^\bullet$ for every
$\lambda\in\Ob(\mathbf{C})$.
\item\label{koszul.4} $\mathrm{K}'\mathsf{I}_{\mathbf{C}^!}(\lambda)^\bullet
\cong\mathsf{L}_{\mathbf{C}}(\lambda)^\bullet$ for every
$\lambda\in\Ob(\mathbf{C})$.
\item\label{koszul.5} The functor $\mathrm{Tot}$ from
\eqref{funnyd2} is dense.
\end{enumerate}
\end{theorem}

\begin{proof}
$\eqref{koszul.1}\Rightarrow\eqref{koszul.2}$. We assume that $\mathbf{C}$ is
Koszul. Since we have an adjoint pair of functors
$({\mathrm{K}'_{\mathbf{C}^!}}, \mathrm{K}_\mathbf{C})$ (see
Theorem~\ref{tkfunctors}) it is enough to show that the adjunction morphisms
are isomorphisms. We even claim that it is enough to show that the
adjunction morphisms are isomorphisms for any simple object. Indeed,
by the definition of $\mathcal{D}^{\downarrow}(\mathbf{C}\defis\fgmod)$,
for any $\mathcal{X}^{\bullet}\in
\mathcal{D}^{\downarrow}(\mathbf{C}\defis\fgmod)$,
$\lambda\in \Ob(\mathbf{C})$ and $i\in\mathbb{Z}$,  the bicomplex
\begin{displaymath}
\mathrm{Hom}_{\mathbf{C}-\fgmod}
(\mathcal{P}_{\lambda}^{\bullet}\langle -i\rangle[i],\mathcal{X}^{\bullet})
\end{displaymath}
has only finitely many non-zero components, each of which is a
finite-dimensional vector space. Hence the claim that the adjunction
morphism is an isomorphism for $\mathcal{X}^{\bullet}$ follows
from the corresponding statement for simple objects by taking the limit as in Corollary~\ref{Pk}.

Now let us prove the statement for simple objects.
We can of course assume that these simple objects are concentrated
in position zero. From Proposition~\ref{primpi} we have $\mathrm{K}\mathsf{L}_{\mathbf{C}}(\lambda)=
\mathsf{I}_{\mathbf{C}^!}(\lambda)^\bullet$.
From Proposition~\ref{primpi} we have an
isomorphism $\mathrm{K}'\mathsf{I}_{\mathbf{C}^!}(\lambda)\cong
\mathrm{S}_{-1}\mathrm{Q}_{0}\mathcal{Q}_{\lambda}^{\bullet}$. The latter one
is isomorphic to $\mathcal{Q}_{\lambda}^{\bullet}$, since $\mathbf{C}$ is
Koszul. Hence
$\mathrm{K}'\mathrm{K}\mathsf{L}_{\mathbf{C}}(\lambda)\cong
\mathsf{L}_{\mathbf{C}}(\lambda)$. Since the adjunction morphism
$\mathrm{K}'\mathrm{K}\mathsf{L}_{\mathbf{C}}(\lambda)\rightarrow
\mathsf{L}_{\mathbf{C}}(\lambda)$ is non-zero, it must be an isomorphism.
Via the duality $\mathbb{D}$ we could also say that $\mathbf{C}$ is
{\em Koszul} provided that the minimal injective resolution of
$\mathsf{L}(\lambda)\in\mathbf{C}\defis\gmod$ is linear for every $\lambda$.
Note that $\mathbf{C}$ is quadratic by \cite[Corollary~2.3.3]{BGS}.
Using again Theorem~\ref{tkfunctors}\eqref{tkfunctors.3} and
Proposition~\ref{primpi} we get, completely analogous to our previous
argument, that the  adjunction morphism $\op{ID}\rightarrow\mathrm{K}\mathrm{K}'$
is an isomorphism. This implies \eqref{koszul.2}.

$\eqref{koszul.2}\Rightarrow\eqref{koszul.1}$. By
Proposition~\ref{primpi} we have
$\mathrm{K}\mathsf{L}_{\mathbf{C}}(\lambda)^\bullet
\cong\mathsf{I}_{\mathbf{C}^!}(\lambda)^\bullet$
for any $\la\in\Ob(\mathbf{C})$. From
Theorem~\ref{tkfunctors}\eqref{tkfunctors.2} and \eqref{tkfunctors.3}
we know that $\mathrm{K}'\mathrm{K}\mathsf{L}_{\mathbf{C}}(\lambda)^\bullet$
is a linear complex of projective $\mathbf{C}$-modules. Since, by assumption,
$\mathrm{K}'\mathrm{K}\mathsf{L}_{\mathbf{C}}(\lambda)^\bullet\cong
\mathsf{L}_{\mathbf{C}}(\lambda)^\bullet$, the module
$\mathsf{L}_{\mathbf{C}}(\lambda)^\bullet$ has a linear projective
resolution, which implies \eqref{koszul.1}.

$\eqref{koszul.1}\Leftrightarrow\eqref{koszul.4}$. From Proposition~\ref{primpi} we get $\mathrm{K}'
\mathsf{I}_{\mathbf{C}^!}(\lambda)
\cong \mathrm{S}_{-1}\mathrm{Q}_{0}
\mathcal{Q}_{\lambda}^{\bullet}$, and the
statement is clear.

$\eqref{koszul.1}\Rightarrow\eqref{koszul.3}$. This follows from
Propositions~\ref{pltcproj} and~\ref{primpi}, since any Koszul category is quadratic (\cite[Corollary 2.3.3]{BGS}).

$\eqref{koszul.3}\Rightarrow\eqref{koszul.1}$. Since we assume $\mathrm{K}\mathsf{P}(\lambda)^\bullet=
\mathsf{L}_{\mathbf{C}^!}(\lambda)^\bullet$ for every
$\lambda\in\Ob(\mathbf{C})$, Proposition~\ref{prcat508}
implies that the minimal injective resolution
of any $\mathsf{L}_{\mathbf{C}^!}(\lambda)^\bullet$ is linear, hence
$(\mathbf{C}^!)$ is Koszul, and therefore so is $(\mathbf{C}^!)^!$. So, it is
enough to show that $\mathbf{C}$ is quadratic. If it is not, then there is
some $\la\in\Ob(\mathbf{C})$ such that the characters of $\mathsf{P}_{\mathbf{C}}(\la)$ and
$\mathsf{P}_{(\mathbf{C}^!)^!}(\la)$ do not agree since
$\mathbf{C}$ is then a proper quotient of $(\mathbf{C}^!)^!$ (it does
have more relations). From
Theorem~\ref{tkfunctors}\eqref{tkfunctors.3} we get that if
$\mathrm{K}\mathsf{P}_{\mathbf{C}}(\la)^\bullet\cong
\mathsf{L}_{\mathbf{C}^!}(\la)$ then
$\mathrm{K}\mathsf{P}_{{(\mathbf{C}^!)}^!}(\la)^\bullet\not
\cong\mathsf{L}_{{(\mathbf{C}^!)^!}^!}(\la)\cong
\mathsf{L}_{\mathbf{C}^!}(\la)$.
This however contradicts Theorem~\ref{tkfunctors}\eqref{tkfunctors.3}.

$\eqref{koszul.1}\Rightarrow\eqref{koszul.5}$. This follows directly from
Proposition~\ref{primpi}.

$\eqref{koszul.5}\Rightarrow\eqref{koszul.1}$. We only have to show that, if
$\mathcal{Q}_{\lambda}^{\bullet}$ is not linear, then
the isomorphism class of the minimal projective resolution
$\mathcal{Q}_{\lambda}^{\bullet}$ of $\mathsf{L}(\lambda)$ in
$\mathcal{D}^{\downarrow}(\mathbf{C}\defis\fgmod)$ does
not intersect the image of $\mathrm{Tot}$. Assuming the contrary we have $\mathrm{K}'\mathcal{X}^{\bullet}\cong
\mathcal{Q}_{\lambda}^{\bullet}$ for some
$\mathcal{X}^{\bullet}\in \mathcal{D}^{\uparrow}(\mathbf{C}^!\defis\fgmod)$
by Proposition~\ref{funnyd2}.
Then $\mathrm{K}\mathrm{K}'\mathcal{X}^{\bullet}\cong
\mathsf{I}_{\mathbf{C}^!}(\lambda)$ by
Theorem~\ref{tkfunctors}\eqref{tkfunctors.3} and
\begin{displaymath}
\mathcal{Y}^{\bullet}:=\mathrm{K}'\mathrm{K}\mathrm{K}'
\mathcal{X}^{\bullet}\cong
\mathrm{S}_1\mathrm{Q}_0\mathcal{Q}_{\lambda}^{\bullet}
\end{displaymath}
by Proposition~\ref{pltcproj}\eqref{concrete}.
The adjunction of $\mathrm{K}'$ and $\mathrm{K}$
(Theorem~\ref{tkfunctors}\eqref{tkfunctors.1})
implies the existence of maps
\begin{equation}\label{eq9087}
\mathcal{Q}_{\lambda}^{\bullet}
\to\mathcal{Y}^{\bullet}\to
\mathcal{Q}_{\lambda}^{\bullet},
\end{equation}
whose composition is the identity map. Since both,
$\mathcal{Q}_{\lambda}^{\bullet}$ and
$\mathcal{Y}^{\bullet}$, are complexes
of projective modules bounded from the right, the maps in \eqref{eq9087} 
can be realized already in the homotopy category (see e.g. 
\cite[Chapter III(2), Lemma~2.1]{Ha}). We obtain that
$\mathcal{Q}_{\lambda}^{\bullet}$ is
a direct summand of $\mathcal{Y}^{\bullet}$, which is
impossible since $\mathcal{Y}^{\bullet}$ is linear and
$\mathcal{Q}_{\lambda}^{\bullet}$ is not. The theorem follows.
\end{proof}

\begin{remark}\label{tstructure}
{\rm
Analogously to \cite[2.13]{BGS}, linear complexes can be interpreted
as objects of the core of a non-standard $t$-structure on the category
$\mathcal{D}^{\downarrow}(\mathbf{C}\defis\fgmod)$ (and
other derived categories we consider). In the case of Koszul
categories, the Koszul duality functors transform the standard
$t$-structure on $\mathcal{D}^{\downarrow}(\mathbf{C}\defis\fgmod)$
into the non-standard $t$-structure on
$\mathcal{D}^{\uparrow}(\mathbf{C}^!\defis\fgmod)$ and vice versa.
}
\end{remark}

We would like to emphasize the following direct consequence:

\begin{corollary}\label{ckoszul}
All projective resolutions of simple $\mathbf{C}$-modules are
linear if and only if they all belong to the image of the functor
$\mathrm{Tot}$ from Proposition~\ref{funnyd2}.
\end{corollary}

\section{Koszul dual functors for the category $\mathcal{O}$}\label{s4}

In this section we apply the results from Section~\ref{section3}
to Koszul algebras associated with the blocks of the classical
Bernstein-Gelfand-Gelfand category $\mathcal{O}$ (see \cite{BGG}, \cite{BGS}). We give an alternative
proof of the result of Ryom-Hansen (\cite{Steen}) on the Koszul duality
of translation and Zuckerman functors on $\mathcal{O}$, and prove
the Koszul duality of twisting/completion and shuffling/coshuffling
functors. In the next section we will describe several applications, in
particular, we will give an alternative proof of the categorification
results of Sussan (\cite{Josh}) by applying Koszul duality to the
corresponding categorification result from \cite{StDuke}.

\subsection{Category $\mathcal{O}$: notation and preliminaries}\label{s4.1}

For any (complex) Lie algebra $\mg$ we denote by $\cU(\mg)$ its universal
enveloping algebra. Let $\mg$ be a complex semisimple Lie algebra with a
fixed Cartan subalgebra $\mh$ inside a Borel subalgebra $\mb$.
Let $\cO=\cO(\mg)$ be the corresponding category
$\cO$ from \cite{BGG} given by all finitely generated $\cU(\mg)$-modules,
which are $\mh$-diagonalizable and locally $\cU(\mb)$-finite. The morphisms
are ordinary $\cU(\mg)$-homomorphisms. The Weyl group $W$ acts naturally on
$\mh^\ast$, via $(x,\la)\mapsto x(\la)$ for any $x\in W$ and $\la\in\mh^\ast$. There is also the so-called ``dot-action'' $x\cdot\la=x(\la+\rho)-\rho$ with
the fixed point $-\rho$, where $\rho$ is the half-sum of positive roots.
It is well-known that the category $\cO$ has enough projectives and
injectives. For $\mu\in\mh^\ast$ let $L(\mu)$ denote the simple module
with the highest weight $\mu$, $P(\mu)$ denote the indecomposable
projective cover and $I(\mu)$ denote the indecomposable injective hull
of $L(\mu)$ in $\cO$.

The action of the center of $\cU(\mg)$ decomposes the category into blocks,
i.e. $\cO=\oplus\cO_\chi$, where (due to the Harish-Chandra isomorphism)
the blocks are indexed by the $W$-orbits under the dot-action. We also
write $\cO_\chi=\cO_\la$, if $\la\in\chi$ is maximal (in the usual ordering
on weights). In particular, $\cO_0$ denotes the principal block containing
the trivial representation, and $P(\mu)$ (resp. $L(\mu)$ or $I(\mu)$) is
an object of $\cO_\la$ if and only if $\mu\in W\cdot\la$. The module
$P_\chi=\oplus_{\mu\in\chi} P(\mu)$ is a minimal projective generator for
$\cO_\chi$, hence $\cO(P_\chi,\bullet)$ defines an equivalence of
categories between $\cO_\chi$ and the category of finitely generated
(which means finite dimensional) right $\END_\cO(P_\chi)$-modules (\cite[Section 2]{Bass}). From \cite{BGS} it is known that
$A(\chi)=\END_\cO(P_\chi)$ can be equipped with a positive $\mZ$-grading such
that the corresponding graded algebra $\mathsf{A}(\chi)$
becomes a Koszul algebra. Since we always worked with left
modules so far we use the duality on $\cO$ to identify $A(\chi)\cong
A(\chi)^{op}$. We denote by $\mathbf{A}(\chi)$ the corresponding
positively graded $\mathbb{C}$-category (recall that the objects of
$\mathbf{A}(\chi)$ can be considered as a minimal system of representatives of
the isomorphism classes of indecomposable
projective modules in $A(\chi)\defis\fgmod$ with head concentrated in degree
zero and the morphisms are the morphisms of graded modules, see
Subsection~\ref{s2.1}). We will identify the objects of the
category $\mathbf{A}(\chi)$ either with  isomorphism classes of
indecomposable projective objects in $A(\chi)\defis\fmod$, or
the  isomorphism classes
of simple modules in $A(\chi)\defis\fmod$ or even just
with the  corresponding highest weights depending on what is the
most convenient way in any particular situation.

Then the category $\mathbf{A}(\chi)\defis\fgmod$ of all finite dimensional
{\it graded} $\mathbf{A}(\chi)$-modules is a ``graded version'' of
$\cO_\chi$. We will also write $\mathbf{A}(\la)$
(resp. $A(\la)$) instead of  $\mathbf{A}(\chi)$
(resp. $A(\chi)$) if $\la\in\chi$ is maximal. In particular, we have $\mathbf{A}(0)\defis\fgmod$, the graded version
of the principal block $\cO_0$.

If $P(\mu)\in\cO_\la$ then we have the corresponding indecomposable
projective $P(\mu)\in A(\la)\defis\fmod$ and
$\mathsf{P}(\mu)=\mathbf{A}(\la)(\mu,_-)\in\mathbf{A}(\la)\defis\fgmod$.
Similarly, $L(\mu)\in\cO_\la$ corresponds to a simple module $L(\mu)\in
A(\la)\defis\fmod$ and to $\mathsf{L}(\mu)\in\mathbf{A}(\la)\defis\fgmod$, the
simple quotient of $\mathsf{P}(\mu)$. Recall that we denoted the injective
hull of $\mathsf{L}(\mu)$ by $\mathsf{I}(\mu)$. The indecomposable projective
modules in $\mathbf{A}(\la)\defis\fgmod$ are exactly the modules of the form
$\mathsf{P}(\mu)\langle j\rangle$ for some $\mu\in W\cdot\la$ and $j\in\mZ$.

For more details  concerning this graded version of category $\cO$ (in the
language of modules over graded algebras) we refer
to  \cite{BGS} and also \cite{Stgrad}.

\subsection{The parabolic categories ${}_\Lambda A(\chi)\defis\fgmod$}

If $\mp\supseteq\mb$ is a parabolic subalgebra of $\mg$, then we denote by
$W_\p\subseteq W$ the corresponding parabolic subgroup, and let $\cO^\p$
denote the full subcategory of $\cO$ given by all locally $\cU(\mp)$-finite
objects.  For a $W$-orbit $\chi$ with maximal weight $\la$
let $\cO^\p_\chi=\cO_\la^\p$ be the full subcategory having as objects all
the objects from $\cO_\chi$, which are locally $\cU(\mp)$-finite. We will
call these categories also ``blocks'', although they are not indecomposable
and of course not even nontrivial in general. The {\it Zuckerman functors}
\begin{eqnarray*}
  \mathrm{Z}^\p_\la:\cO_\la&\rightarrow&\cO_\la^\p
\end{eqnarray*}
are defined as the functors of  taking the maximal $\p$-locally
finite quotient. These functors are right exact.
Let $\mathrm{Z}^\p=\mathrm{Z}_0^\p$ and let $\mathrm{i}^\p$
denote its right adjoint,
i.e. $\mathrm{i}^\p$ is nothing else than the inclusion functor
$\cO_0^\p\rightarrow \cO_0$.
Note that $\mathrm{Z}^\p P(x\cdot0)\not=0$ if and only if $x\in W^\p$,
the set of shortest coset representatives of $W_\p\backslash W$.
The module $\mathrm{Z}^\p_\la P_{\chi}$ is a minimal projective
generator for $\cO_\la^\p=\cO_\chi^\p$. Let $A(\chi)^\p$ denote its endomorphism
ring, which is the quotient of $A(\chi)$ modulo the homogeneous ideal
generated by all idempotents corresponding to the simple modules which
are in $\cO_\chi$ but not in $\cO^\mp_\chi$. In particular,
$A(\chi)^\p$ inherits a (positive) grading from $A(\chi)$. We will
consider the positively graded category corresponding to $A(\chi)^\p$
via the correspondence~\eqref{introequiv} and denote it by
$\mathbf{A}(\chi)^\p$. Using the language from Section~\ref{s3.5} we have

\begin{lemma}
\label{parab}
There is a canonical isomorphism of categories, $\mathbf{A}(\chi)^\p
\cong {}_\Lambda(\mathbf{A}(\chi))$, where
$\Lambda=\Lambda(\chi,\p)$ is the set of idempotents corresponding
to simple modules in $\cO_\chi$, which are contained in $\cO_\chi^\p$.
\end{lemma}

\begin{proof}
This follows directly from the definitions.
\end{proof}

\subsection{The category $\mathbf{A}(0)_{\Lambda(\p)'}\defis\fgmod$}

Let $W_\p$ be a parabolic subgroup of $W$. Let
$W(\p)=\{x\in W\mid x^{-1}w_0\in W^\p\}$, where $w_0$ is the
longest element in $W$. Let $\Lambda'=\Lambda(\p)'$
denote the set of weights of the form $x\cdot0$, $x\in W(\p)$.
We consider $\Lambda(\p)'$ as a subset of $\Ob(\mathbf{A}(0))$.
Then the following holds

\begin{lemma}
\label{AL}
A complete system of indecomposable projective  objects in the category
$\mathbf{A}(0)_{\Lambda'}\defis\fgmod$ is given by restricting the
modules $\mathsf{P}(\la)\langle j\rangle\in \mathbf{A}(0)\defis\fgmod$,
where  $\la\in\Lambda'$, $j\in\mZ$, to objects in $\mathbf{A}(0)_{\Lambda'}\defis\fgmod$.
\end{lemma}

\begin{proof}
This follows directly from the definition of $\mathbf{A}(0)_{\Lambda'}$.
\end{proof}

\subsection{Koszul duality of translation and Zuckerman functors}\label{s4.2}

Let $\p\supseteq\mb$ a parabolic subalgebra and let
$\la=\la(\p)\in\mh^\ast$ be such
that its stabilizer under the dot-action is $W_\p$ and it is maximal in its
orbit. In \cite{BGS} (and \cite{ErikKoszul}), it is proved that $A(\chi)^\p$
is always a Koszul algebra. More precisely (\cite[Corollary~3.7.3]{BGS}),
there is an isomorphism of graded algebras
$\mathsf{A}(0)^!\rightarrow \mathsf{A}(0)$,
which induces an isomorphism of categories $\mathbf{A}(0)^!\rightarrow \mathbf{A}(0)$, such that the object ${x\cdot0}$ is
mapped to the object ${x^{-1}w_0\cdot 0}$, where $w_0$ is the longest element
in $W$. More generally, (\cite[Proposition~3.1]{ErikKoszul}), there
is an isomorphism of categories $(\mathbf{A}(\la(\p)))^!\cong\mathbf{A}(0)^\p$ mapping  the object
${x\cdot\la(\p)}$ to the object ${x^{-1}w_0\cdot 0}$. For any $\p$ we fix
such an isomorphism and the induced isomorphism of categories
$\sigma^\p:\mathcal{D}^{\uparrow}
\big(\mathbf{A}{(\la(\p))}^!\defis\fgmod\big)\cong
\mathcal{D}^{\uparrow}({\mathbf{A}(0)^\p}\defis\fgmod)$.
Set $\sigma=\sigma^\mb$. We have the following Koszul duality functors:
\begin{eqnarray*}
\xymatrix{
\mathcal{D}^{\downarrow}(\mathbf{A}(0)\defis\fgmod)
\ar@/^/[rrr]^{\mathrm{K}_{\mathbf{A}(0)}}
&&&
\ar@/^/[lll]^{\mathrm{K}'_{\mathbf{A}(0)}}
\mathcal{D}^{\uparrow}({\mathbf{A}(0)}^!\defis\fgmod)
\stackrel{\sigma}{\cong}\mathcal{D}^{\uparrow}({\mathbf{A}(0)}\defis\fgmod),
}
\end{eqnarray*}
such that
$\mathrm{K}\mathsf{L}(x\cdot0)\cong\mathsf{I}(x^{-1}w_0\cdot0)$. More
generally,
\small
\begin{eqnarray*}
\xymatrix{
\mathcal{D}^{\downarrow}\big({\mathbf{A}({\la(\p)})}\defis\fgmod\big)
\ar@/^/[rrr]^{\mathrm{K}_{\mathbf{A}(\la(\p))}}
&&&
\ar@/^/[lll]^{\mathrm{K'}_{\mathbf{A}(\la(\p))}}
\mathcal{D}^{\uparrow}
\big(\mathbf{A}(\la(\p))^!\defis\fgmod\big)
\stackrel{\sigma^\p}{\cong}
\mathcal{D}^{\uparrow}(\mathbf{A}(0)^\p\defis\fgmod)
},
\end{eqnarray*}
\normalsize
such that $\mathrm{K}_{\mathbf{A}(\la(\p))}\mathsf{L}(x\cdot\la(\p))
\cong I(x^{-1}w_0\cdot 0)$. 

As the algebra $\mathsf{A}$ is finite-dimensional, we have that
the bounded derived category $\mathcal{D}^b(\mathbf{A}(0)\defis\fdgmod)$ is by definition contained in $\mathcal{D}^{\downarrow}(\mathbf{A}(0)\defis\fgmod)$ as well as in $\mathcal{D}^{\uparrow}(\mathbf{A}(0)\defis\fgmod)$. Hence it makes sense to restrict the functors to this subcategory. Since the Koszul functor sends simple modules to injective module (Theorem~\ref{tkfunctors}) and the involved algebra has finite global dimension, we obtain functors as follows (see \cite[Theorem 2.12.6]{BGS} for details): 

\small
\begin{eqnarray}
\label{restrict}
\xymatrix{
\mathcal{D}^{b}(\mathbf{A}(0)\defis\fdgmod)
\ar@/^/[rrr]^{\mathrm{K}_{\mathbf{A}(0)}}
&&&
\ar@/^/[lll]^{\mathrm{K}'_{\mathbf{A}(0)}}
\mathcal{D}^{b}({\mathbf{A}(0)}^!\defis\fdgmod)
\stackrel{\sigma}{\cong}\mathcal{D}^{b}({\mathbf{A}(0)}\defis\fdgmod).
}
\end{eqnarray}
\normalsize

With the notation from Lemma~\ref{parab} the Zuckerman
functors induce functors
\begin{eqnarray*}
\mathrm{Z}^\p:\mathbf{A}(0)\defis\fdgmod&\rightarrow&\mathbf{A}(0)^\p
\defis\fdgmod\cong{}_\Lambda \mathbf{A}(0)\defis\fdgmod\\
\mathrm{i}^\p:{}_\Lambda \mathbf{A}(0)^\p\defis\fdgmod\cong \mathbf{A}(0)^\p\defis\fdgmod&\rightarrow&\mathbf{A}(0)\defis\fdgmod.
\end{eqnarray*}

On the other hand, for any block $\cO_\la$, where $\la$ is integral, we have
the translation functors
\begin{eqnarray*}
  \theta_0^\la:&& \cO_0\rightarrow \cO_\la\\
  \theta_\la^0:&& \cO_\la\rightarrow\cO_0,
\end{eqnarray*}
given by translation onto and out of the wall (for details see for example
\cite{Ja1}, \cite{GJ}). They induce functors
\begin{eqnarray*}
  \theta_0^\la:&&\mathbf{A}(0)\defis\fdmod\rightarrow
  \mathbf{A}(\la)\defis\fdmod\\
  \theta_\la^0:&&\mathbf{A}(\la)\defis\fdmod\rightarrow
  \mathbf{A}(0)\defis\fdmod.
\end{eqnarray*}
In \cite{Stgrad} it is proved that the latter have graded lifts
\begin{eqnarray*}
  \tilde\theta_0^\la:&&\mathbf{A}(0)\defis\fdgmod\rightarrow
  \mathbf{A}(\la)\defis\fdgmod\\
  \tilde\theta_\la^0:&&\mathbf{A}(\la)\defis\fdgmod\rightarrow
  \mathbf{A}(0)\defis\fdgmod,
\end{eqnarray*}
which give rise to the original functors if we forget the grading. We are mostly interested in the case when $\la$ is integral and the stabilizer
$W_\la$ of $\la$ is generated by a simple reflection $s$, that is $\la$ is
``lying on exactly one wall''. To avoid even more notation we restrict from
now on to this case. We fix a {\it
standard lift} $\tilde\theta_0^\la$ such that $\theta_0^\la$ maps
$\mathsf{P}(0)$ to $\mathsf{P}(\la)$. We fix a {\it standard lift}
$\tilde\theta_\la^0$ of $\theta_\la^0$ such that the adjunction morphism
$\op{ID}\rightarrow\tilde\theta_\la^0\tilde\theta_0^\la$
is homogeneous of degree
$1$. This means
$\tilde\theta_\la^0\mathsf{P}(\la)\cong\mathsf{P}(s\cdot0)$. For more
details we refer to \cite[Section 1, Section 3.2]{Stgrad}. \\

As an application of our general setup we get the following
result, conjectured in~\cite{BGS}, and originally proved in
\cite{Steen}, concerning the restrictions of the Koszul functors as given in \eqref{restrict}:

\begin{theorem}
\label{Steen}
Let $\p\supseteq\mb$ be a parabolic subalgebra of $\mg$ such that
$W_{\p}=\{1,s\}$ for some simple reflection $s$. The following
diagrams commute up to isomorphisms of functors:
\begin{eqnarray}
\label{thetaZ}
\xymatrix{
\mathcal{D}^{b}(\mathbf{A}(0)\defis\fdgmod)
\ar[rrrr]^{\sigma\mathrm K_{\mathbf{A}(0)}}
\ar[dd]_{\tilde\theta_0^\la\langle 1\rangle}
&&&&
\mathcal{D}^{b}(\mathbf{A}(0)\defis\fdgmod)
\ar[dd]^{\cL{Z}^\p}
\\\\
\mathcal{D}^{b}(\mathbf{A}(\la(\p))\defis\fdgmod)
\ar[rrrr]^{\sigma^\p\mathrm{K}_{\mathbf{A}(\la(\p))}}
&&&&
\mathcal{D}^{b}(\mathbf{A}(0)^\p\defis\fdgmod)
}
\end{eqnarray}
\begin{eqnarray}
  \label{thetai}
\xymatrix{
\mathcal{D}^{b}(\mathbf{A}(\la(\p))\defis\fdgmod)
\ar[rrrr]^{\sigma^\p\mathrm K_{\mathbf{A}(\la(\p))}}
\ar[dd]_{\tilde\theta^0_\la}
&&&&
\mathcal{D}^{b}(\mathbf{A}(0)^\p\defis\fdgmod)
\ar[dd]^{\mathrm{i}_\p}
\\\\
\mathcal{D}^{b}(\mathbf{A}(0)\defis\fdgmod)
\ar[rrrr]^{\sigma\mathrm{K}_{\mathbf{A}(0)}}
&&&&
\mathcal{D}^{b}(\mathbf{A}(0)\defis\fdgmod)
}
\end{eqnarray}
The category $\mathcal{D}^{b}(\mathbf{A}(0)_{\Lambda(\p)'}\defis\fdgmod)$,
considered as a subcategory of the category
$\mathcal{D}^{b}(\mathbf{A}(0)\defis\fdgmod)$,
is exactly the image of the translation functor $\tilde\theta^0_\la$.
\end{theorem}

\begin{corollary}\label{cor345}
\begin{enumerate}[(i)]
\item\label{cor345.1}
The functor $\mathrm{i}_\p\mathcal{L} Z^\p\langle 1\rangle[-1]$ and
the translation functor
$\tilde\theta_s=\tilde\theta_0^\la\tilde\theta_0^\la$ through
the $s$-wall  are Koszul dual to each other.
\item\label{cor345.2}
The functor $\mathrm{i}_\p\cL Z^\p\langle 1\rangle[-1]$ is both
left  and right adjoint to itself.
\end{enumerate}
\end{corollary}

\begin{proof}
\eqref{cor345.1} follows directly from
Theorem~\ref{Steen} and
Theorem~\ref{tkfunctors}~\eqref{tkfunctors.2}.
The statement \eqref{cor345.2} is then clear, since
$\tilde\theta_s$ is self-adjoint (\cite[Corollary 8.5]{Stgrad}).
\end{proof}

\begin{remark}
{\rm
Our result differs from the one in \cite{Steen} by a shift in the
grading. This is because in \cite{Steen} the graded lift of the
$\theta_0^\la$ is chosen such that it maps a simple module concentrated in
degree $k$ to zero or to a simple module concentrated in degree $k$. We
chose the lift such that it maps a simple module concentrated in
degree $k$ to zero or to a simple module concentrated in degree $k-1$ (see
\cite[Theorem 8.1]{Stgrad}).
}
\end{remark}

To prove Theorem~\ref{Steen} we use the following auxiliary statement:

\begin{lemma}
\label{lemmakos}
Let $\p\supseteq\mb$ be a fixed parabolic subalgebra of $\mg$.
  For any $x\in W(\p)$ there is an isomorphism
  $\beta_x:\mathsf{P}(x\cdot0)\cong\tilde\theta_\la^0
  \mathsf{P}(x\cdot\la(\p))\in
  \mathbf{A}(0)\defis\fdgmod$.
\end{lemma}

\begin{proof}
 Using \cite[Theorem 8.4]{Stgrad} and \cite[4.12 (3)]{Ja2} we get for any $j\in\mZ$ that
\begin{eqnarray*}
&&{\mathbf{A}(0)\defis\fdgmod}(\tilde\theta_\la^0
P(x\cdot\la(\p)),\mathsf{L}(y\cdot0)\langle j\rangle)\\
&\cong&
{\mathbf{A}(0)\defis\fdgmod}(\mathsf{P}(x\cdot\la(\p)),
\tilde\theta_0^\la\mathsf{L}(y\cdot0)\langle j-1\rangle)
\end{eqnarray*}
is only nonzero if $y=x$ and $j=0$, in which case it is isomorphic to
$${\mathbf{A}(0)\defis\fdgmod}(\mathsf{P}(x\cdot\la(\p)),
\mathsf{L}(x\cdot\la(\p)))=\mC.$$ Since the translation
functors map projective objects to projective objects, the
statement follows.
\end{proof}

\begin{proof}[Proof of Theorem~\ref{Steen}.]
By adjointness it is enough to prove the commutativity of the second diagram. We start with some general statements. Let $x, y\in W(\p)$ and
$\mathsf{P}=\mathsf{P}(x\cdot\la(\p))$,
$\mathsf{Q}=\mathsf{P}(y\cdot\la(\p))\in
\mathbf{A}(\la(\p))\defis\fdgmod$. The functor
$\mathrm{T}:=\tilde\theta_\la^0$ is exact and does not
annihilate any submodule of a given projective module
(\cite[4.13 (5) or (3')]{Ja2}), hence it induces a natural inclusion
\begin{eqnarray}
\label{identificationhom}
\mathbf{A}(\la(\p))\defis\fdgmod(\mathsf{P}
\langle -1\rangle,\mathsf{Q})&\hookrightarrow&
{\mathbf{A}(0)\defis\fdgmod}(\tilde\theta_\la^0
\mathsf{P}\langle-1\rangle,\tilde\theta_\la^0
\mathsf{Q})
\end{eqnarray}
of graded vector spaces. We claim that this is even an isomorphism. By
\cite[Theorem 8.4, Proposition 6.7 (2)]{Stgrad} we have
\begin{eqnarray*}
&&  {\mathbf{A}(0)\defis\fdgmod}(\tilde\theta_\la^0
  \mathsf{P}\langle-1\rangle,\tilde\theta_\la^0 \mathsf{Q})\\
&\cong&
    {\mathbf{A}(0)\defis\fdgmod}(\mathsf{P},
    \tilde\theta^\la_0\tilde\theta_\la^0 \mathsf{Q})\\
&\cong& {\mathbf{A}(0)\defis\fdgmod}(\mathsf{P},
\mathsf{Q}\langle 1\rangle\oplus \mathsf{Q}\langle
  -1\rangle)\\
&\cong& {\mathbf{A}(0)\defis\fdgmod}(\mathsf{P}\langle
  -1\rangle,\mathsf{Q})\oplus
  {\mathbf{A}(0)\defis\fdgmod}(\mathsf{P}\langle 1\rangle,\mathsf{Q})\\
&\cong&{\mathbf{A}(0)\defis\fdgmod}(\mathsf{P}\langle-1\rangle,\mathsf{Q}),
\end{eqnarray*}
the latter follows from the positivity of the grading. Hence the spaces in
\eqref{identificationhom} have the same dimension and the map has to be an
isomorphism. Together with Lemma~\ref{AL}, the functor $\mathrm{T}$
induces an isomorphism $$\alpha:
\mathbf{A}(\la(\p))_1\cong (\mathbf{A}(0)_{\Lambda(\p)'})_1.$$
From Lemma~\ref{lemmakos} and Lemma~\ref{AL} we know that
$\mathrm{T}$ induces a functor \begin{eqnarray*}
  \mathrm{T}:&&\mathscr{LC}
(\mathsf{P}_{\mathbf{A}(\la(\p))})\rightarrow
\mathscr{LC}(\mathsf{P}_{\mathbf{A}(0)}).
\end{eqnarray*}
To show that the
second diagram in Theorem~\ref{Steen} commutes, it is enough (by
Theorem~\ref{koszul} and Proposition~\ref{funnyd2}) to show that the
following diagram commutes:
\begin{eqnarray}
\xymatrix{
\mathscr{LC}(\mathsf{P}_{\mathbf{A}(0)})
\ar@{<-}[r]^{\epsilon^{-1}_{\mathbf{A}(0)}}
&\mathbf{A}(0)\defis\fdgmod\ar@{->}[r]^{\sigma}
\ar@{<-}[d]^{\sigma^{-1}\mathrm{i}^\p\sigma^\p}
&\mathbf{A}(0)^!\defis\fdgmod
\ar@{<-}[d]^{\mathrm{i}^\p}
\\
\mathscr{LC}(\mathsf{P}_{\mathbf{A}(\la(\p))})
\ar@{->}[u]^{\mathrm{T}}
&\mathbf{A}(\la(\p))^!\defis\fdgmod
\ar@{->}[l]_{\epsilon^{-1}_{\mathbf{A}(\la(\p))}}
\ar@{->}[r]^{\sigma^\p}
&\mathbf{A}(0)^\p\defis\fdgmod
}\nonumber
\end{eqnarray}

The right hand square commutes by definition. We have the isomorphisms
$\beta_x:\mathsf{P}(x\cdot0)\cong \mathrm{T}\mathsf{P}(x\cdot\la(\p))$ from
Lemma~\ref{lemmakos}. The explicit description of
$\epsilon^{-1}$ in the proof of Theorem~\ref{tqdual} implies therefore that
for any $\mathsf{M}\in\mathbf{A}(\la(\p))^!\defis\fdgmod$, the components of the complexes $\mathrm{T}\epsilon^{-1}_{\mathbf{A}(\la(\p))}(\mathsf{M})$ and
$\epsilon^{-1}_{\mathbf{A}(0)}\sigma^{-1}\mathrm{i}^\p\sigma^\p(\mathsf{M})$ are
isomorphic via the isomorphism $\beta_x$. Moreover, the isomorphism $\alpha$
implies that we even have an isomorphism of complexes. This isomorphism is
natural by the definition of morphisms in the category of linear complexes of
projective modules. Hence the diagram commutes and implies Theorem~\ref{Steen}.
\end{proof}

\subsection{Koszul duality of twisting and shuffling functors}\label{s4.3}

For any simple reflection $s$ let $T_s:\cO_0\rightarrow \cO_0$ be the
twisting functors described for example in \cite{AS}. Let
$\mathrm{T}_s:\mathbf{A}(0)\defis\fdgmod\rightarrow\mathbf{A}(0)\defis\fdgmod$
be the graded version of $T_s$ such that $\mathrm{T}_s\mathsf{P}(0)$ has
head $\mathsf{L}(s\cdot0)$ (\cite[Proposition 5.1]{FKS}) and let
$\mathrm{G}_s$ be its right adjoint. This functor is a graded version
of Joseph's completion functor (\cite{JEnright}, \cite[Theorem~4]{MS}). Let $\mathrm{C}_s:\mathbf{A}(0)\defis\fdgmod\rightarrow\mathbf{A}(0)\defis\fdgmod$
denote the graded version of Irving's shuffling functor, which is given
by taking the cokernel of the adjunction morphism  $\op{ID}\langle
-1\rangle\rightarrow\tilde\theta_\la^0\tilde\theta_0^\la$. Let $\mathrm{D}_s$ be
its right adjoint, which is given by taking the kernel of the adjunction
morphism $\tilde\theta_\la^0\tilde\theta_0^\la\rightarrow\op{ID}\langle 1\rangle$.
In this section we will prove that twisting functors and shuffling
functors are Koszul dual to each other:

\begin{theorem}
\label{Twist}
For any simple reflection $s$, the following diagrams commute up
to isomorphism of functors
\begin{displaymath}
\xymatrix{
\mathcal{D}^{b}(\mathbf{A}(0)\defis\fdgmod)
\ar[r]^{\sigma\mathrm K_{\mathbf{A}(0)}}
\ar[d]_{\cL\mathrm{C}_s}
&
\mathcal{D}^{b}(\mathbf{A}(0)\defis\fdgmod)
\ar[d]^{\cL\mathrm{T}_s}
\\
\mathcal{D}^{b}(\mathbf{A}(0)\defis\fdgmod)
\ar[r]^{\sigma\mathrm{K}_{\mathbf{A}(0)}}
&
\mathcal{D}^{b}(\mathbf{A}(0)\defis\fdgmod)
}
\end{displaymath}
\begin{displaymath}
\xymatrix{
\mathcal{D}^{b}(\mathbf{A}(0)\defis\fdgmod)
\ar[r]^{\sigma\mathrm K_{\mathbf{A}(0)}}
\ar[d]_{\mathcal{R}\mathrm{D}_s}
&\mathcal{D}^{b}(\mathbf{A}(0)\defis\fdgmod)
\ar[d]^{\mathcal{R}\mathrm{G}_s}
\\
\mathcal{D}^{b}(\mathbf{A}(0)\defis\fdgmod)
\ar[r]^{\sigma\mathrm{K}_{\mathbf{A}(0)}}
&
\mathcal{D}^{b}(\mathbf{A}(0)\defis\fdgmod)
}
\end{displaymath}
\end{theorem}

\begin{remark}\label{rem94256}
{\rm
Applying Proposition~\ref{qdualq}, \cite[Lemma~5.2]{MS2} and
\cite[Corollary~6]{KM}, from Theorem~\ref{Twist}
it also follows that  the following diagrams commute up
to isomorphism of functors
\begin{displaymath}
\xymatrix{
\mathcal{D}^{b}(\mathbf{A}(0)\defis\fdgmod)
\ar[r]^{\sigma\mathrm K_{\mathbf{A}(0)}}
\ar[d]_{\cL\mathrm{T}_s}
&
\mathcal{D}^{b}(\mathbf{A}(0)\defis\fdgmod)
\ar[d]^{\cL\mathrm{C}_s}
\\
\mathcal{D}^{b}(\mathbf{A}(0)\defis\fdgmod)
\ar[r]^{\sigma\mathrm{K}_{\mathbf{A}(0)}}
&
\mathcal{D}^{b}(\mathbf{A}(0)\defis\fdgmod)
}
\end{displaymath}
\begin{displaymath}
\xymatrix{
\mathcal{D}^{b}(\mathbf{A}(0)\defis\fdgmod)
\ar[r]^{\sigma\mathrm K_{\mathbf{A}(0)}}
\ar[d]_{\mathcal{R}\mathrm{G}_s}
&\mathcal{D}^{b}(\mathbf{A}(0)\defis\fdgmod)
\ar[d]^{\mathcal{R}\mathrm{D}_s}
\\
\mathcal{D}^{b}(\mathbf{A}(0)\defis\fdgmod)
\ar[r]^{\sigma\mathrm{K}_{\mathbf{A}(0)}}
&
\mathcal{D}^{b}(\mathbf{A}(0)\defis\fdgmod).
}
\end{displaymath}
}
\end{remark}

The rest of the section will be devoted to the proof of Theorem~\ref{Twist}.
We start with the following definition: Let $\mathscr{A}$, $\mathscr{B}$
be categories and assume  $\mathscr{B}$ is abelian. Let
\begin{eqnarray}
\label{eq:func}
0\rightarrow \mathrm{F}_1\rightarrow \mathrm{F}_2\rightarrow
\mathrm{F}_3\rightarrow 0
\end{eqnarray}
be a complex of functors  $\mathrm{F}_i:\mathscr{A}\rightarrow
\mathscr{B}$, ($1\leq i\leq 3$).  The complex~\eqref{eq:func} is
{\it exact} if it gives rise to a short exact sequence
$0\rightarrow \mathrm{F}_1(M)\rightarrow
\mathrm{F}_2(M)\rightarrow \mathrm{F}_3(M)\rightarrow 0$ in
$\mathscr{B}$ for any object  $M\in\mathscr{A}$. Analogously, if
$\mathscr{B}$ is a triangulated category then for
functors $\mathrm{F}_i:\mathscr{A}\rightarrow \mathscr{B}$ we say that
$\mathrm{F}_1\rightarrow \mathrm{F}_2\rightarrow
\mathrm{F}_3\rightarrow \mathrm{F}_1[1]$ is a
{\em distinguished triangle} if it  gives rise to a distinguished
triangle in $\mathscr{B}$ when evaluated at any  object in $\mathscr{A}$.

\begin{lemma}
\label{lemma1}
Let $s$ be a simple reflection and
$\mathscr{P}\subseteq \mathbf{A}(0)\defis\fdgmod$  be the
full additive category given by all projective objects.
Then there is an exact sequence of functors from $\mathscr{P}$
to $\mathbf{A}(0)\defis\fdgmod$ of the form
\begin{eqnarray*}
0\rightarrow \mathrm{T}_s\rightarrow\op{id}\langle 1\rangle\rightarrow
\mathrm{i_\p}\mathrm{Z}^\p\langle 1\rangle\rightarrow 0,
\end{eqnarray*}
where $\p$ is the parabolic subalgebra of $\mathfrak{g}$, associated
with $s$.
\end{lemma}

\begin{proof}
  We have $\mathrm{T}_s\mathsf{P}(w_0\cdot0)\cong\mathsf{P}(w_0\cdot0)$ by
  definition. From the proof of \cite[Theorem 6]{MS} it follows that
  $\HOM(\mathrm{T}_s, \op{ID}\langle 1\rangle)\cong C(\mh)$, where $C(\mh)$
  is the coinvariant
  algebra as in \cite[1.2]{Sperv}. In particular, there is a unique up to scalar
  natural transformation $\op{can}$ of lowest degree. It must be non-trivial on
  $\mathsf{P}(0)$, otherwise it would be trivial anywhere, since
  $\mathrm{T}_s$ commutes with translation functors through
  walls (\cite[Section~3]{AS}).
  Then the cokernel of $\op{can}$ is $\mathrm{i_\p}\mathrm{Z}^\p\langle
  1\rangle$ (\cite[Proposition 5.4]{AS}).
\end{proof}

\begin{corollary}
\label{cor}
There exists a morphism $\phi$ of functors such that
\begin{eqnarray*}
\mathrm{i_\p}\cL\mathrm{Z}^\p\langle
1\rangle[-1]\stackrel{\phi}{\rightarrow}\cL\mathrm{T}_s
\rightarrow\op{id}\langle 1\rangle\rightarrow
\mathrm{i_\p}\mathrm{Z}^\p\langle 1\rangle
\end{eqnarray*}
is a triangle of functors.
\end{corollary}

\begin{proof}
This follows immediately from Lemma~\ref{lemma1} (see for example
\cite[Proposition~1.8.8]{KS}).
\end{proof}

\begin{proof}[Proof of Theorem~\ref{Twist}]
Since $\mathrm{K}=\mathrm{K}_{\mathbf{A}(0)}$ is a functor of
triangulated categories and an equivalence by Theorem~\ref{koszul},
from Corollary~\ref{cor} we get the triangle
\begin{eqnarray}
\label{triangle}
\mathrm{K}^{-1}\,\mathrm{i_\p}\cL\mathrm{Z}^\p\langle 1\rangle[-1]
\mathrm{K} \stackrel{\phi}{\rightarrow}
\mathrm{K}^{-1}\,\cL\mathrm{T}_s\mathrm{K}
\rightarrow\mathrm{K}^{-1}\op{ID}\langle 1\rangle\mathrm{K}\rightarrow
\mathrm{K}^{-1}\mathrm{i_\p}\mathrm{Z}^\p\langle 1\rangle\mathrm{K}.
  \end{eqnarray}
From Theorem~\ref{tkfunctors} we have
$\mathrm{K}\langle 1\rangle\cong\langle-1\rangle[1]\mathrm{K}$.
Together with Corollary~\ref{cor345} the triangle
\eqref{triangle} gives rise to a triangle
\begin{eqnarray}
\tilde\theta_s\stackrel{\phi'}{\rightarrow}
\mathrm{K}^{-1}\cL \mathrm{T}_s\mathrm{K}\rightarrow\op{ID}
\langle -1\rangle [1]\rightarrow
\tilde\theta_s[1]
\end{eqnarray}
and therefore to the triangle
\begin{eqnarray}
\mathrm{K}^{-1}\cL\mathrm{T}_s\mathrm{K}[-1]\rightarrow\op{ID}
\langle -1\rangle \stackrel{\phi''}{\rightarrow}
\tilde\theta_s\rightarrow
\mathrm{K}^{-1}\cL\mathrm{T}_s\mathrm{K}
\end{eqnarray}
Note that the map $\phi''$ is graded (homogeneous of degree zero).
Since the graded vector space
$\HOM(\op{ID}\langle -1\rangle,\tilde\theta_s)$
is one dimensional in degree zero (\cite[Theorem~4.9]{Backelinhom}),
$\phi''$ must be the adjunction morphism up to a scalar.
Hence we have $\mathrm{K}^{-1}\cL\mathrm{T}_s\mathrm{K}\cong\cL\mathrm{C}_s$.
Therefore, the first diagram of Theorem~\ref{Twist} commutes. The
commutativity of the second follows by adjointness.
\end{proof}

\section{Applications}\label{s5}

Finally we would like to indicate applications of our results.

\subsection{A categorical version of the quantized Schur-Weyl duality}

In \cite{FKS} a categorification of finite dimensional quantum
$\mathfrak{sl}_2$-modules was obtained using certain graded versions
of blocks of the category of Harish-Chandra bimodules for
$\mathfrak{sl}_n$ and translation  functors. The quantized
Schur-Weyl duality was categorified using certain singular blocks
of the category $\cO$ together with the action of twisting functors
and translation functors through walls (see \cite[Section 5]{FKS} based
on \cite[Corollary~1]{BFK}). The standard and the dual canonical bases
were realized using graded versions of Verma modules and simple modules.
Now Theorems~\ref{Steen} and \ref{Twist} provide the Koszul dual version
of it: The Schur-Weyl duality can be categorified using the bounded
derived categories of certain parabolic blocks of $\cO$ (as suggested
in \cite[Section 4]{BFK}) together with the action of shuffling and
derived Zuckerman and inclusion functors. From
Theorem~\ref{tkfunctors}\eqref{tkfunctors.3} and
\cite[Theorem 5.3 (e)]{FKS} it follows directly that the standard
and canonical basis can be realized using graded versions of dual
Verma modules and injective modules.

\subsection{A functorial tangle invariant}

J.~Sussan proved in \cite{Josh} that the categorification from
\cite[3.2.3]{BFK} of the Temperley-Lieb algebra using singular
blocks of category $\cO$ together with Zuckerman functors and
inclusion functors can be extended to a functorial tangle invariant
using derived twisting and derived completion functors.
Theorem~\ref{Twist} shows that the functorial invariants of
\cite{Josh} and of \cite{StDuke} are Koszul dual to each other.

\subsection{A "Koszul dual" for Harish-Chandra bimodules}

The categorification of finite dimensional quantum
$\mathfrak{sl}_2$-modules from \cite{BFK} was obtained using certain
graded versions of blocks of the category of Harish-Chandra bimodules
and translation functors. In general, these graded blocks  are not
Koszul, hence it does not make sense to speak about a Koszul dual
version at all. However, we propose the following  alternative to
the "Koszul dual" of  the graded version of the category
${}_\la\mathcal{H}_\mu^1$ of Harish-Chandra bimodules with
generalized central character $\chi_\la$ from the left hand side
the and central character $\chi_\mu$ from the right hand side:
There is the well-known equivalence from
\cite{BG} which identifies ${}_\la\mathcal{H}_\mu^1$ with a
certain subcategory of $\cO_\la$ (see \cite{BG} or \cite[Section 6]{Ja2}).
By \cite[6.17]{Ja2} the graded version of
${}_\la\mathcal{H}_\mu^1$ is equivalent to
$\mathbf{A}(\la)_{\Lambda}\defis\fdgmod$ for some $\Lambda$. Hence
Lemma~\ref{technical} provides the quadratic dual, namely
${}_\Lambda{\mathbf{A}^\p}\defis\fdgmod$, where $\mathbf{A}^\p\defis\fdgmod$
is the Koszul dual of $\mathbf{A}(\la)\defis\fdgmod$. Using the Koszul
duality of translation and Zuckerman functors (Theorem~\ref{Steen} and
\cite{Steen}), we get directly from Theorem~\ref{tqdfunctors} a quadratic dual version of the results in \cite{FKS}.

\subsection{A Koszul duality for Kac-Moody Lie algebras}

In our opinion, one advantage of our setup using graded categories in
comparison with the setup in \cite{BGS} is the fact that the categories are
allowed to have infinitely many objects. Instead of considering the principal
block of the category $\cO$ for a semisimple Lie algebra we could consider
the category $\cO$ for a symmetrizable complex Kac-Moody algebra and in
there any regular block outside the critical hyperplanes. Translation
functors through walls are defined in \cite{CT}. The generalization
of Soergel's structure theorem (\cite{Sperv}) holds (see
\cite{Fiebigcomb} for the deformed case, and \cite{St4TQFT} for
the non-deformed case). In analogy with \cite{Sperv} and \cite{Stgrad}
the morphism spaces between indecomposable projective objects
(if they exist) can be equipped with a positive grading giving rise to a
positively graded category $\mathbf{C}$ as defined in Section~\ref{s2}. In
case projective objects do not exist, they can be replaced by tilting objects
(see \cite{Fiebigcomb}). Theorem~\ref{tkfunctors} and Theorem~\ref{koszul}
then provide an analogue of the Koszul duality for regular block outside
the critical hyperplanes for the category $\cO$ of a symmetrizable
complex Kac-Moody algebra.

\appendix
\section{Appendix: An abstract generalized Koszul complex}\label{s6}

\subsection*{Some linear algebra}\label{s6.1}

Let $V$ be a finite-dimensional $\Bbbk$-vector space and
$\mathbf{v}=\{v_i:i=1,\dots,n\}$ be a basis in $V$. We denote by
$\{v^i:i=1,\dots,n\}$ the dual basis in $V^*:=\mathbf{d}(V)$.
Then we have a canonical isomorphism,
\begin{equation}\label{apex}
\begin{array}{rrcl}
\varphi:&
V\otimes_{\Bbbk}V^*&\overset{\sim}{\rightarrow}&
\Hom_{\Bbbk}(V,V)\\
&v\otimes w^*&\mapsto& w^*({}_-)v.
\end{array}
\end{equation}
We have $I_V:=\varphi^{-1}(\mathrm{Id}_V)=\sum_{i=1}^k v_i\otimes v^i$,
in particular, the element $I_V$ does not depend on the choice of the
basis $\mathbf{v}$. Let $H\subset V$ be a subspace and
$H^{\perp}=\{f\in V^*\,|\,f(h)=0\text{ for all }h\in H\}$ be the
corresponding orthogonal complement in $V^*$. Then we have
\begin{displaymath}
\begin{array}{lcl}
\varphi(H\otimes_{\Bbbk}V^*) & = &
\{f:V\to V\,|\, \mathrm{Im}(f)\subset H\}\\
\varphi(V\otimes_{\Bbbk}H^{\perp}) & = &
\{f:V\to V\,|\, H\subset\mathrm{Ker}(f)\}.\\
\end{array}
\end{displaymath}

\begin{lemma}\label{la.1}
Let $H\subset V$. Then $I_V\in H\otimes_{\Bbbk}V^*+
V\otimes_{\Bbbk}H^{\perp}$.
\end{lemma}

\begin{proof}
Let $\mathfrak{p}:V\to H$ be any projector on $H$. Then
$\mathrm{Id}_V=\mathfrak{p}+(\mathrm{Id}_V-\mathfrak{p})$,
$\mathrm{Im}(\mathfrak{p})=
\mathrm{Ker}(\mathrm{Id}_V-\mathfrak{p})=H$.
Hence $\varphi^{-1}(\mathfrak{p})\in H\otimes_{\Bbbk}V^*$
and $\varphi^{-1}(\mathrm{Id}_V-\mathfrak{p})\in
V\otimes_{\Bbbk}H^{\perp}$.
\end{proof}

\subsection*{A semi-simple analogue}\label{s6.2}

Let $\mathbf{C}_0$ be as in Section~\ref{s3.1}, and let
$V_{\mathbf{C}_0}$ be an arbitrary right $\mathbf{C}_0$-module.
Define
\begin{displaymath}
{}_{\mathbf{C}_0} V^*={\mathbf{C}_0}\defis\Mod
(V_{\mathbf{C}_0},(\mathbf{C}_0)_{\mathbf{C}_0})
\end{displaymath}
(note that the authors of \cite[2.7]{BGS} use the notation
${}^*V$ for the same object). The formula \eqref{apex}
defines a canonical isomorphism,
\begin{displaymath}
\varphi:V_{\mathbf{C}_0}\otimes_{\mathbf{C}_0}
{}_{\mathbf{C}_0} V^* \overset{\sim}{\rightarrow}
\Hom_{\mathbf{C}_0}(V_{\mathbf{C}_0},V_{\mathbf{C}_0}).
\end{displaymath}
Let $H_{\mathbf{C}_0}\subset V_{\mathbf{C}_0}$ be a (right) submodule.
Then ${}_{\mathbf{C}_0} H^{\perp}$ is  a (left) submodule of
${}_{\mathbf{C}_0} V^*$, and, analogously to
Lemma~\ref{la.1}, we obtain
\begin{equation}\label{ape1}
I_V:=\varphi^{-1}(\mathrm{Id}_V)\in
H_{\mathbf{C}_0}\otimes_{\mathbf{C}_0}{}_{\mathbf{C}_0} V^*+
V_{\mathbf{C}_0}\otimes_{\mathbf{C}_0}{}_{\mathbf{C}_0} H^{\perp}.
\end{equation}

\subsection*{A differential vector space for quadratic duals}\label{s6.3}

Now let $\mathbf{C}$ and $\mathbf{C}^!$ be as in Section~\ref{se3}.
Let $M_{\mathbf{C}}$ be a right $\mathbf{C}$-module and
${}_{\mathbf{C}^!}N$ be a left $\mathbf{C}^!$-module. Let
$\{a_i:i=1,\dots,k\}$ be a basis of $\mathbf{C}_1$ and
$\{a^i:i=1,\dots,k\}$ the corresponding dual basis of $\mathbf{C}^!_1$.

\begin{proposition}\label{prap1}
The linear transformation
\begin{displaymath}
\begin{array}{rrcl}
\delta: & M_{\mathbf{C}}\otimes_{\mathbf{C}_0}
{}_{\mathbf{C}^!}N& \to & M_{\mathbf{C}}\otimes_{\mathbf{C}_0}
{}_{\mathbf{C}^!}N,\\
&m\otimes n&\mapsto& \sum_{i=1}^k m a_i\otimes a^i n
\end{array}
\end{displaymath}
satisfies $\delta^2=0$. Moreover, if both, $M_{\mathbf{C}}$ and
${}_{\mathbf{C}^!}N$, are graded modules, then
$M_{\mathbf{C}}\otimes_{\mathbf{C}_0}
{}_{\mathbf{C}^!}N$ has a canonical bigrading, and
$\delta$ is a homogeneous map of bidegree $(1,1)$.
\end{proposition}

\begin{proof}
That $\delta$ is a homogeneous map of bidegree $(1,1)$ in the
graded situation is clear from the definition. What we have to
prove is that $\delta^2=0$. Let $\mathbf{m}$ denote the
multiplication in $\mathbf{C}$ (see Subsection~\ref{s3.2}), and
$\mathbf{m}^!$ denote the multiplication in $\mathbf{C}^!$.
We have
\begin{displaymath}
\begin{array}{rcl}
\delta^2(m\otimes n) & =&
\sum_{i=1}^k\sum_{j=1}^k ma_ia_j\otimes a^ja^in\\
& =&
m\left((\mathbf{m}\otimes \mathbf{m}^!)
\left(\sum_{i=1}^k\sum_{j=1}^k (a_i\otimes
a_j)\otimes (a^j\otimes a^i)\right)\right)n\\
& =&
m\left((\mathbf{m}\otimes \mathbf{m}^!)
I_{\mathbf{C}_1\otimes_{\mathbf{C}_0} \mathbf{C}_1}\right)n.
\end{array}
\end{displaymath}
Let now $R\subset \mathbf{C}_1\otimes_{\mathbf{C}_0} \mathbf{C}_1$
be the set of quadratic relations of $\mathbf{C}$. Then
$R^{\perp}$ is the set of defining quadratic relations of
$\mathbf{C}^!$ by definition. Then, by \eqref{ape1} we have
$I_{\mathbf{C}_1\otimes_{\mathbf{C}_0} \mathbf{C}_1}=X+Y$,
where
\begin{displaymath}
X\in R\otimes_{\mathbf{C}_0}
\mathbf{C}_1^*\otimes_{\mathbf{C}_0} \mathbf{C}_1^*,\quad\quad
Y\in \mathbf{C}_1\otimes_{\mathbf{C}_0} \mathbf{C}_1
\otimes_{\mathbf{C}_0} R^{\perp}.
\end{displaymath}
Hence $\mathbf{m}\otimes \mathbf{m}^!(X)=0$ and
$\mathbf{m}\otimes \mathbf{m}^!(Y)=0$ and thus
$\delta^2(m\otimes n)=0$.
This completes the proof.
\end{proof}

Consider the vector space $\mathcal{C}^{\bullet}=
\mathcal{C}^{\bullet}(M,N)$ defined
via $\mathcal{C}^{i}=(M_{\mathbf{C}}\otimes_{\mathbf{C}_0}
{}_{\mathbf{C}^!}N)$ for all $i\in\Z$ (which means that
we just place a copy of $M_{\mathbf{C}}\otimes_{\mathbf{C}_0}
{}_{\mathbf{C}^!}N$ in each position).

\begin{corollary}\label{cap1}
\begin{enumerate}[(i)]
\item\label{cap1.1} The linear transformation
\begin{displaymath}
\begin{array}{rccc}
\delta: & \mathcal{C}^{\bullet}& \to & \mathcal{C}^{\bullet},\\
&\mathcal{C}^{i}\ni(m\otimes n)&\mapsto& (\sum_{i=1}^k m a_i\otimes a^i n)\in \mathcal{C}^{i+1}
\end{array}
\end{displaymath}
satisfies $\delta^2=0$, in particular,
$\mathcal{C}^{\bullet}$ is a complex.
\item\label{cap1.2} If $V$ is a $\mathbf{C}$-bimodule and
$W$ is a $\mathbf{C}^!$-bimodule, then
$\mathcal{C}^{\bullet}(V,W)$ is a complex of
$\mathbf{C}\defis\mathbf{C}^!$-bimodules.
\end{enumerate}
\end{corollary}

Finally, assume that both $\mathsf{M}$ and $\mathsf{N}$ are
graded modules and
define $\mathcal{C}^{i}=\mathsf{M}_{\mathbf{C}}\langle i\rangle
\otimes_{\mathbf{C}_0}{}_{\mathbf{C}^!}\mathsf{N}\langle i\rangle$
for all $i\in\mathbb{Z}$.

\begin{corollary}\label{cap2}
\begin{enumerate}[(i)]
\item\label{cap2.1} The linear transformation
from Corollary~\ref{cap1}\eqref{cap1.1}
defines on $\mathcal{C}^{\bullet}$ the structure of
a complex of graded vector spaces (i.e. the differential is
a homogeneous map of degree $0$).
\item\label{cap2.2} If $\mathsf{V}$ is a graded
$\mathbf{C}$-bimodule and $\mathsf{W}$ is a graded $\mathbf{C}^!$-bimodule,
then $\mathcal{C}^{\bullet}(\mathsf{V},\mathsf{W})$ is a complex of
bigraded $\mathbf{C}\defis\mathbf{C}^!$-bimodules.
\end{enumerate}
\end{corollary}

\subsection*{Generalized Koszul complexes}\label{s6.4}

Several known complexes can be obtained by this technique,
for example:

\begin{itemize}
\item The complex $\mathcal{C}^{\bullet}(\mathbf{C},\mathbf{C}^!)$,
given by Corollary~\ref{cap2}\eqref{cap2.2}, is isomorphic to
$\mathbb{P}^{\bullet}$ from Section~\ref{s3.2new} by construction.
\item The complex $\mathcal{C}^{\bullet}(\mathbf{C},(\mathbf{C}^!)^*)$,
given by Corollary~\ref{cap2}\eqref{cap2.2} contains the classical
Koszul complex (as in \cite[2.8]{BGS}) as a subcomplex
(of $\mathbf{C}$-modules). In particular,
$\mathcal{C}^{\bullet}(\mathbf{C},(\mathbf{C}^!)^*)$ can be considered
as a natural bimodule extension of the Koszul complex.
\end{itemize}

Because of the last example it is natural to call the complexes,
given by  Corollary~\ref{cap2}\eqref{cap2.2} {\em generalized Koszul
complexes}.

\noindent
Volodymyr Mazorchuk, Department of Mathematics, Uppsala University,
Box 480, 751 06, Uppsala, SWEDEN,\\
e-mail: {\tt mazor\symbol{64}math.uu.se},
web: {``http://www.math.uu.se/$\tilde{\hspace{1mm}}$mazor/''}.
\vspace{0.2cm}

\noindent
Serge Ovsienko, Department of Mathematics, Kyiv University,
64, Vo\-lo\-dy\-myr\-ska st., 01033, Kyiv, Ukraine,\\
e-mail: {\tt ovsienko\symbol{64}zeos.net},
{\tt ovsko\symbol{64}voliacable.net}.
\vspace{0.2cm}

\noindent
Catharina Stroppel, Department of Mathematics,
University of Glasgow, University Gardens,
Glasgow G12 8QW, UK,\\
e-mail: {\tt cs\symbol{64}maths.gla.ac.uk},
web: {``http://www.maths.gla.ac.uk/$\tilde{\hspace{1mm}}$cs/}.
\vspace{0.2cm}


\begin{thebibliography}{99999999}
\bibitem[AS]{AS}
H.~H.~Andersen, C.~Stroppel, {\em Twisting functors on $\cO$},
Represent. Theory {\bf 7} (2003), 681--699 (electronic).
\bibitem[Au]{Au}
M.~Auslander, {\em Representation theory of Artin algebras. I, II.}
Comm. Algebra  {\bf 1}  (1974), 177--268; {\bf 1} (1974), 269--310.
\bibitem[AR]{AR}
M.~Auslander, I.Reiten, {\em Stable equivalence of dualizing
$R$-varieties.}  Advances in Math.  {\bf 12}  (1974), 306--366.
\bibitem[Ba1]{ErikKoszul}
E.~Backelin, {\em Koszul duality for parabolic and singular category
$\mathcal{O}$},  Represent. Theory  {\bf 3}  (1999), 139--152
(electronic).
\bibitem[Ba2]{Backelinhom}
E.~Backelin, {\em The Hom-spaces between projective functors},
Represent. Theory {\bf 5} (2001), 267--283 (electronic).
\bibitem[Bass]{Bass}
H.~Bass, {\em Algebraic $K$-theory}, Benjamin,
New York-Amsterdam, 1968.
\bibitem[BGS]{BGS}
A.~Beilinson, V.~Ginzburg, and W.~Soergel, {\em Koszul duality patterns
in representation theory},  J. Amer. Math. Soc.  {\bf 9}  (1996),
no. 2, 473--527.
\bibitem[BFK]{BFK}
J.~Bernstein, I.~Frenkel, M.~Khovanov, {\em A categorification of
the Temperley-Lieb algebra and Schur quotients of $U(\mathfrak{sl}_2)$
via projective and Zuckerman functors}, Selecta Math. (N.S.) {\bf 5}
(1999), no. 2, 199--241.
\bibitem[BG]{BG}
J.~Bernstein, S.~Gelfand, {\em Tensor products of finite-
and infinite-dimensional representations of semisimple Lie
algebras},  Compositio Math.  {\bf 41}  (1980), no. 2, 245--285.
\bibitem[BGG1]{BGGKoszul}
I.~Bern{\v{s}}te{\u\i}n, and I.~Gelfand, and
S.~Gelfand,
{\em Algebraic vector bundles on {${\bf P}\sp{n}$} and problems of linear
  algebra}, Funktsional. Anal. i Prilozhen., {\bf 12}, (1978), no. 3, 66--67.
\bibitem[BGG2]{BGG}
I.~Bern{\v{s}}te{\u\i}n, and I.~Gelfand, and
S.~Gelfand, {\em A certain category of {${\mathfrak g}$}-modules},
Funkcional. Anal. i Prilo\v zen. {\bf 10} (1976), no. 2, 1--8.
\bibitem[BoGa]{BoGa}
K.~Bongartz, P.~Gabriel, {\em Covering spaces in representation theory.}
Invent. Math. {\bf 65} (1981/82), no. 3, 331--378.
\bibitem[Br]{Bredon}
G. Bredon,  {\em Equivariant cohomology theories}, {Lecture Notes
in Math}, {\bf 34}, 1967.
\bibitem[CM]{CM}
C.~Cibils, E.~Marcos, {\em Skew category, Galois covering and smash
pro\-duct of a $k$-category}, Proc. AMS, {\bf 134} (2006), no.1, 39--50.
\bibitem[De]{Deligne}
P.~Deligne, {\em Cohomologie a support propre et construction du
foncteur $f^!$}, Lecture Notes in Mathematics {\bf 20}, 1966, p. 404--423.
\bibitem[tDi]{Dieck}
T. tom Dieck,  {\em \"Uber projektive Moduln und Endlichkeitshindernisse
bei Transformationsgruppen}, Manuscripta Mathematica {\bf 34}, 135--155,
(1981).
\bibitem[Fi1]{CT}
P.~Fiebig, {\em  Centers and translation functors for the category
$\cO$ over Kac-Moody algebras}, Math. Z. {\bf 243} (2003), no. 4,
689--717.
\bibitem[Fi2]{Fiebigcomb}
P.~Fiebig, {\em The combinatorics of category $\cO$ for
symmetrizable Kac-Moody algebras}, Transf. Groups. 11 (2006), no. 1, 29-49.
\bibitem[Fl]{Fl}
G.~Fl{\o}ystad, {\em Koszul duality and equivalences of categories},
math.RA/0012264.
\bibitem[FKS]{FKS}
I.~Frenkel, M.~Khovanov, C.~Stroppel,
{\em A categorification of finite-dimensional irreducible
representations  of quantum $sl(2)$ and their tensor products},
math.QA/0511467, to appear in Slecta Mathematica
\bibitem[GJ]{GJ}
O.~Gabber, A.~Joseph, {\em Towards the Kazhdan-Lusztig conjecture},
Ann. Sci. {\'E}cole Norm. Sup. (4) {\bf 14} (1981), no. 3, 261--302.
\bibitem[Ga]{Ga}
P.~Gabriel, {\em Des cat{\'e}gories ab{\'e}liennes.} Bull. Soc. Math.
France  {\bf 90}  (1962), 323--448.
\bibitem[GM]{GM}
S.~Gelfand, Y.~Manin, {\em Methods of homological algebra.}
Second edition, Springer Monographs in Mathematics.
Springer-Verlag, Berlin, 2003.
\bibitem[GK]{GK}
V.~Ginzburg, M.~Kapranov, {\em Koszul duality for operads.}  Duke
Math. J.  {\bf 76}  (1994),  no. 1, 203--272.
\bibitem[GMRSZ]{GMVRSZ}
E.~Green, R.~Mart{\'{\i}}nez-Villa, I.~Reiten, {\O}.~Solberg, D.~Zacharia,
{\em On modules with linear presentations.}  J. Algebra  {\bf 205}
(1998),  no. 2, 578--604.
\bibitem[GRS]{GRS}
E.~Green, I.~Reiten, {\O}.~Solberg, {\em Dualities on generalized
Koszul algebras.} Mem. Amer. Math. Soc. {\bf 159} (2002), no. 754.
\bibitem[Ha]{Ha}
D.~Happel, {\em Triangulated categories in the representation
theory of finite-dimensional algebras}, London Mathematical Society
Lecture Note Series, {\bf 119}, Cambridge University Press, Cambridge,
1988.
\bibitem[HI]{HI}
J.~Herzog, S.~Iyengar, {\em Koszul modules.} J. Pure Appl. Algebra
{\bf 201} (2005), no. 1-3, 154--188.
\bibitem[Ir]{Irvingshuffle}
R.~Irving, {\em Shuffled Verma modules and principal series modules
over complex semisimple Lie algebras}, J. London Math. Soc. (2) {\bf 48}
(1993), no. 2, 263--277.
\bibitem[Ja1]{Ja1}
J.~C.~Jantzen, {\em Moduln mit einem h{\"o}chsten Gewicht},
Lecture Notes in Mathematics, {\bf 750}, Springer, Berlin, 1979.
\bibitem[Ja2]{Ja2}
J.~C.~Jantzen, {\em Einh{\"u}llende Algebren halbeinfacher Lie-Algebren}, Ergebnisse der Mathematik und ihrer Grenzgebiete (3), {\bf 3}.
Springer-Verlag, Berlin, 1983.
\bibitem[Jo]{JEnright}
A.~Joseph, {\em The Enright functor on the Bernstein-Gelfand-Gelfand
category $\cO$}, Invent. Math. {\bf 67} (1982), no. 3, 423--445.
\bibitem[KS]{KS}
M.~Kashiwara, P.~Schapira, {\em Sheaves on manifolds},
Grundlehren der Mathematischen Wissenschaften, {\bf 292}.
Springer-Verlag, Berlin, 1994.
\bibitem[Ke1]{KellerinKZ}
B.~Keller, {\em On the construction of triangle equivalences}, {in:
Derived equivalences for group rings}. {Lecture Notes in Math.},
{\bf 1685}, 155--176, 1998.
\bibitem[Ke2]{Ke2}
B.~Keller, {\em Deriving DG categories},  Ann. Sci. {\'E}cole Norm.
Sup. (4) {\bf 27} (1994), no. 1, 63--102.
\bibitem[Ke3]{Ke3}
B.~Keller, {\em Koszul duality and coderived categories
(after K. Lef{\`e}vre)}, Preprint 2003.
\bibitem[KM]{KM}
O.~Khomenko, V.~Mazorchuk, {\em On Arkhipov's and Enright's
functors}, Math. Z. {\bf 249} (2005), no. 2, 357--386.
\bibitem[McL]{McL}
S.~Mac Lane, {\em Homology}. Reprint of the 1975 edition.
Classics in Mathematics. Springer-Verlag, Berlin, 1995.
\bibitem[MVS]{MVS}
R.~Mart{\'{\i}}nez Villa, M.~Saor{\'{\i}}n, {\em Koszul equivalences
and  dualities.} Pacific J. Math. {\bf 214} (2004), no. 2, 359--378.
\bibitem[MVZ]{MVZ}
R.~Mart{\'{\i}}nez Villa, D.~Zacharia, {\em Approximations with
modules having linear resolutions.}  J. Algebra  {\bf 266}  (2003),
no. 2, 671--697.
\bibitem[Ma]{Ma1}
V.~Mazorchuk, {\em Applications of the category of linear complexes
of tilting modules associated with the category $\mathcal{O}$},
math.RT/0501220, to appear in Alg. Rep. Theory.
\bibitem[MO]{MO}
V.~Mazorchuk and S.~Ovsienko, {\em A pairing in homology and the category
of linear complexes of tilting modules for a quasi-hereditary
algebra}, J. Math. Kyoto Univ. 45 (2005), no. 4, 711-741.
\bibitem[MS1]{MS}
V.~Mazorchuk and C.~Stroppel,
{\em On functors associated to a simple root}, math.RT/0410339,
to appear in J. Algebra.
\bibitem[MS2]{MS2}
V.~Mazorchuk and C.~Stroppel, {\em Translation and shuffling of
projectively presentable modules and a categorification of a parabolic
Hecke module},  Trans. Amer. Math. Soc.  {\bf 357}  (2005),  no. 7,
2939--2973.
\bibitem[Mi]{Mi}
B.~Mitchell, {\em Rings with several objects.}  Advances in Math.
{\bf 8} (1972), 1--161.
\bibitem[Ri]{Ri}
J.~Rickard, {\em Morita theory for derived categories.}
J. London Math. Soc. (2) {\bf 39} (1989), no. 3, 436--456.
\bibitem[RH]{Steen}
S.~Ryom-Hansen, {\em Koszul duality of translation and Zuckerman
functors.} J. Lie Theory {\bf 14} (2004), no. 1, 151--163.
\bibitem[Sc]{Schubert}
H. Schubert, {\em Kategorien. {I}}, {Heidelberger Taschenb\"ucher},
{\bf 66}, Springer, Berlin, (1970).
\bibitem[Sh]{Sh}
U. Shukla, {\em On the projective cover of a module and related results.}  Pacific J. Math.  {\bf 12}  1962 709--717.
\bibitem[So]{Sperv}
W.~Soergel, Kategorie $\cO$, perverse Garben und Moduln {\"u}ber
den Koinvarianten zur Weylgruppe, J. Amer. Math. Soc. {\bf 3}
(1990), no. 2, 421--445.
\bibitem[St1]{Stgrad}
C. Stroppel, {\em Category ${\mathcal{O}}$: gradings and translation
functors},  J. Algebra  {\bf 268}  (2003),  no. 1, 301--326.
\bibitem[St2]{StDuke}
C. Stroppel, {\em Categorification of the Temperley-Lieb category,
tangles, and cobordisms via projective functors},  Duke Math. J.
{\bf 126}  (2005),  no. 3, 547--596.
\bibitem[St3]{St4TQFT}
C. Stroppel, {\em TQFT with corners and  tilting functors in
the Kac-Moody case}, to appear.
\bibitem[Su]{Josh}
J. Sussan, {\em Ph.D. thesis} (in preparation), Yale university
(New Haven).
\end{thebibliography}
\end{document}